%&Plain
%%%%%%%% lezioni di pisa %%%%%
\input amssym.def
\input amssym.tex
%%%%% macros %%%%%
\magnification=1200 

\catcode`\@=11

\hsize=125 mm   \vsize =177mm %normalmente metto 187
\hoffset=4mm    \voffset=10mm
\pretolerance=500 \tolerance=1000 \brokenpenalty=5000

\catcode`\;=\active
\def;{\relax\ifhmode\ifdim\lastskip>\z@
\unskip\fi\kern.2em\fi\string;}

\catcode`\:=\active
\def:{\relax\ifhmode\ifdim\lastskip>\z@\unskip\fi
\penalty\@M\ \fi\string:}

\catcode`\!=\active
\def!{\relax\ifhmode\ifdim\lastskip>\z@
\unskip\fi\kern.2em\fi\string!}

\catcode`\?=\active
\def?{\relax\ifhmode\ifdim\lastskip>\z@
\unskip\fi\kern.2em\fi\string?}

\def\^#1{\if#1i{\accent"5E\i}\else{\accent"5E #1}\fi}
\def\"#1{\if#1i{\accent"7F\i}\else{\accent"7F #1}\fi}

%\frenchspacing

\catcode`\@=12

%\newif\ifpagetitre      \pagetitretrue
%\newtoks\hautpagetitre  \hautpagetitre={\hfil}
%\newtoks\baspagetitre   \baspagetitre={\hfil}

%\newtoks\auteurcourant  \auteurcourant={\hfil}
%\newtoks\titrecourant   \titrecourant={\hfil}

%\newtoks\hautpagegauche \newtoks\hautpagedroite
%\hautpagegauche={\hfil\the\auteurcourant\hfil}
%\hautpagedroite={\hfil\the\titrecourant\hfil}

%\newtoks\baspagegauche  \baspagegauche={\hfil\tenrm\folio\hfil}
%\newtoks\baspagedroite  \baspagedroite={\hfil\tenrm\folio\hfil}

%\headline={\ifpagetitre\the\hautpagetitre
%\else\ifodd\pageno\the\hautpagedroite
%\else\the\hautpagegauche\fi\fi}

%\footline={\ifpagetitre\the\baspagetitre
%\global\pagetitrefalse
%\else\ifodd\pageno\the\baspagedroite
%\else\the\baspagegauche\fi\fi}

%\def\nopagenumbers{\def\folio{\hfil}}
%%%%%%%%
%\hautpagetitre={\hfill\tenbf preliminary version: Not for diffusion\hfill}
%\hautpagetitre={\hfill\tenbf preliminary version: Not for diffusion\hfill}
%\hautpagegauche={\tenbf\folio\hfill\tenrm\the\auteurcourant}
%\hautpagedroite={\tenrm\the\titrecourant\hfill\tenbf\folio}
%\baspagegauche={\hfil} \baspagedroite={\hfil}
%\auteurcourant{Stefano Marmi}
%\titrecourant{Introduction to Small Divisors}
%%%%%%%%%
%\def\mois{\ifcase\month\or January\or February\or March\or April\or
%May\or June\or July\or August\or September\or October\or November\or
%December\fi}
%\def\Date{\rightline{\mois\ /\ \the\day\ /\/ \the\year}}
%\hfuzz=0.3pt
\font\tit=cmb10 scaled \magstep2
\font\titsec=cmb10 scaled \magstep1
\font\titrm=cmr10 scaled \magstep1
\def\dst{\displaystyle}

\def\IM{\mathop{\Im m}\nolimits}
\def\RE{\mathop{\Re e}\nolimits}
\def\H{\Bbb H}
\def\R{\Bbb R}
\def\T{\Bbb T}
\def\Z{\Bbb Z}
\def\Q{\Bbb Q}
\def\C{\Bbb C}
\def\Cbar{\overline{\Bbb C}}
\def\N{{\Bbb N}}
\def\D{{\Bbb D}}
\def\E{{\Bbb E}}
\def\S{{\Bbb S}}
\def\hH+{\hat{\H}^{+}}
\def\hHZ+{\widehat{\H^{+}/\Z}}

\def\be{\beta}

\def\CD#1#2{\hbox{CD}\,(#1,#2)}
\def\cd#1{\hbox{CD}\,(#1)}
\def\abcd{\left(\matrix{ a & b\cr c & d \cr}\right)}

\def\cotg{\hbox{\rm cotg}\,}

\def\proof{\noindent{\it Proof.\ }}
\def\qed{\hfill$\square$\par\smallbreak} 
\def\Proc#1#2\par{\medbreak \noindent {\bf #1\enspace }{\sl #2}%
\par\ifdim \lastskip <\medskipamount \removelastskip \penalty 55\medskip \fi}%
\def\Def#1#2\par{\medbreak \noindent {\bf #1\enspace }{#2}%
\par\ifdim \lastskip <\medskipamount \removelastskip \penalty 55\medskip 
\fi}
\def\qed{\hfill$\square$\par\smallbreak}
%\long\def\Proc#1#2\par{\medbreak \noindent {\bf #1\enspace }{\sl #2}%
%\par\ifdim \lastskip <\medskipamount \removelastskip \penalty 55\medskip 
%\fi}
\def\hfl#1#2{\smash{\mathop{\hbox to 12mm{\rightarrowfill}}
\limits^{\scriptstyle#1}_{\scriptstyle#2}}}
\def\vfl#1#2{\llap{$\scriptstyle #1$}
\left\downarrow\vbox to 6mm{}\right.\rlap{$\scriptstyle #2$}}
\def\diagram#1{\def\normalbaselines{\baselineskip=0pt
\lineskip=10pt\lineskiplimit=1pt} \matrix{#1}}

%%% fin macros %%%%%%%%

%%%% new fonts %%%%%%%%

\baselineskip=14 pt
%%%%% fine macros %%%%%
\vfill\eject
%\nopagenumbers
%%%%% pagina 0: da togliere. Taglia da qui%%%%%
\pageno=0
\vskip 1. truecm
\centerline{UNIVERSIT\`A DI PISA}
\centerline{DIPARTIMENTO DI MATEMATICA}
\vskip 1. truecm
\centerline{Dottorato di ricerca in matematica}
\vskip 2. truecm
\centerline{\bf AN INTRODUCTION TO SMALL DIVISORS PROBLEMS}
\vskip .5 truecm
\centerline{by}
\vskip .5 truecm
\centerline{Stefano Marmi}
\vfill\eject
%%%%% fino a qui %%%%%
\pageno=1
%%%%% prefazione %%%%%
\noindent
{\tit Preface }

\vskip 1. truecm
\noindent
The material treated in this book was brought together 
for a PhD course I taught at the University of Pisa
in the spring of 1999. 
It is  intended to 
be an introduction to small divisors problems. The book is 
divided in two parts. In the first one 
I discuss in some detail 
the theory of linearization of germs of
analytic diffeomorphisms of one complex variable. This is a part of 
the theory where many complete results are known. 
The second part is more informal. It 
deals with Nash--Moser's implicit function theorem in Fr\'echet spaces
and Kolmogorov--Arnol'd--Moser theory. Many results (and even 
some statements) are just briefly sketched but I always refer the 
reader to a choice of the huge original literature on the subject.

I am particularly fond of the topics described in the first part, especially 
because of their interplay with complex analysis and number theory. 
The second part is also fascinating both  because of its generality and because 
it leads to applications to Hamiltonian systems. 
Both are  the object of major active research. 

These lectures contain many problems (some of which may challenge 
the reader): they should be considered as an essential 
part of the text. The proof of many useful and important facts is left
as an exercise. 

I hope that the reader will find these notes a useful introduction 
to the subject. However the reason of the 
long list of references at the end of these 
notes is my belief that the best way to learn a subject
is to study directly the papers of those who invented it: 
Poincar\'e, Siegel, 
Kolmogorov, Arnol'd, Moser, Herman, Yoccoz, etc.

I am very grateful to Mariano Giaquinta for his invitation 
to give this series of lectures. I also wish to thank 
Carlo Carminati, whose enthusiasm
is also at the origin of this project, and whose remarks have been 
essential in correcting some mistakes. 

\vskip .3 truecm\noindent
Udine, December 8, 1999.

\rightline{\it Stefano Marmi}
%%%%% fine prefazione %%%%%
\vfill\eject
%%%%% contenuti %%%%%
%%%%%%%%% brj1.tex chap. 1 %%%%%%%%%%%%%%
%\Date
%\rightline{IHES, Bures--Sur--Yvette, France. September 17, 1999}
%\medskip
%\centerline{\tit An Introduction to Small Divisors Problems}
%\bigskip
%\centerline{S. Marmi\footnote{$^1$}{Dipartimento di Matematica e 
%Informatica, 
%Universit\`a di Udine, Via delle Scienze 206, Localit\`a Rizzi, 33100 
%Udine, Italy; 
%e-mail: marmi@dimi.uniud.it}}
%\vskip 2. truecm
%\centerline{\bf  Summary}

\bigskip
\centerline{\tit Table of Contents}

\vskip 1. truecm
{\titrm

\noindent
{\titsec PART I. One--dimensional Small Divisors. Yoccoz's Theorems}

\line{}

\noindent
\item{1.} Germs of Analytic Diffeomorphisms. Linearization 
\item{2.} Topological Stability vs. Analytic Linearizability
\item{3.} The Quadratic Polynomial: Yoccoz's Proof of the Siegel 
Theorem
\item{4.} Douady--Ghys' Theorem. Continued Fractions and the Brjuno 
Function
\item{5.} Siegel--Brjuno Theorem, Yoccoz's Theorem. Some Open Problems
\item{6.} Small divisors and loss of differentiability

\line{}

\noindent
{\titsec PART II. Implicit Function Theorems and KAM Theory}

\line{}

\noindent
\item{7.} Hamiltonian Systems and Integrable Systems 
\item{8.} Quasi--integrable Hamiltonian Systems
\item{9.} Nash--Moser's Implicit Function Theorem
\item{10.} From Nash--Moser's Theorem to KAM: Normal 
Form of Vector Fields on the Torus

\line{}

\noindent
{\titsec Appendices}

\line{}

\noindent
\item{A1.} Uniformization, Distorsion and Quasi--conformal maps
\item{A2.} Continued Fractions
\item{A3.} Distributions, Hyperfunctions, Formal Series. 
Hypoellipticity and Diophantine Conditions

\line{}

\noindent
References

\line{}

\noindent
Analytical index

\line{}

\noindent
List of symbols}
%\vskip 1. truecm
%%%%% fine contenuti %%%%%
\vfill\eject
%%%%% capitolo 1: germi %%%%%
\noindent
{\tit Part I. One--Dimensional Small Divisors. Yoccoz's 
Theorems}

\vskip 2. truecm
\noindent
{\titsec 1. Germs of Analytic Diffeomorphisms. Linearization }

\vskip .3 truecm
\noindent
A {\it dynamical system} is the action of a group (or a semigroup)
on some space. In looking for the simplest cases 
we are led to ask for the lowest possible dimension of the ambient 
space together with the highest possible regularity of the action. 
A remarkably rich but 
elementary situation is obtained considering the group of 
germs of holomorphic local diffeomorphisms of $\C$ which leave the 
point $z=0$ fixed. In what follows we will omit the symbol $\circ$
for the composition of two germs (unless some confusion may be 
possible). 

Let $\C[[z]]$ denote  the ring of formal power series and 
 $\C\{z\}$ denote  the ring 
of convergent power series. 

Let $G$ denote the group of germs of holomorphic diffeomorphisms of 
$(\C,0)$ and let $\hat{G}$ denote the group of formal germs of holomorphic 
diffeomorphisms of $(\C,0)$: $G=\{f\in z\C\{z\}\, , f'(0)\not= 0\}$, 
$\hat{G}=\{\hat{f}\in z\C[[z]]\, , \hat{f}_1\not= 0\}$.
One has the
trivial fibrations 
$$
\diagram{
G=\cup_{\lambda\in\C^{*}}G_{\lambda}& 
\phantom{\hfl{}{}} & \hat{G}=\cup_{\lambda\in\C^{*}}\hat{G}_{\lambda}
    \cr
\vfl{\pi}{}&&\vfl{\hat{\pi}}{} 
    \cr
\C^{*} & 
\phantom{\hfl{}{}} & \C^{*}
    \cr}
    \eqno(1.1)
$$
where 
$$
\eqalignno{
\hat{G}_{\lambda} &= \{\hat{f}(z)=\sum_{n=1}^\infty\hat{f}_{n}z^n
\in\C[[z]]\, , \; \hat{f}_{1}=\lambda\}\; , \; &(1.2)\cr
G_{\lambda} &= \{f(z)=\sum_{n=1}^\infty f_{n}z^n
\in\C\{z\}\, , \; f_{1}=\lambda\}\; . \; &(1.3)\cr
}
$$

\line{}

\vskip .3 truecm
\noindent
{\titsec 1.1 Conjugation, Symmetries}

\vskip .3 truecm
\noindent
Let $\hbox{Ad}_{g}\, f$ denote the adjoint action of $g$ on $f$: 
$\hbox{Ad}_{g}\, f = g^{-1}fg$. 

\vskip .3 truecm
\noindent
\Proc{Definition 1.1}{Let $f\in G$ (resp. $\hat{f}\in \hat{G}$). We 
say that a germ $g$ (resp. a formal germ $\hat{g}\in \hat{G}$)
is {\rm equivalent} or {\rm conjugate} to $f$ (resp. $\hat{f}$) if it 
belongs to the orbit of $f$ (resp. $\hat{f}$) under the adjoint action 
of $G_{1}$ (resp. $\hat{G}_{1}$):}
$$
\eqalign{
f\sim g &\iff \exists h\in G_{1}\, : \; g=h^{-1}fh\; , \cr
\hat{f}\sim \hat{g} &\iff \exists \hat{h}\in \hat{G}_{1}\, : \; 
\hat{g}=\hat{h}^{-1}\hat{f}\hat{h}\; . \cr
}
$$

\vskip .3 truecm
\noindent
The set of germs equivalent to $f$ obviously forms an equivalence class, the 
{\it orbit} of $f$ under the adjoint action of $G_{1}$: 
$$
[f] = \hbox{Ad}_{G_{1}} f = 
\{ g\in G\, , \; \exists h\in G_{1}\, : g=\hbox{Ad}_{h} f = 
h^{-1}fh\}\; . 
$$
The same holds in the formal case. 

\vskip .3 truecm
\noindent
\Proc{Definition 1.2}{A germ $g\in G$ is a {\rm symmetry} of $f\in G$ 
if $g\in \hbox{Cent}\, (f)$, i.e. if $\hbox{Ad}_{g}f=f$. We will 
denote by $\widehat{\hbox{Cent}}\, (\hat{f})$ the formal analogue of
$\hbox{Cent}\, (f)$.}

\vskip .3 truecm
\noindent
{\bf Exercise 1.3} Let $f\in G_{\lambda}$ (resp. 
$\hat{f}\in\hat{G}_{\lambda}$) and assume $g\sim f$
(resp. $\hat{g}\sim\hat{f}$), i.e. $f=h^{-1}gh$ for some $h\in 
G_{1}$. Then show that 
\item{(1)} $g\in G_{\lambda}$ (resp. $\hat{g}\in\hat{G}_{\lambda}$) 
thus $f_{1}=f'(0)=\lambda$ is invariant under conjugation.
\item{(2)} $\hbox{Cent}\, (f)$ is conjugated to $\hbox{Cent}\, (g)$, 
i.e. $\hbox{Cent}\, (f)=h^{-1}\hbox{Cent}\, (g)h$;
\item{(3)} $f^\Z = \{f^{ n}\, , \; n\in \Z\}\subset
\hbox{Cent}\, (f)$.

\line{}

\vskip .3 truecm
\noindent
{\titsec 1.2 Linearization}

\vskip .3 truecm
\noindent
Let $R_{\lambda}$ denote the germ $R_{\lambda}(z)=\lambda z$. This is 
the simplest element of $G_{\lambda}$. It is easy to check that, if 
$\lambda$ is not a root of unity, its 
centralizer is $\hbox{Cent}\,(R_{\lambda})=\{R_{\mu}\, , 
\;\mu\in\C^{*}\}$. 

\vskip .3 truecm
\noindent
{\bf Exercise 1.4} Let $f\in G_{\lambda}$ and assume that $\lambda$ 
is not a root of unity. The morphism 
$$
\eqalign{
\mu\, : \hbox{Cent}\, (f) &\rightarrow \C^{*}\cr
g &\mapsto \mu (g):=g_{1}=g'(0)\cr
}
$$
is injective. [Hint: this is equivalent to showing that $g\in G_{1}$, 
$g\in \hbox{Cent}\, (f)\Rightarrow g=\hbox{id}\,$. On the other hand
if $g\in G_{\mu}$ and $g\in \hbox{Cent}\, (f)$ one can recursively 
determine the power series coefficients of $g$:  one 
has
$$
(\lambda^n-\lambda )g_{n} = (\mu^n-\mu )f_{n}+\sum_{j=2}^{n-1}f_{j}
\sum_{n_{1}+\ldots n_{j}=n}g_{n_{1}}\cdots g_{n_{j}}-
\sum_{j=2}^{n-1}g_{j}\sum_{n_{1}+\ldots n_{j}=n}f_{n_{1}}\cdots 
f_{n_{j}}\; , 
$$
for all $n\ge 2$.]

\vskip .3 truecm
\noindent
\Proc{Definition 1.5}{A germ $f\in G_{\lambda}$ is {\rm linearizable}
if there exists $h_{f}\in G_{1}$ (a linearization of $f$) such that 
$h_{f}^{-1}fh_{f}=R_{\lambda}$, i.e. $f$ is conjugate to (its linear 
part) $R_{\lambda}$. $f$ is {\rm formally linearizable} if there 
exists $\hat{h}_{f}\in \hat{G}_{1}$ such that 
$\hat{h}_{f}^{-1}f\hat{h}_{f}=R_{\lambda}$
(note that in this case this is a functional equation in the ring 
$\C [[z]]$ of formal power series). }

\vskip .3 truecm
\noindent
From Exercise 1.4 it follows that when $\lambda$ is not a root of 
unity the linearization (if it exists) is unique: 
if $h_{1}$ and $h_{2}$ are two linearizations of the same $f\in 
G_{\lambda}$ then $h_{1}h_{2}^{-1}\in \ker\mu$. 

Our first  result on the existence of linearizations will concern the 
case when $\lambda$ is a root of unity. 

\vskip .3 truecm
\noindent
\Proc{Proposition 1.6}{Assume $\lambda$ is a primitive root of unity 
of order $q$. A germ $f\in G_{\lambda}$ is linearizable if and only if 
$f^{q} = \hbox{id}\,$. The same holds for a formal germ $\hat{f}\in
\hat{G}_{\lambda}$.}

\vskip .3 truecm
\noindent
\proof
Assume that $f$ is linearizable. Then 
$z=\lambda^q z= (h_{f}^{-1}\circ f\circ h_{f})^q (z) = 
(h_{f}^{-1}\circ f^q\circ h_{f})(z)$ from which one gets 
$f^q(z)= (h_{f}\circ \hbox{id}\,\circ h_{f}^{-1}) (z) = z$. 

Conversely if $f^q=\hbox{id}\,$ then defining 
$h_{f}^{-1} := {1\over q}\sum_{j=0}^{q-1}\lambda^{-j}f^j$ one 
immediately checks that $h_{f}^{-1}\in G_{1}$ if $f\in G_{\lambda}$
(resp. $h_{f}^{-1}\in \hat{G}_{1}$ if $f\in \hat{G}_{\lambda}$) 
and $h_{f}^{-1}\circ f\circ h_{f}=R_{\lambda}$. \qed

\line{}

\vskip .3 truecm
\noindent
{\titsec 1.3 Formal Conjugacy Classes}

\vskip .3 truecm
\noindent
In the formal case, all conjugacy classes of germs whose 
linear part is a root of 
unity are well known:

\vskip .3 truecm
\noindent
\Proc{Proposition 1.7}{Let $\lambda$ be a primitive root of unity 
of order $q$. Let $\hat{f}\in \hat{G}_{\lambda}$ and assume that 
$\hat{f}^{q} \not= \hbox{id}\,$. Then there exists a unique integer 
$n\ge 1$ and two complex numbers $a,b\in \C$, 
$a\not= 0$, such that $\hat{f}$ is formally 
conjugated to}
$$
P_{n,a,b,\lambda}(z)=\lambda z(1+az^{nq}+a^{2}bz^{2nq})\; . 
$$

\vskip .3 truecm
\noindent
{\bf Exercise 1.8 } Prove Proposition 1.7. Note that if one allows to 
conjugate also with homoteties then $\hat{f}$ is formally conjugated 
to $P_{n,c,\lambda} (z)=\lambda z(1+z^{nq}+cz^{2nq})$. 
[Hint: the idea of the proof is to iterate conjugations by polynomials
$\varphi_{j}(z)=z+\beta_{j}z^j$ with $j\ge 2$ and suitably chosen 
$\beta_{j}$. See also [Ar3], [Be].]

\vskip .3 truecm
\noindent
But in the formal case everything is very simple: 

\vskip .3 truecm
\noindent
\Proc{Proposition 1.9}{Assume that $\lambda$ is not a root of unity. Then 
$\hat{G}_{\lambda}$ is a conjugacy class and $\hat{G}_{1}$ acts freely and 
transitively on $\hat{G}_{\lambda}$.}

\vskip .3 truecm
\noindent
\proof
To see that any $f\in \hat{G}_{\lambda}$ is conjugate to $R_\lambda$
we look for $\hat{h}_f\in \hat{G}_{1}$ such that 
$\hat{f}\hat{h}_f = \hat{h}_fR_\lambda$. We develop and solve this 
functional equation by recurrence: we get, for $n\ge 2$
(denoting $\hat{h}_f(z) =\sum_{n=1}^\infty \hat{h}_nz^n$, 
$\hat{h}_1=1$)
$$
\hat{h}_n={1\over \lambda^n-\lambda}\sum_{j=2}^nf_j
\sum_{n_1+\ldots +n_j=n}\hat{h}_{n_1}\cdots\hat{h}_{n_j}\; 
. \eqno(1.4)
$$
The action of $\hat{G}_{1}$ on $\hat{G}_{\lambda}$
is free. This  follows from the fact that the only germ tangent 
to the identity 
belonging to the centralizer of $\hat{f}$ is the identity
(see Exercise 1.4). Transitivity of the action is trivial:
given two formal germs $\hat{f}_1$ and $\hat{f}_2$ both in 
$\hat{G}_{\lambda}$ there exist two formal linearizations 
$\hat{h}_1$ and $\hat{h}_2$ and clearly $\hat{f}_1=
\hbox{Ad}_{\hat{h}_2\hat{h}_1^{-1}}(\hat{f}_2)$. \qed

\vskip .3 truecm
\noindent
Collecting propositions 1.6, 1.7 and 1.9 together we have a complete 
classification of the conjugacy classes of $\hat{G}$:
\item{(I)} if $\lambda$ is not a root of unity then 
$\hat{G}_{\lambda}$ is a conjugacy class;
\item{(II)} if $\lambda =e^{2\pi i p/q}$, $q\ge 1$, 
$(p,q)=1$ then the conjugacy classes in $\hat{G}_{\lambda}$
are $[R_\lambda]$ and $\{[P_{n,a,b,\lambda}]\}_{a\in \C^{*}\, ,\, 
b\in\C\, ,\, n\ge 1}$.

\line{}

\vskip .3 truecm
\noindent
{\titsec 1.4 Koenigs--Poincar\'e Theorem}

\nobreak
\vskip .3 truecm
\noindent
In the holomorphic case the problem of a {\it complete}
classification of the conjugacy classes is still open and, as Yoccoz showed,
perhaps unreasonable. 
The first important result in the holomorphic case is the 
Koenigs--Poincar\'e Theorem:

\vskip .3 truecm
\noindent
\Proc{Theorem 1.10 (Koenigs--Poincar\'e)}{If $|\lambda|\not= 1$
then $G_\lambda$ is a conjugacy class, i.e. all $f\in G_\lambda$
are linearizable.}

\vskip .3 truecm
\noindent
\proof
Since $f$ is holomorphic around $z=0$ there exists $c_1>1$
and $r\in (0,1)$ such that $|f_j|\le c_1r^{1-j}$
for all $j\ge 2$. Since  $|\lambda|\not= 1$ there exists $c_2>1$
such that $|\lambda^n-\lambda|^{-1}\le c_2$ for all 
$n\ge 2$.

Let $(\sigma_n)_{n\ge 1}$ be the following recursively defined sequence:
$$
\sigma_1 = 1\; , \;\; 
\sigma_n=\sum_{j=2}^n\sum_{n_1+\ldots +n_j=n}
\sigma_{n_1}\cdots\sigma_{n_j}\; . \eqno(1.5)
$$
The generating function $\sigma (z)=\sum_{n=1}^\infty
\sigma_nz^n$ satisfies the functional equation 
$$
\sigma (z) = z+{\sigma (z)^2\over 1-\sigma (z)}\; , 
\eqno(1.6)
$$
thus 
$\sigma (z)={1+z-\sqrt{1-6z+z^2}\over 4}$
is analytic in the disk $|z|<3-2\sqrt{2}$
and bounded and continuous on its closure. By Cauchy's 
estimate one has 
$\sigma_n\le c_3(3-2\sqrt{2})^{1-n}$ for some $c_3>0$.

Since $\lambda$ is not a root of unity, $f$ is formally linearizable
and the power series coefficients of 
its formal linearization $\hat{h}_f$ satisfy (1.4). 
By induction one can check that $|\hat{h}_n|\le 
(c_1c_2r^{-1})^{n-1}\sigma_n$, thus $\hat{h}_f\in 
\C\{z\}$. \qed

\vskip .3 truecm
\noindent
{\bf Remark 1.11} Since the bound  $|\lambda^n-\lambda|^{-1}\le c_2$
is uniform w.r.t $\lambda\in D(\lambda_0,\delta )$, where $\lambda_0\in
\C^*\setminus\S^1$ and $\delta < \hbox{dist}\,(\lambda_0,\S^1)$, 
the above given proof of the Poincar\'e--Koenigs Theorem shows that the 
map 
$$
\eqalign{
\C^*\setminus\S^1 &\rightarrow G_1\cr
\lambda &\mapsto h_{\tilde{f}} (\lambda )\cr}
$$
is analytic\footnote{$^{1}$}{This notion needs a little comment
since $\C\{ z\}$ is a rather wild space: it is an inductive limit 
of Banach spaces, thus it is a locally convex topological vector space
and it is complete but it is not metrisable, thus it is not a 
Fr\'echet space (see Section 9.1). Here we simply mean that if 
$\lambda$ varies in some relatively compact open connected subset
of $\C^*\setminus\S^{1}$ then $h_{\tilde{f}} (\lambda )$
belongs to some fixed Banach space of holomorphic functions
(e.g. the Hardy space $H^\infty (\D_{r})$ of bounded analytic
functions on the disk $\D_{r}=\{z\in\C\, , \, |z|<r\}$, where
$r>0$ is fixed and small enough) and depends analytically on $\lambda$ in 
the usual sense.}
for all $\tilde{f}\in z^2\C\{z\}$, 
where $h_{\tilde{f}} (\lambda )$ is the linearization of 
$\lambda z+ \tilde{f}(z)$.

\vskip .3 truecm
\noindent
The Poincar\'e--Koenigs Theorem has the following 
straightforward generalization:

\vskip .3 truecm
\noindent
\Proc{Theorem 1.12 (Koenigs--Poincar\'e with parameters)}
{Let $r>0$, let $f\, : \D_r^n\times\D_r\subset \C^n\times\C
\rightarrow \C$, $(t,z)\mapsto f(t,z)=f_t(z)$ be an holomorphic 
map such that $f_0(z)=\lambda (0)z+\hbox{O}\, (z^2)$, 
with $|\lambda (0)|\not\in\{ 0,1\}$. Then there exists 
$r_0\in (0,r)$, a unique holomorphic function 
$z_0\,:\D_{r_0}^n\rightarrow \C$ and a unique 
$h\,:\D_{r_0}^n\times\D_{r_0}\rightarrow \C$, 
$(t,z)\mapsto h_t(z)=h(t,z)$ holomorphic such that for 
$t\in \D_{r_0}^n$ one has the following properties:
\item{(i)} $f_t(z_0(t))=z_0(t)$, 
$f_t(z)=\lambda (t)(z-z_0(t))+\hbox{O}\, ((z-z_0(t))^2)$,
$|\lambda (t)|\not\in\{ 0,1\}$;
\item{(ii)} $h_t(0)=z_0(t)$, $h_t'(0)={\partial\over\partial z}
h_t|_{z=0}=1$;
\item{(iii)} $h_t^{-1}\circ f\circ h_t=R_{\lambda (t)}$.}

\vskip .3 truecm
\noindent
\proof
(sketch) The existence of $z_0$ and (i) follows easily from the 
implicit function theorem applied to $F(t,z)=f(t,z)-z$
at the point $(t,z)=(0,0)$ 
(note that $F(0,0)=0$ and ${\partial\over\partial z}f_t(z)|_{(t,z)=
(0,0)}=\lambda (0)-1\not= 0$).
Therefore there exists a unique fixed point 
for $f_t$ close to $z=0$ when $t$ is close to $0$  
depending analytically on $t$ as $t$ varies in a neighborhood of 
$(t,z)=(0,0)$. Then one can consider $g_t(z)=f_t(z+z_0(t))-z_0(t)$
and apply the proof given above of the  Koenigs--Poincar\'e Theorem to 
$g_t (z)$. It is easy to convince oneself that the linearizing 
map depends analytically on $t$. \qed

\line{}

\vskip .3 truecm
\noindent
{\titsec 1.5 Centralizers and Linearizations}

\vskip .3 truecm
\noindent
The study of centralizers generalizes the study of 
linearizability as the following exercises show:

\vskip .3 truecm
\noindent
{\bf Exercise 1.13 } Prove that if $f\in G_\lambda$ is 
linearizable and $\lambda$ is not a root of unity then 
$\hbox{Cent}\, (f)\simeq \C^*$. [Hint: use the fact that the 
centralizer of $f$ is conjugate to the centralizer of $R_\lambda$
which is  completely known.]

\vskip .3 truecm
\noindent
{\bf Exercise 1.14}  Prove that if $g\in \hbox{Cent}\, (f)$, 
$g\in G_\mu$ is linearizable and $\mu$ is not a root of unity
then $f$ is linearizable. [Hint: use that   $f\in 
\hbox{Cent}\, (g)=\{h_gR_\nu h_g^{-1}\, , \, \nu\in\C^*\}$
and that $\nu$ is invariant under conjugacy.]

\vskip .3 truecm
\noindent
{\bf Exercise 1.15}  Prove that if $f\in G_\lambda$
 and $\lambda$ is not a root of unity then $f$ is linearizable 
if and only if $\hbox{Cent}\, (f)\simeq \C^*$. [Hint: apply
exercises 1.13, 1.4 and the Koenigs--Poincar\'e Theorem]

\line{}

\vskip .3 truecm
\noindent
{\titsec 1.6 Cremer's Non--Linearizable Germs}

\vskip .3 truecm
\noindent
When $|\lambda |=1$ and $\lambda$ is not a root of unity we can write
$$
\lambda = e^{2\pi i \alpha}\;\;\hbox{with}\;
\alpha\in\R\setminus\Q\cap (-1/2,1/2)\; ,
\eqno(1.7)
$$
and whether  $f\in G_\lambda$ is linearizable or not depends crucially 
on the arithmetical properties of $\alpha$. Let $\{x\}$ denote the fractional 
part of a real number $x$: $\{x\}=x-[x]$, where $[x]$ is the integer 
part of $x$.

\vskip .3 truecm
\noindent
\Proc{Theorem 1.16 (Cremer)}{If $\;\limsup_{n\rightarrow +\infty} 
|\{n\alpha\}|^{-1/n}=+\infty$ then there exists $f\in G_{e^{2\pi i \alpha}}$
which is not linearizable.}

\vskip .3 truecm
\noindent
\proof
First of all note that 
$\limsup_{n\rightarrow +\infty} 
|\{n\alpha\}|^{-1/n}=+\infty$ if and only if 
$$
\limsup_{n\rightarrow +\infty} 
|\lambda^n-1|^{-1/n}=+\infty
$$ 
since 
$$
|\lambda^n-1|=2|\sin (\pi n\alpha )|\in 
(2|\{n\alpha\}|,\pi |\{n\alpha\}|)\; .
$$
Then we construct $f$ in the following manner: for $n\ge 2$
we take $|f_n|=1$ and we choose inductively $\arg f_n$
such that 
$$
\arg f_n = \arg \sum_{j=2}^{n-1}f_j
\sum_{n_1+\ldots +n_j=n}\hat{h}_{n_1}\cdots\hat{h}_{n_j}\; ,
\eqno(1.8)
$$
(recall the induction formula (1.4) for the coefficients of the 
formal linearization of $f$ and note that the r.h.s. of (1.8) is a 
polynomial in $n-2$ variables $f_2,\ldots ,f_{n-1}$ with 
coefficients in the field $\C (\lambda )$).
Thus 
$$
|\hat{h}_{n}|\ge {|f_n|\over |\lambda^n-1|} = 
{1\over |\lambda^n-1|}
$$
and $\limsup_{n\rightarrow +\infty} 
|\hat{h}_{n}|^{1/n}=+\infty$: the formal linearization 
$\hat{h}$ is a divergent series. \qed

\vskip .3 truecm
\noindent
{\bf Exercise 1.17} Write the decimal expansion of an irrational 
number $\alpha$ satisfying  the assumption of Cremer's Theorem.

\vskip .3 truecm
\noindent
{\bf Exercise 1.18} Show that the set of irrational numbers satisfying the 
assumption of Cremer's Theorem is a dense $G_\delta$ with zero Lebesgue 
measure (following Baire, 
a set is a dense $G_\delta$ if it is a countable intersection of 
dense open sets. These sets are ``big'' from the point of view of topology).

\vskip .3 truecm
\noindent
In the next Chapter we will continue our study of the problem of 
the existence of a linearization of germs of holomorphic diffeomorphisms. 
To this purpose the following ``normalization'' will be useful.

\line{}

\vskip .3 truecm
\noindent
{\titsec 1.7 Normalized Germs }

\vskip .3 truecm
\noindent
Let us note that there is an obvious action of 
$\C^*$ on $G$ by homotheties:
$$
(\mu ,f)\in \C^*\times G\mapsto \hbox{Ad}_{R_\mu}f=R_\mu^{-1}fR_\mu\; .
\eqno(1.9)
$$
Note that this action leaves the fibers $G_\lambda$ invariant 
by Exercise 1.3. Also, 
$f\in G_\lambda$ is linearizable if and only if $ \hbox{Ad}_{R_\mu}f$
is also linearizable for all $\mu\in \C^*$
(indeed if $h_f$ linearizes $f$ then $\hbox{Ad}_{R_\mu}h_f$
linearizes $\hbox{Ad}_{R_\mu}f$). Therefore, 
in order to study the problem of 
the existence of a linearization, it is enough to consider 
$G/\C^*$, i.e. we identify two germs of holomorphic diffeomorphisms 
which are conjugate by a homothety.

Consider the space $S$ of univalent maps $F\,:\D\rightarrow\C$
such that $F(0)=0$ and the projection 
$$
\eqalign{
G &\rightarrow S\cr
f &\mapsto F=\cases{
f &if $f$ is univalent in $\D$\cr
 \hbox{Ad}_{R_r}f &if $f$ is univalent in $\D_r$\cr}
\cr}
$$
This map is clearly onto and two germs have the same image only if 
they coincide or if they are conjugate by some homothety. Thus this 
projection induces a bijection from $G/\C^*$ onto $S$. 

In what follows we will always consider the topological 
space $S$ of germs 
of holomorphic diffeomorphisms $f\,:\D\rightarrow\C$
such that $f(0)=0$ and $f$ is univalent in $\D$. We will denote 
\item{$\bullet$}
$S_\lambda$ the subspace of $f$ such that $f'(0)=\lambda$;
\item{$\bullet$}
$S_{\T}$ the subspace of $f$ such that $|f'(0)|=1$.
 
Clearly the projection above induces a bijection between 
$G_\lambda/\C^*$ and $S_\lambda$.
%%%%% fine capitolo 1 %%%%%
\vfill\eject
%\vskip 1. truecm
%%%%% capitolo 2: stabilita' %%%%%
\noindent
{\tit 2. Topological Stability vs. Analytic
Linearizability}

\vskip .3 truecm
\noindent
The purpose of this Chapter is to connect the study of 
the conjugacy classes of germs of holomorphic diffeomorphisms 
to the theory of one--dimensional conformal dynamical systems 
and in particular to the notion of stability of a fixed point.
The extremely remarkable fact is that stability, which is a topological 
property, will turn out to be equivalent to linearizability, which 
is an analytic property.

\line{}

\vskip .3 truecm
\noindent
{\titsec 2.1 Dynamics of Rational Maps}

\vskip .3 truecm
\noindent
Let us first of all recall the notion of normal family on an
open subset $U$ of the Riemann sphere $\overline{\C}=\C\cup
\{\infty\}$. To this purpose we recall the usual system 
of coordinates on $\overline{\C}$ determined by the stereographic projection:
$z\,:\overline{\C}\setminus\{\infty\}\rightarrow\C$, $z(0)=0$ and
$w\,:\overline{\C}\setminus\{0\}\rightarrow\C$, $w(\infty )=0$,
related by $zw=1$.
The {\it spherical metric} on $\overline{\C}$ is defined as follows:
$$
ds_{\overline{\C}}=\cases{
{2|dz|\over 1+|z|^2} &in the $z$--chart;\cr
{2|dw|\over 1+|w|^2} &in the $w$--chart;\cr
}
\eqno(2.1)
$$
Let $U\subset \overline{\C}$ be open and ${\cal F}_U=\{f\,:U\rightarrow
\overline{\C}\, , \; f \;\hbox{meromorphic}\}$. We endow 
$\overline{\C}$ with the spherical metric and ${\cal F}_U$ with the 
topology of uniform convergence on compact subsets of $U$. It is a 
classical result of Weierstrass that the limit of a convergent sequence in 
${\cal F}_U$ still belongs to ${\cal F}_U$ (note that the 
constant function $f\equiv \infty$ is considered meromorphic).

\vskip .3 truecm
\noindent
\Proc{Definition 2.1}{A family ${\cal F}\subset {\cal F}_U$ is 
{\rm normal} if it is relatively compact in ${\cal F}_U$, 
i.e. any sequence $\{f_n\}\subset {\cal F}$ contains a subsequence
which converges uniformly in the spherical metric on compact subsets
of $U$.}

\vskip .3 truecm
\noindent
{\it Warning!} If $\{f_n\}$ is a normal family 
then $\{f_n'\}$ needs not  be normal: e.g.
$f_n(z)=n(z^2-n)$ on $\C$.

By means of the Ascoli--Arzel\`a theorem one gets:

\vskip .3 truecm
\noindent
\Proc{Proposition 2.2}{
\item{(I)} A family of meromorphic functions on $U$ is normal on $U$
if and only if it is equicontinuous on every compact subset of $U$;
\item{(II)} A family of analytic functions on $U$ is normal on $U$
if and only if it is locally uniformly bounded (i.e. uniformly 
bounded on every compact subset of $U$).}

\vskip .3 truecm
\noindent
\proof
The first statement is obvious since the compactness of 
$\overline{\C}$ guarantees that the family is uniformly bounded.
The second statement follows from Cauchy's integral theorem.
\qed

\vskip .3 truecm
\noindent
The notion of normal family allows us to introduce the basic notions 
of one--dimensional holomorphic dynamics. Here we are interested in 
studying the dynamics of a discrete dynamical system 
(i.e. an action of $\N$) on the Riemann sphere $\overline{\C}$
generated by a holomorphic transformation $R\,:
\overline{\C}\rightarrow \overline{\C}$, i.e. and element of 
$\hbox{End}\,(\overline{\C})$. 

Let $d$ denote the topological 
degree of $R$. We will assume $d\ge 2$ thus $R$ is a $d$--fold 
branched covering of the Riemann sphere and can be written 
in a unique way in the form 
$R(z)={P(z)\over Q(z)}$, where $P(z)\in\C [z]$,
$Q(z)\in\C [z]$ have no common factors and $d=\max (
\hbox{deg}\, P,\hbox{deg}\, Q)$. In fact every 
$d:1$ conformal branched covering of 
$\overline{\C}$ comes from some such rational function and 
$$
\hbox{End}\,(\overline{\C})=\{
R\,:\overline{\C}\rightarrow \overline{\C} \,
\hbox{holomorphic}\} = \{
R\,:\overline{\C}\rightarrow \overline{\C} \, , 
\, R(z)=P(z)/Q(z)\}\; .
$$
Note that 
$$
R(z)=\cases{
{P(z)\over Q(z)} & if $Q(z)\not= 0$,\cr
\infty & if $Q(z)=0$,\cr
\lim_{z\rightarrow\infty}{P(z)\over Q(z)} & if $z=\infty$.\cr
}
$$
We define the iterates $R^n$ of $R$ as usual: $R^n=R\circ R^{n-1}$.
Note that $R^n$ has degree $d^n$.

Given a point $z_0\in\overline{\C}$ the sequence of points 
$\{z_n\}_{n\ge 0}$ defined by $z_{n+1}=R(z_n)$
is called the orbit of $z_0$. A point $z_0$ is a {\it fixed point}
of $R$ if $R(z_0)=z_0$, {\it periodic} if $z_n=R^n(z_0)=z_0$
for some $n$ (the minimal $n$ is the period). The orbit 
$\{z_1,\dots ,z_n=z_0\}$ is called a {\it cycle}. The point 
$z_0$ is called {\it preperiodic} if $z_k$ is periodic for some $k>0$.

The fundamental dichotomy of $\overline{\C}$ associated to the 
dynamics of $R$ is the following:

\vskip .3 truecm
\noindent
\Proc{Definition 2.3}{The {\rm Fatou set} $F(R)$ of $R$ is the set of 
points $z_0\in \overline{\C}$ such that $\{R^n\}_{n\ge 0}$ is a 
normal family in some disk $D(z_0,r)$ (w.r.t. the spherical 
metric). The complement of the Fatou set is the {\rm Julia set}
$J(R)$.}

\vskip .3 truecm
\noindent
{\bf Exercise 2.4} Show that $J(R)=J(R^n)$ for all $n\ge 1$ and that 
$J(R)$ is nonempty and closed. [Hint: if $J(R)=\emptyset$
then $\{R^n\}_{n\ge 0}$ is a normal family on all $\overline{\C}$
thus $R^{n_j}\rightarrow S\in\hbox{End}\,(\overline{\C})$. Compare 
degrees.]

\line{}

\vskip .3 truecm
\noindent
{\titsec 2.2 Stability}

\vskip .3 truecm
\noindent
\Proc{Definition 2.5}{A point $z_0\in\overline{\C}$ is {\rm stable}
if for all $\delta >0$ there exists a neighborhood 
$W$ of $z_0$ such that for all $z\in W$ and for all $n\ge 0$
one has $d(R^n(z),R^n(z_0))\le \delta$
(here, as usual, $d$ denotes the spherical metric).}

\vskip .3 truecm
\noindent
{\bf Exercise 2.6} Show that a point is stable if and only 
if it belongs to the Fatou set. 
%[Hint: the less trivial part is to 
%show that all the points in the Fatou set are stable. 

\vskip .3 truecm
\noindent
{\bf Exercise 2.7} Let $R\in \hbox{End}\, (\overline{\C})$
and assume that $R(0)=0$.  If $R$ is linearizable and 
$R'(0)=\lambda$, $|\lambda |\le 1$ then $0$ is a stable fixed point.

\vskip .3 truecm
\noindent
If we consider the more general situation of a germ $f\in S$, 
i.e. $f\,:\D\rightarrow \C$ injectively and holomorphically, $f(0)=0$, 
the definition of stability must be slightly generalized so as to take
into account the fact that the iterates of $f$ are not necesaarily 
defined for all $n$.

\vskip .3 truecm
\noindent
\Proc{Definition 2.8}{$0$ is {\rm stable} if and only if 
there exists a neighborhood $U$ of $0$ such that 
$f^n$ is defined on $U$ for all $n\ge 0$ and for all 
$z\in U$ and $n\ge 0$ one has $|f^n(z)|<1$.}

\vskip .3 truecm
\noindent
{\bf Exercise 2.9} Show that if $f$ is a rational map with a fixed point $0$
then definitions 2.8 and  2.5 are equivalent.

\vskip .3 truecm
\noindent
{\bf Exercise 2.10} If $f'(0)=\lambda$ and $|\lambda |>1$ 
then $0$ is not stable.

\vskip .3 truecm
\noindent
To each germ $f\in S$, $|f'(0)|\le 1$, 
one can associate a natural $f$--invariant compact set 
$$
0\in K_f := \bigcap_{n\ge 0} f^{-n}(\D)\; . \eqno(2.2)
$$
Let $U_f$ denote the connected component of the interior of 
$K_f$ which contains $0$. Then $0$ is stable if and only if
$U_f\not=\emptyset$, i.e. if and only if $0$ belongs to the interior 
of $K_f$.

\vskip .3 truecm
\noindent
{\bf Exercise 2.11} Show that  if $f\in S$ and $|f'(0)|<1$ then $0$
is stable. [Hint: consider a small disk around $0$ on which 
the inequality $|f(z)|\le\rho |z|$ with $\rho <1$ holds.]

\line{}

\vskip .3 truecm
\noindent
{\titsec 2.3 Stability vs. Linearizability}

\vskip .3 truecm
\noindent
The main result of this section is the equivalence (for $|f'(0)|\le 1$)
of stability (a topological notion) and linearizability (an analytical 
notion):

\vskip .3 truecm
\noindent
\Proc{Theorem 2.12}{Let $f\in S$, $|f'(0)|\le 1$. $0$ is stable if and 
only if $f$ is linearizable.}

\vskip .3 truecm
\noindent
\proof
The statement is non--trivial only if $\lambda=f'(0)$ has unit modulus.
If $f$ is linearizable then the linearization $h_f$ maps a small disk 
$\D_r$ around zero conformally into $\D$. Since $h_f(0)=0$
and $|f^n(z)|<1$ for all $z\in h_f(\D_r)$ one 
sees that $0$ is stable. 

Conversely assume now that $0$ is stable. Then $U_f\not=\emptyset$
and one can easily see that it must also be simply connected 
(otherwise, if it had a hole $V$, surrounding it with some closed 
curve $\gamma$ contained in $U_f$ since 
$|f^n(z)|<1$ for all $z\in\gamma$ and $n\ge 0$ the maximum principle 
leads to the same conclusion for all the points in $V$ thus 
$V\subset U_f$). Applying the Riemann mapping theorem to $U_f$
one sees that  by conjugation with the Riemann map $f$ induces 
a univalent map $g$ of the disk into itself with the same linear part
$\lambda$. By Schwarz' Lemma one must have $g(z)=\lambda z$
thus $f$ is analytically linearizable.
\qed

\vskip .3 truecm
\noindent
When $\lambda =f'(0)$ has modulus one, is not a root of unity and $0$ is stable
then $U_f$ is conformally equivalent to a disk and is called 
the {\it Siegel disk} of $f$ (at $0$). 
Thus the Siegel disk of $f$ is the maximal connected open set containing
$0$ on which $f$ is conjugated to $R_\lambda$. 
The conformal representation 
$\tilde{h}_f\,:\D_{c(f)}\rightarrow U_f$ of $U_f$ which 
satisfies $\tilde{h}_f(0)=0$, $\tilde{h}_f'(0)=1$ linearizes
$f$ thus the power series of $\tilde{h}_f$ and $h_f$ coincide.
If $r(f)$ denotes the radius of convergence of the linearization 
$h_f$ (whose power series coefficients are recursively determined 
as in (1.4)), recalling the definition of conformal capacity
(Exercise A1.4, Appendix 1) we see that:
\item{(i)} if $r(f)>0$ then $0<c(U_f,0)=c(f)\le r(f)$;
\item{(ii)} if $r(f)=0$ then $c(f)=r(f)=0$.

\vskip .3 truecm
\noindent
{\bf Exercise 2.13} Show that the map $c\,:S_{\T}\rightarrow [0,1]$
which associates to each germ $f$ the conformal capacity of $U_f$ 
w.r.t.~$0$ is upper semicontinuous.

\vskip .3 truecm
\noindent
When $0$ is not stable and $\lambda$ is not a root of unity 
$K_f$ is called a {\it hedgehog}.

We conclude this introduction to Siegel disks with two results on 
their conformal capacity.

\vskip .3 truecm
\noindent
\Proc{Proposition 2.14}{One has $c(f)=r(f)$ when at least one of the
two following conditions is satisfied:
\item{(i)} $U_f$ is relatively compact in $\D$;
\item{(ii)} each point of $\S^1$ is a singularity of $f$.}

\vskip .3 truecm
\noindent
\Proc{Proposition 2.15}{Let $\lambda\in\T$ and assume that $\lambda$
is not a root of unity. Then}
$$
\inf_{f\in S_\lambda}c(f)=\inf_{f\in S_\lambda}r(f)\; .
$$

\vskip .3 truecm
\noindent
For the proofs see [Yo2, p.19].
%%%%% fine capitolo 2 %%%%%
\vfill\eject
%\vskip 1. truecm
%%%%% capitolo 3: polinomio quadratico %%%%%
\noindent
{\tit 3. The Quadratic Polynomial. Yoccoz's Proof of 
Siegel Theorem} \par
\vskip .3 truecm
\noindent
In this Chapter we will show the special role played by the quadratic 
polynomial 
$$
P_\lambda (z)=\lambda\left(z-{z^2\over 2}\right)\; . 
\eqno(3.1)
$$
Indeed $P_\lambda$ is the ``worst possible perturbation
of the linear part $R_\lambda$'' as the following theorem 
shows 

\vskip .3 truecm
\noindent
\Proc{Theorem 3.1 (Yoccoz)}{Let $\lambda=e^{2\pi i\alpha}$, 
$\alpha\in\R\setminus\Q$. If $P_\lambda$ is linearizable then 
every germ $f\in G_\lambda$ is also linearizable.}

\vskip .3 truecm
\noindent
\proof
Let $f\in G_\lambda$, $f(z)=\lambda z+\sum_{n=2}^\infty
f_nz^n$. By conjugating with some homothety one has 
$|f_n|\le 10^{-3}4^{-n}$. We now consider the one--parameter 
family $f_b(z)=\lambda z+bz^2+\sum_{n=2}^\infty
f_nz^n$. Note that $f_0=f$. Since $\lambda$ is not a root 
of unity there exists a unique formal germ $\hat{h}_b\in
\hat{G}_1$ such that $\hat{h}_b^{-1}f_b\hat{h}_b=R_\lambda$. 
Its power series expansion is 
$\hat{h}_b(z)=z+\sum_{n=2}^\infty h_n(b)z^n$ with 
$h_n(b)\in\C[b]$. Thus by the maximum principle one has 
$|h_n(0)|\le\max_{|b|=1/2}|h_n(b)|$.
If $|b|=1/2$ then (possibly after conjugation with a
rotation) $f_b(z)=P_\lambda (z)+\sum_{n=2}^\infty
f_nz^n=P_\lambda (z)+\psi (z)$ and it is immediate to 
check that $\sup_{z\in\D_3}|\psi (z)|<10^{-2}$. 
From Douady--Hubbard's theorem on the stability of the 
quadratic polynomial (Appendix 1) it follows that 
$f_b$ is quasiconformally conjugated to $P_\lambda$. If 
$P_\lambda$ is linearizable then $0$ is stable for 
$P_\lambda$ thus also for $f_b$ since a quasiconformal conjugacy
is in particular a topological conjugacy. But we know that 
this implies that $f_b$ is linearizable. Therefore there exists 
two positive constants $C$ and $r$ such that 
$|h_n(b)|\le Cr^{-n}$ for all $b$ of modulus $1/2$, 
thus $|h_n(0)|\le Cr^{-n}$. Then $\hat{h}_0$ converges
and $f_0=f$ is linearizable. \qed

\line{}

\vskip .3 truecm
\noindent
{\titsec 3.1 Yoccoz's Linearization Theorem for the Quadratic Polynomial}

\vskip .3 truecm
\noindent
Once one has established that the linearizability 
of the quadratic polynomial for a certain $\lambda$
implies that $G_\lambda$ is a conjugacy class the following 
remarkable theorem of Yoccoz shows that $G_\lambda$ is a 
conjugacy class for almost all $\lambda\in\T$.

\vskip .3 truecm
\noindent
\Proc{Theorem 3.2}{Let  $\lambda=e^{2\pi i\alpha}$, 
$\alpha\in\R\setminus\Q$. For almost all $\lambda\in\T$
the quadratic polynomial $P_\lambda$ is linearizable.}

\vskip .3 truecm
\noindent
This statement deserves a comment. As we will see in Chapter 5
this theorem of Yoccoz is indeed weaker than the Siegel [S] and Brjuno 
[Br] theorems
which date respectively to 1942 and 1970. What is {\it very remarkable}
is the proof of Theorem 3.2 which does not need any subtle estimate 
on the growth of the coefficients of the formal linearization 
as provided by (1.4): compare with the proof of the  Siegel--Brjuno Theorem
given in Section 5.1.

Let us note that $P_\lambda$ has a unique critical point
$c=1$ (apart from $z=\infty$) and that the corresponding 
critical value is $v_\lambda=P_\lambda (c)=\lambda /2$. 
If $|\lambda|<1$ by Koenigs--Poincar\'e theorem we know 
that there exists a unique analytic linearization 
$H_\lambda$ of $P_\lambda$ and that it depends analytically 
on $\lambda$ as $\lambda$ varies in $\D$. Let 
$r_2(\lambda )$ denote the radius of convergence of $H_\lambda$. 
One has the following

\vskip .3 truecm
\noindent
\Proc{Proposition 3.3}{Let $\lambda\in\D$. Then:
\item{(1)} $r_2(\lambda )>0$;
\item{(2)} $r_2(\lambda )<+\infty$ and $H_\lambda$ has a 
continuous extension to $\overline{\D_{r_2(\lambda )}}$. Moreover 
the map $H_\lambda\,:\overline{\D_{r_2(\lambda )}}
\rightarrow\C$ is conformal and verifies $P_\lambda\circ
H_\lambda =H_\lambda\circ R_\lambda$.
\item{(3)} On its circle of convergence  
$\{z\, ,\, |z|=r_2(\lambda )\}$, $H_\lambda$ has a unique singular
point which will be denoted $u(\lambda )$.
\item{(4)} $H_\lambda (u(\lambda ))=1$ and 
$(H_\lambda (z)-1)^2$ is holomorphic in 
$z=u(\lambda )$.}

\vskip .3 truecm
\noindent
\proof
The first assertion is just a consequence of Koenigs--Poincar\'e theorem.

The functional equation $P_\lambda (H_\lambda (z))=H_\lambda (\lambda z)$
is satisfied for all $z\in\D_{r_2(\lambda )}$. Moreover $H_\lambda\,:
\D_{r_2(\lambda )}\rightarrow\C$ is univalent (if one had 
$H_\lambda (z_1)=H_\lambda (z_2)$ with $z_1\not=z_2$ and 
$z_1,z_2\in\D_{r_2(\lambda )}$ one would have $H_\lambda
(\lambda^nz_1)=H_\lambda
(\lambda^nz_2)$ for all $n\ge 0$ which is impossible 
since $|\lambda |<1$ and $H_\lambda'(0)=1$). Thus 
$r_2(\lambda )<+\infty$. On the other hand if $H_\lambda$ is 
holomorphic in $\D_r$ for some $r>0$ and the critical value 
$v_\lambda\not\in H_\lambda (\D_r)$ the functional equation 
allows to continue analytically $H_\lambda$ to the disk 
$\D_{|\lambda |^{-1}r}$. Therefore there exists $u(\lambda )\in
\C$ such that $|u(\lambda )|=r_2(\lambda )$ and 
$H_\lambda (\lambda u(\lambda ))=v_\lambda$.
Such a $u(\lambda )$ is unique since $H_\lambda$ is 
injective on $\D_{r_2(\lambda )}$. If 
$|w|=|\lambda |r_2(\lambda )$ and $w\not=\lambda u(\lambda )$
one has 
$H_\lambda (w)=P_\lambda (H_\lambda (\lambda^{-1}w))$
and 
$$
H_\lambda (\lambda^{-1} w)=1-\sqrt{1-2\lambda^{-1}H_\lambda
(w)}\eqno(3.2)
$$
which shows how to extend continuously and injectively $H_\lambda$
to $\overline{\D_{r_2(\lambda )}}$. By construction the 
functional equation is trivially verified. This completes 
the proof of (2).

To prove (3) and (4) note that from $H_\lambda (\lambda u(\lambda ))
= P_\lambda (H_\lambda (u(\lambda )))$ it follows that 
$H_\lambda (u(\lambda ))=1$. Formula (3.2) shows that all points 
$z\in\C$, $|z|=r_2(\lambda )$ are regular except for 
$z=u(\lambda )$. Finally one has $(H_\lambda (z)-1)^2=
1-2\lambda^{-1}H_\lambda (\lambda z)$ which is holomorphic 
also at $z=u(\lambda )$. \qed

\vskip .3 truecm
\noindent
The fact that $H_\lambda$ is injective on $\overline{\D}_{r_2(\lambda )}$
implies that $r_2 (\lambda )<+\infty$ (otherwise it would be 
a biholomorphism of $\C$ thus an affine map). A more precise upper 
bound is provided by the following 

\vskip .3 truecm
\noindent
\Proc{Lemma 3.4 (a priori estimate of $r_2(\lambda )$).}
{$r_2(\lambda )\le 2$.}

\vskip .3 truecm
\noindent
\proof
It is an easy consequence of Koebe $1/4$--Theorem. Indeed if 
$\tilde{f}\in S_1$ and $t>0$ then $f=R_t^{-1}\tilde{f}R_t\,:
\D_{t^{-1}}\rightarrow\C$ is univalent and $f(\D_{t^{-1}})=
R_{t^{-1}}\tilde{f}(\D )$. By Koebe $1/4$--Theorem one has 
$\D_{1/4}\subset\tilde{f}(\D )$ thus $\D_{t^{-1}/4}\subset
f(\D_{t^{-1}})$. But we know that $v_\lambda\not\in
H_\lambda (\D_{|\lambda |r_2(\lambda )})$ thus
$|v_\lambda |={|\lambda |\over 2}\ge {|\lambda |r_2(\lambda )
\over 4}$. \qed

\vskip .3 truecm
\noindent
{\bf Exercise 3.5} Show that $r_2(\lambda )\le 8/7$. 
[Hint: apply (A1.1) to the function 
$$
\tilde{H}_\lambda (z)\left(1+{2r_2(\lambda )\over\lambda}
\tilde{H}_\lambda (z)\right)^{-1}\in S_1\; , 
$$
where $\tilde{H}_\lambda (z)={H_\lambda (r_2(\lambda )z)\over
r_2(\lambda )}$.]

\vskip .3 truecm
\noindent
{\bf Exercise 3.6} Show that the image by $H_\lambda$ of its circle 
of convergence is a Jordan curve, analytic except at $H_\lambda (u(\lambda
))=1$ where it has a right angle.

\vskip .3 truecm
\noindent
\Proc{Proposition 3.7}{$u\,:\D^*\rightarrow\C$
has a bounded analytic extension to $\D$. Moreover it is 
the limit of the sequence of polynomials 
$u_n(\lambda )=\lambda^{-n}P_\lambda^n(1)$
uniformly on compact subsets of $\D$. One has $u(0)=1/2$.}

\vskip .3 truecm
\noindent
\proof
By Proposition 3.3 one has $P_\lambda^n(1)=H_\lambda (\lambda^n
u(\lambda ))$. From Lemma 3.4 and Koebe's distorsion 
estimates (specifically (A1.4) ) applied to 
$\tilde{H}_\lambda (z)={H_\lambda (u(\lambda )z)\over u(\lambda )}$
one has 
$$
|P_\lambda^n (1)|=|u(\lambda )\tilde{H}_\lambda (\lambda^n)|
\le r_2(\lambda ){|\lambda |^n\over (1-|\lambda|^n)^2}\le
2{|\lambda |^n\over (1-|\lambda|^n)^2}\; , 
$$
thus $|u_n(\lambda )|\le 2(1-|\lambda |)^{-2}$ 
for all $\lambda\in\D$ and the polynomials 
$u_n$ verify the recurrence relation 
$$
u_0(\lambda )=1\; , \;\;\;
u_{n+1}(\lambda )=u_n(\lambda )-{\lambda^n\over 2}
(u_n(\lambda ))^2\; . \eqno(3.3)
$$
This shows that $u_n$ converges uniformly on compact subsets 
of $\D$. The limit is $u$ since 
$$
\lim_{n\rightarrow +\infty}u_n(\lambda )=
\lim_{n\rightarrow +\infty}\lambda^{-n}H_\lambda (\lambda^nu(
\lambda ))=u(\lambda )\; . 
$$
Finally from $|u(\lambda )|=r_2(\lambda )$ and Lemma 3.4 one has 
$|u(\lambda )|\le 2$ on $\D$. 
\qed

\vskip .3 truecm
\noindent
The function $u\, :\D\rightarrow\C$ will be called {\it Yoccoz's function}. 
It has many remarkable properties and it is the object of various 
conjectures (see Section 5.3). 

\vskip .3 truecm
\noindent
{\bf Exercise 3.8} Check that:
\item{1)} $u(\lambda )-u_{n}(\lambda )=\hbox{O}\,(\lambda^n)$;
\item{2)} $u(\lambda )={1\over 2}-{\lambda\over 8}-{\lambda^2\over
8}-{\lambda^3\over 16}-{9\lambda^{4}\over 128}-{\lambda^{5}\over 128}
-{7\lambda^{6}\over128}+{3\lambda^{7}\over 256}-{29\lambda^{8}\over 
1024}-{\lambda^{9}\over 256}+{25\lambda^{10}\over 
2048}+{559\lambda^{11}\over 32768}+
\ldots$;
\item{3)} $u(\lambda )\in\Q\{\lambda \}$ and all 
the denominators are a power of $2$.

\noindent
Write a computer program to calculate the power series 
expansion of $u$ and use it to design the level 
sets of $\log |u|$ and $\arg u$. 
Try to compute the graph of $\theta\mapsto\arg u(re^{2\pi i\theta })$
as $r\rightarrow 1-$. You may use some formulas given in [Yo2, pp. 
70--71] and to compare with [MMY2]. If you get nice pictures I would 
like to get a copy of them. 

\line{}

\vskip .3 truecm
\noindent
{\titsec 3.2 Radial Limits of Yoccoz's Function. Conclusion of the Proof}

\vskip .3 truecm
\noindent
\Proc{Proposition 3.9}{Let $\lambda_0\in\T$ and assume that 
$\lambda_0$ is not a root of unity. Then 
$r_2(\lambda_0 )\ge \limsup_{\D\ni\lambda\rightarrow
\lambda_0}|u(\lambda )|$.}

\vskip .3 truecm
\noindent
\proof
Let $r= \limsup_{\D\ni\lambda\rightarrow
\lambda_0}|u(\lambda )|$. It is not restrictive to assume 
$r>0$. Let $(\lambda_n )_{n\ge 1}\subset\D$ such that $\lambda_n
\rightarrow\lambda_0$ and $|u(\lambda_n)|\rightarrow r$. 
Since the linearizations $H_{\lambda_n}$ are univalent on 
their disks of convergence $\D_{r_2(\lambda_n )}$
one can extract a subsequence uniformly convergent on the 
compact susbets of $\D_r$. The limit function $H$ verifies 
$H(0)=0$, $H'(0)=1$ and $H(\lambda_0z)=P_{\lambda_0}(H(z))$
(this is immediate by taking the limit of the corresponding 
equations for $\lambda_n$). Thus $H_{\lambda_0}=H$
and $r_2(\lambda_0)\ge r$. \qed

\vskip .3 truecm
\noindent
Yoccoz has indeed proved the following stronger 
result [Yo2, pp. 65-69]

\vskip .3 truecm
\noindent
\Proc{Theorem 3.10}{For all $\lambda_0\in\T$, $|u(\lambda )|$
has a non--tangential limit in $\lambda_0$ which is equal to the 
radius of convergence $r_2(\lambda_0)$ of $H_{\lambda_0}$.}

\vskip .3 truecm
\noindent
Of course, if $\lambda_0$ is a root of unity then $P_{\lambda_0}$
is not even formally linearizable and one poses $r_2(\lambda_0)=0$.

Collecting Propositions 3.3, 3.7 and 3.9 together one can finally 
prove Theorem 3.2.

\vskip .3 truecm
\noindent
\proof of Theorem 3.2. Applying Fatou's Theorem to $u\,:\D
\rightarrow\C$ one finds that there exists $u^*\in L^\infty 
(\T,\C)$ such that for almost all $\lambda_0\in\T$ one has 
$|u^*(\lambda_0)|>0$ and $u(\lambda )\rightarrow u(\lambda_0)$
as $\lambda\rightarrow\lambda_0$ non tangentially. From 
Proposition 3.7 one concludes that for almost all $\lambda_0\in
\T$ one has $r_2(\lambda_0)>0$. \qed

\vskip .3 truecm
\noindent
{\bf Remark 3.11} Continuing the above argument of Yoccoz, L. 
Carleson and P. Jones prove that for almost all $\lambda\in\T$
the critical point $z=1$ of $P_{\lambda}$ belongs to the boundary of 
the Siegel disk (see, for example, [CG]). 
This has also been proved directly by M. Herman
under the assumption that $\alpha$ is diophantine [He3]. 
M. Herman has also shown that there are $\lambda$'s for which 
the critical point is not on the boundary of the Siegel disk 
(even though the boundary is a quasicircle) [Do]. 
%%%%% fine capitolo 3 %%%%%
\vfill\eject
%\vskip 1. truecm
%%%%% capitolo 4: Douady-Ghys %%%%%
\noindent
{\tit 4. Douady--Ghys' Theorem. Continued Fractions and the Brjuno 
Function} 

\vskip .3 truecm
\noindent
From Yoccoz's theorem it follows that $G_\lambda$, $\lambda =
e^{2\pi i \alpha}$, $\alpha\in\R\setminus\Q$, is a conjugacy
class for almost all values of $\alpha$. Let ${\cal Y}$ denote the
set of $\alpha\in\R\setminus\Q$ such that $G_{e^{2\pi i \alpha}}$ is
a conjugacy class. Then we already know that ${\cal Y}$
has full measure but that its complement in $\R\setminus\Q$
is a $G_\delta$--dense (Exercise 1.18). The goal of this Section 
is to prove a result due to Douady and Ghys on the structure of 
${\cal Y}$ (Section 4.1) and to introduce various sets of irrational numbers
(Sections 4.3 and 4.4)  
which have the same properties of ${\cal Y}$. Our main tool 
will be the use of continued fractions
(see Appendix 2 for a short introduction).

\line{}

\vskip .3 truecm
\noindent
{\titsec 4.1 Douady--Ghys' Theorem }

\vskip .3 truecm
\noindent
We recall that $\hbox{SL}\,(2,\Z )$ is the group of 
matrices $g=\left(\matrix{ a & b\cr c & d\cr}\right)$
with integer coefficients $a,b,c,d$ such that $ad-bc=1$. 
It acts on $\R\cup\{\infty\}$ (thus on ${\cal Y}$ too) 
as usual: $g\cdot\alpha = 
{a\alpha +b\over c\alpha +d}$. $\hbox{SL}\,(2,\Z )$
is generated by $T=\left(\matrix{ 1 & 1 \cr 0 & 1\cr}\right)$,
$T\cdot\alpha =\alpha +1$, and $U=
\left(\matrix{ 0 & -1 \cr 1 & 0\cr}\right)$, 
$U\cdot\alpha = -1/\alpha$. 

Further information on the structure of 
${\cal Y}$ is provided by the following 

\vskip .3 truecm
\noindent
\Proc{Theorem 4.1 (Douady--Ghys)}{${\cal Y}$ is 
$\hbox{SL}\,(2,\Z )$--invariant.}

\vskip .3 truecm
\noindent
\proof (sketch). ${\cal Y}$ is clearly invariant 
under $T$, thus we only need to show that if $\alpha\in {\cal Y}$
then also $U\cdot\alpha = -1/\alpha\in {\cal Y}$. 

Let $f\in S_{e^{2\pi i \alpha}}$ and consider a domain $V'$ bounded
by 
\item{1)} a segment $l$ joining $0$ to $z_0\in\D^*$, $l\subset\D$;
\item{2)} its image $f(l)$;
\item{3)} a curve $l'$ joining $z_0$ to $f(z_0)$.

\noindent
We choose $l'$ and $z_0$ (sufficiently close to $0$) so that $l$, 
$l'$ and $f(l)$ do not intersect except at their extremities.
Note that $l$ and $f(l)$ form an angle of $2\pi\alpha$ at $0$. Then 
glueing $l$ to $f(l)$ one obtains a topological manifold $\overline{V}$
with boundary which is homeomorphic to $\overline{\D}$. With the 
induced complex structure its interior is biholomorphic to 
$\D$. Let us now consider the {\it first return map} $g_{V'}$
to the domain $V'$ (this is well defined if $z$ is choosen with 
$|z|$ small enough): if $z\in V'$ (and $|z|$ is small enough) we define
$g_{V'}(z)=f^n(z)$ where $n$ (depends on $z$) is defined asking that 
$f(z),\ldots ,f^{n-1}(z)\not\in V'$ and $f^n(z)\in V'$, i.e.
$n=\inf\{k\in\N\, , \, k\ge 1\, , \, f^k(z)\in V'\}$. Then it is 
easy to check that $n=\left[{1\over\alpha}\right]$ or 
$n=\left[{1\over\alpha}\right]+1$. 
The first return map $g_{V'}$ induces a map $g_{\overline{V}}$
on a neighborhood of $0\in\overline{V}$ and finally a germ 
$g$ of holomorphic diffeomorphism at $0\in\D$ ($gp=pg_{\overline{V}}$,
where $p$ is the projection from $\overline{V}$ to the disk 
$\D$). It is easy to check that $g(z)=e^{-2\pi i/\alpha}z+\hbox{O}(z^2)$
(note that in the passage from $V'$ to $\D$ through $\overline{V}$
the angle $2\pi\alpha$ at the origin is mapped in $2\pi$). 

To each orbit of $f$ near $0$ corresponds an orbit of $g$ near $0$. In 
particular 
\item{$\bullet$} $f$ is linearizable if and only if $g$ is linearizable;
\item{$\bullet$} if $f$ has a periodic orbit near $0$ then also $g$ 
has a periodic orbit;
\item{$\bullet$} if $f$ has a point of instability (i.e. a point which 
does not belong to $K_f$) then also $g$ has a point of instability
(which, after having normalized $g$ so as to be univalent on $\D$, will
leave the unit disk even more rapidly). 

In particular these statements show that $\alpha\in {\cal Y}$
if and only if $-1/\alpha
\in {\cal Y}$. \qed

\line{}

\vskip .3 truecm
\noindent
{\titsec 4.2 SL(2,Z) and Continued Fractions}

\vskip .3 truecm
\noindent
To better understand  the action of $\hbox{SL}\,(2,\Z )$ on $\R\setminus\Q$
we can introduce a fundamental domain $[0,1)$ for one of the two generators
(the translation $T$) and restrict our attention to the inversion 
$\alpha\mapsto 1/\alpha$ restricted to $[0,1)$. This gives us a 
``microscope'' since $\alpha\mapsto 1/\alpha$ is expanding on $[0,1)$, 
i.e. its derivative is always greater than $1$. Our microscope magnifies more 
and more as $\alpha\rightarrow 0+$ and leads to the introduction of 
continued fractions discussed in Appendix A2. 

\vskip .3 truecm
\noindent
{\bf Exercise 4.2} Show that given any pair $x,y\in\R\setminus\Q$
there exists $g\in \hbox{SL}\,(2,{\Bbb Z})$ such that $x=g\cdot y$
{\it if and only if}
$x=[a_0,a_1,\ldots ,a_m,c_0,c_1,\ldots ]$ and 
$y=[b_0,b_1,\ldots ,b_n,c_0,c_1,\ldots ]$.
[Hint: it is easy to check that the condition is sufficient
for having $x=g\cdot y$; necessity is more tricky, see 
[HW] pp. 141--143.]

\vskip .3 truecm
\noindent
{\bf Exercise 4.3} Show that if $x$ is a quadratic irrational, 
i.e. $x\in\R\setminus\Q$ is a zero of a monic quadratic 
polynomial equation with coefficients in $\Q$, then there exists $N\in\N$
such that the partial fractions $a_n$ of $x$ are bounded
$a_n\le N$ for all $n\ge 0$. 

\vskip .3 truecm
\noindent
The two main results which make continued fractions so useful 
in the study of one--dimensional small divisors problems 
are the following

\vskip .3 truecm
\noindent
\Proc{Theorem 4.4 (Best approximation)}{Let $x\in\R\setminus\Q$
and let $p_n/q_n$ denote its $n$--th convergent. If $0<q<q_{n+1}$
then $|qx-p|\ge |q_nx-p_n|$ for all $p\in\Z$ and equality 
can occur only if $q=q_n$, $p=p_n$.}

\vskip .3 truecm
\noindent
\Proc{Theorem 4.5}{If $\left|x-{p\over q}\right|<{1\over 2q^2}$
then ${p\over q}$ is a convergent of $x$.}

\vskip .3 truecm
\noindent
For the proofs see [HW], respectively Theorems 182, p. 151 and 184, p. 153. 

\line{}

\vskip .3 truecm
\noindent
{\titsec 4.3 Classical Diophantine Conditions}

\vskip .3 truecm
\noindent
Let $\gamma >0$ and $\tau\ge 0$ be two real numbers.

\vskip .3 truecm
\noindent
\Proc{Definition 4.6}{
$x\in\R\setminus\Q$ is {\rm diophantine} of exponent $\tau$
and constant $\gamma$ if and only if for all $p,q\in\Z$, 
$q> 0$, one has $\left|x-{p\over q}\right|\ge\gamma q^{-2-\tau}$.}

\vskip .3 truecm
\noindent
{\bf Remark 4.7} Note that Theorem 4.5 implies that 
given any irrational number there are infinitely many 
solutions to $\left|x-{p\over q}\right|<{1\over q^2}$
with $p$ and $q$ coprime. This explains why the previous 
definition would never be satisfied if $\tau <0$. 

\vskip .3 truecm
\noindent
We denote $\CD{\gamma}{\tau}$ the set of all irrationals 
$x$ such that $\left|x-{p\over q}\right|\ge\gamma q^{-2-\tau}$
for all $p,q\in\Z$, $q>0$. $\cd{\tau}$ will denote 
the union $\cup_{\gamma >0}\CD{\gamma}{\tau}$ and 
$\hbox{CD} \,=\cup_{\tau\ge 0}\cd{\tau}$. 

\vskip .3 truecm
\noindent
{\bf Exercise 4.8} Show that 
$$
\eqalign{
\cd{\tau} &=\{x\in\R\setminus\Q\mid\; q_{n+1}=\hbox{O}(q_n^{1+\tau})\}
=\{x\in\R\setminus\Q\mid\; a_{n+1}=\hbox{O}(q_n^{\tau})\}\cr
&= \{x\in\R\setminus\Q\mid\; x_{n}^{-1}=\hbox{O}(\beta_{n-1}^{-\tau})\}
= \{x\in\R\setminus\Q\mid\; \beta_{n}^{-1}=\hbox{O}(\beta_{n-1}^{-1-\tau})\}
\cr
}
$$

\vskip .3 truecm
\noindent
{\bf Exercise 4.9} Show that if $x$ is an algebraic number of degree 
$n\ge 2$, i.e. $x\in\R\setminus\Q$ is a zero of a monic 
polynomial with coefficients in $\Q$ and degree $n$, then 
$x\in\cd{n-2}$ (Liouville's theorem). Thue improved this result
in 1909 showing that  $x\in\cd{\tau -1+n/2}$ for all $\tau >0$ 
(see [ST], Chapter V, for a very nice discussion of the proof
in the cubic case). 
Actually one can prove that if $x$ is algebraic 
then $x\in\cd{\tau}$ for all $\tau >0$ regardless of the degree, but 
this is difficult (Roth's theorem). 

\vskip .3 truecm
\noindent
{\bf Exercise 4.10} Using the fact that the continued fraction of
$e=\sum_{n=0}^\infty {1\over n!}$ is 
$$
[2,1,2,1,1,4,1,1,6,1,1,8,1,1,10,\ldots ]
$$
show that $e\in\cap_{\tau >0}\cd{\tau }$. A proof of the continued 
fraction expansion of $e$, which is due to L. Euler, can be found in 
[L], Chapter V. Perhaps you may like to try to obtain it yourself 
starting from the knowledge of the continued 
fraction of ${e+1\over e-1}=[2,6,10,14,\ldots ]$. 

\vskip .3 truecm
\noindent
{\bf Exercise 4.11} Use the result of Exercise 4.9 to exhibit 
explicit examples of trascendental numbers, e.g. $x=\sum_{n=0}^\infty
10^{-n!}$. 

\vskip .3 truecm
\noindent
The complement in $\R\setminus\Q$ of $\hbox{CD}$ is called the set 
of Liouville numbers. 

\vskip .3 truecm
\noindent
{\bf Exercise 4.12} Show that $\cd{\tau}$ and $\hbox{CD}$ are both 
$\hbox{SL}\,(2,\Z)$--invariant.

\vskip .3 truecm
\noindent
\Proc{Proposition 4.13}{For all $\gamma >0$ and $\tau >0$ the Lebesgue measure 
of $\CD{\gamma}{\tau} (\hbox{mod}\, 1)$ is at least 
$1-2\gamma\zeta (1+\tau)$, where $\zeta$ denotes the Riemann zeta function.}

\vskip .3 truecm
\noindent
\proof
The complement of $\CD{\gamma}{\tau} (\hbox{mod}\, 1)$ is contained in 
$$
\cup_{p/q\in\Q\cap [0,1]}\left(
{p\over q}-\gamma q^{-2-\tau},{p\over q}+\gamma q^{-2-\tau}
\right)
$$
whose Lebesgue measure is bounded by 
$$
\sum_{q=1}^\infty \sum_{p=1}^q 2\gamma q^{-2-\tau}\le 2\gamma\sum_{q=1}^\infty
q^{-1-\tau}\; . 
$$
\qed

\vskip .3 truecm
\noindent
From the point of view of dimension one has (see [Fa], p. 142 for a proof)

\vskip .3 truecm
\noindent
\Proc{Theorem 4.14 (Jarnik)}{Let $\tau >0$ and let $F_\tau$ be the set of real 
numbers $x\in [0,1]$ such that $\{qx\}\le q^{-1-\tau}$ for infinitely many 
positive integers $q$. The Hausdorff dimension of $F_\tau$ is $2/(2+\tau )$.}

\vskip .3 truecm
\noindent
{\bf Exercise 4.15} The set $\cd{0}$ is also called the set of numbers of 
{\it constant type} since $x\in\cd{0}$ if and only if the sequence of its
partial fractions is bounded. Show that $\cd{0}$ has Hausdorff dimension $1$
and zero Lebesgue measure.

\vskip .3 truecm
\noindent
{\bf Exercise 4.16} Show that the set of Liouville numbers has zero Lebesgue 
measure, zero Hausdorff dimension but it is a dense $G_\delta$--set
%(i.e. it is a countable intersection of open sets: dense $G_\delta$--sets
%are ``big'' from a topologist's point of view.

\line{}

\vskip .3 truecm
\noindent
{\titsec 4.4 Brjuno Numbers and the Brjuno Function} 

\vskip .3 truecm
\noindent
Let $x\in\R\setminus\Q$, let $\left({p_n\over q_n}\right)_{n\ge 0}$
denote the sequence of its convergents and let $(\beta_n)_{n\ge -1}$
be defined as in (A2.14).

\vskip .3 truecm
\noindent
\Proc{Definition 4.17}{$x$ is a {\rm Brjuno number} if 
$B(x):=\sum_{n=0}^\infty\beta_{n-1}\log x_n^{-1}<+\infty$. 
The function $B\,: \R\setminus\Q\rightarrow (0,+\infty ]$
is called the {\rm Brjuno function}.}

\vskip .3 truecm
\noindent
{\bf Exercise 4.18} Show that all diophantine numbers are Brjuno numbers.

\vskip .3 truecm
\noindent
{\bf Exercise 4.19} Show that there exists $C>0$ such taht for all 
Brjuno numbers $x$ one has 
$$
\left| B(x)-\sum_{n=0}^\infty {\log q_{n+1}\over q_n}\right|\le C\; .
$$

\vskip .3 truecm
\noindent
{\bf Exercise 4.20 (see [MMY])} Show that the Brjuno function satisfies 
$$
\eqalign{
B(x) &= B(x+1)\; , \;\;\;\forall  x\in\R\setminus\Q\cr
B(x) &= -\log x+xB\left({1\over x}\right)\; ,\;\;\;x\in\R\setminus\Q
\cap (0,1)\cr}
\eqno(4.1)
$$
Deduce from this that the set of Brjuno numbers is $\hbox{SL}\,(2,\Z)$--invariant.
Use the above given functional equation to compute $B(x_p)$, where 
$x_p={\sqrt{p^2+4}-p\over 2}$, $p\in\N$. 

\vskip .3 truecm
\noindent
{\bf Exercise 4.21 (see [MMY])}
Show that the linear operator 
$(Tf)(x)=xf\left({1\over x}\right)$, $x\in (0,1)$, 
 acting on periodic functions which 
belong to $L^p\left(\T,{dx\over (1+x)\log 2}\right)$
has spectral radius bounded by $g={\sqrt{5}-1\over 2}$. 
Conclude that the Brjuno function $B\in\cap_{p\ge 1}L^p(\T )$. 
Note that $B\not\in L^\infty (\T )$. 

\vskip .3 truecm
\noindent
{\bf Exercise 4.22} Write the continued fraction expansion of a Brjuno number 
which is not a diophantine number. The same for the decimal expansion.
Is $\sum_{n=0}^\infty 10^{-n!}$ a Brjuno number? What about 
$\sum_{n=0}^\infty 10^{-10^{n!}}$? 

\vskip .3 truecm
\noindent
{\bf Exercise 4.23} Let $\sigma >0$. Use the results of Exercise 4.21
to study the functions 
$$
\eqalign{
B^{(\sigma )} (x) &= B^{(\sigma )}(x+1)\; , \;\;\;\forall  x\in\R\setminus\Q\cr
B^{(\sigma )} (x) &= x^{-1/\sigma}+
xB^{(\sigma )}\left({1\over x}\right)\; ,\;\;\;x\in\R\setminus\Q
\cap (0,1)\cr}
$$
Show that if $B^{(\sigma )}(x)<+\infty$ then 
$x\in \cd{\sigma}$. Viceversa, if $x\in\cd{\tau}$ then 
$B^{(\sigma )}(x)<+\infty$ for all $\sigma >\tau$. 
%%%%% fine capitolo 4 %%%%%
\vfill\eject
%\vskip 1. truecm
%%%%% capitolo 5: teoremi di Brjuno, Siegel e Yoccoz %%%%%
\noindent
{\tit 5. Siegel--Brjuno Theorem, Yoccoz's Theorem and Some
Open Problems}

\vskip .3 truecm
\noindent

Recall that  ${\cal Y}$ denotes the
set of $\alpha\in\R\setminus\Q$ such that $G_{e^{2\pi i \alpha}}$ is
a conjugacy class. Here we list what we already know about it
\item{$\bullet$} ${\cal Y}$
has full measure (Chapter 3); 
\item{$\bullet$} the complement of ${\cal Y}$ in $\R\setminus\Q$
is a $G_\delta$--dense (Exercise 1.18);
\item{$\bullet$} ${\cal Y}$ is invariant under the action of 
$\hbox{SL}\,(2,\Z )$ (Douady--Ghys' Theorem, Chapter 4). 

The purpose of this Chapter is to prove the classical results 
of Siegel [S] and Brjuno [Br] which show that:

{\it all 
Brjuno numbers belong to ${\cal Y}$}

Moreover we will state the Theorems of Yoccoz [Yo2] and in particular 
his celebrated result:

{\it ${\cal Y}$ is equal to the set of  Brjuno numbers.}

We will also mention several open problems. 

For the sake of brevity, starting with this section  
all the proofs will be only sketched: the reader can try autonomously
to fill in the details but we will always refer to the original 
literature where complete proofs are given. 

\line{}

\vskip .3 truecm\noindent
{\titsec 5.1 Siegel--Brjuno Theorem}

\vskip .3 truecm\noindent
The theorem of Siegel and Brjuno says that the set of Brjuno numbers 
is a subset of ${\cal Y}$. Indeed in 1942 C.L. Siegel was the first 
to show that ${\cal Y}$ is not empty showing that $\hbox{CD}\,\subset
{\cal Y}$. 

We will sketch the proof of a more precise result which follows from 
the Theorem of Yoccoz which we will discuss in the next section but 
which can also be proved following the classical majorant series 
method (see [CM]). 
Let us recall that $S_{\lambda}$ denotes the topological space of all 
germs of holomorphic diffeomorphisms $f\, :\D\rightarrow\C$ such that 
$f(0)=0$, $f$ is univalent on $\D$ and $f'(0)=\lambda$. By Theorem 
A1.19 it is a compact space. Given a germ $f\in S_{\lambda}$ 
we let $r(f)$ indicate the radius of convergence of the linearization 
$h_{f}$ of $f$ and we set 
$$
r(\alpha )= \inf_{f\in S_{e^{2\pi i\alpha}}}r(f)\; . \eqno(5.1)
$$

\vskip .3 truecm\noindent
\Proc{Theorem 5.1 (Yoccoz's lower bound)}{
$$
\log r(\alpha )\ge -B(\alpha )-C\eqno(5.2)
$$
where $C>0$ is a universal constant (independent of $\alpha$)
and $B$ is the Brjuno function.}

\vskip .3 truecm\noindent
Before of sketching the proof let us briefly mention what is the main 
difficulty which was first overcome by Siegel and which was clearly 
well--known among mathematicians at the end of the $19$th and at the 
beginning of the $20$th century (in 1919 Gaston Julia even claimed, in 
an incorrect paper, to disprove Siegel's theorem). 

Assume that $\alpha\in\cd{\tau}$ for some $\tau\ge 0$. Recalling the 
recurrence (1.2) for the power series coefficients of the 
linearization $h_{f}(z)=\sum_{n=1}^\infty h_{n}z^n$
$$
h_{1}=1\; , \;\;\;
h_n={1\over \lambda^n-\lambda}\sum_{j=2}^nf_j
\sum_{n_1+\ldots +n_j=n}h_{n_1}\cdots h_{n_j}\; \; n\ge 2\; \; 
 \eqno(5.3)
$$
one sees that $h_{n}$ is a polynomial in $f_{2},\ldots ,f_{n}$ with
coefficients which are rational functions of $\lambda$: 
$h_{n}\in\C (\lambda )[f_{2},\ldots ,f_{n}]$ for all $n\ge 2$.

Let us compute explicitely the first few terms of the recurrence 
$$
\eqalign{
h_{2} &= (\lambda^{2}-\lambda )^{-1}f_{2}\; , \cr
h_{3} &= (\lambda^{3}-\lambda )^{-1}[f_{3}+2f_{2}^{2}(\lambda^{2}
-\lambda )^{-1}]\; , \cr
h_{4} &= (\lambda^{4}-\lambda )^{-1}[f_{4}+3f_{3}f_{2}(\lambda^{2}-
\lambda )^{-1}+2f_{2}f_{3}(\lambda^{3}-\lambda )^{-1}\cr
&\phantom{= (\lambda^{4}-\lambda )^{-1}[} 4f_{2}^{3}
(\lambda^{3}-\lambda )^{-1}(\lambda^{2}-\lambda )^{-1}
+f_{2}^{3}(\lambda^{2}-\lambda )^{-2}]\; , \cr
}
\eqno(5.4)
$$
and so on. It is not difficult to see that among all contributes to 
$h_{n}$ there is always a term of the form $2^{n-2}f_{2}^{n-1}
[(\lambda^{n}-\lambda )\ldots (\lambda^{3}-\lambda )
(\lambda^{2}-\lambda )]^{-1}$. If one then tries to estimate 
$|h_{n}|$ by simply summing up the absolute values of each 
contribution  then one term will be
$$
2^{n-2}|f_{2}|^{n-1}
[|\lambda^{n}-\lambda |\ldots |\lambda^{3}-\lambda |
|\lambda^{2}-\lambda |]^{-1}\le 
2^{n-2}|f_{2}|^{n-1}(2\gamma )^{(n-1)\tau}[(n-1)!]^\tau
\eqno(5.5)
$$
if $\alpha\in\CD{\gamma}{\tau}$ and one obtains a divergent bound. 
Note the difference with the case $|\lambda |\not= 1$: in this case 
the bound would be $|\lambda |^{-(n-1)}2^{n-2}|f_{2}|^{n-1}c^{n-1}$
for some positive constant $c$ independent of $f$. Thus one must 
use a more subtle majorant series method. 

The key point is that the estimate (5.5) is far too pessimistic: 

\vskip .3 truecm\noindent
{\bf Exercise 5.2} Show that the series $\sum_{n=1}^\infty {z^n
\over (\lambda^n-1)\ldots (\lambda -1)}$, with $\lambda =
e^{2\pi i\alpha}$, has positive radius of convergence whenever 
$\limsup_{n\rightarrow\infty}{\log q_{k+1}\over q_{k}}<+\infty$. 
[Hint: see [HL] for a proof.]

\vskip .3 truecm\noindent
Indeed when a small divisor is really small then, for a certain time, 
all other small divisors cannot be too small. This vague idea is made 
clear by the two following lemmas of A.M. Davie [Da] which extend and 
improve previous results of A.D. Brjuno. 

Let $x\in\R$, $x\not= 1/2$, we denote 
$\Vert x\Vert_{\Z}=\min_{p\in\Z}|x+p|$. 

\vskip .3 truecm\noindent
\Proc{Lemma 5.3}{Let $\alpha\in\R\setminus\Q$, $(p_{j}/q_{j})_{j\ge 
0}$ denote the sequence of its convergents, 
$k\in\N$, $n\in\N$, $n\not= 0$, and assume 
that $\Vert n\alpha \Vert_{\Z}\le 1/(4q_{k})$. 
Then $n\ge q_{k}$ and either $q_{k}$ divides $n$ or $n\ge 
q_{k+1}/4$.}

\vskip .3 truecm\noindent
\proof
From Theorem 4.5 it follows that if $r$ is an integer and 
$0<r<q_{k}$ then $\Vert r\alpha\Vert_{\Z}\ge (2q_{k})^{-1}$. Thus 
$n\ge q_{k}$. Assume that $q_{k}$ does not divide $n$ and that 
$n<q_{k+1}/4$. Then $n=mq_{k}+r$ where $0<r<q_{k}$ and $m<q_{k+1}/
(4q_{k})$. Since $\Vert q_{k}\alpha\Vert_{\Z}\le q_{k+1}^{-1}$
one gets $\Vert mq_{k}\alpha\Vert_{\Z}\le 
mq_{k+1}^{-1}<(4q_{k})^{-1}$.  But $\Vert r\alpha\Vert_{\Z}\ge 
(2q_{k})^{-1}$ thus $\Vert n\alpha\Vert_{\Z}>(4q_{k})^{-1}$. \qed

\vskip .3 truecm\noindent
Using this information on the sequence $(\Vert 
n\alpha\Vert_{\Z})_{n\ge 0}$ Davie shows the following: 
Let $A_k = \left\{ n \ge 0 \mid \| n \alpha \| \le {1\over 8q_k} 
\right\}$, $E_k=\max \left( q_k , q_{k+1}/4 \right)$ and $\eta_k = 
q_k / E_k$.
Let $A_k^{*}$ be the set of non negative integers 
$j$ such that either $j \in A_k$ or for some $j_1$ and $j_2$ in 
$A_k$, with $j_2-j_1 < E_k$, one has $j_1 < j < j_2$ and $q_k$ 
divides $j-j_1$.
For any non negative integer $n$ define:
$$
l \left( n \right) = \max \left\{ \left( 1+\eta_k \right) 
{n\over q_k}-2 , \left( m_n \eta_k+n \right) {1\over q_k}-1 \right\}
\eqno(5.6)
$$
where $m_n = \max \{ j \mid 0 \le j \le n , j \in A_k^{*} \}$.
We then define a function $h_k \, : \N\rightarrow\R_{+}$
as follows
$$
h_k \left( n \right)= \cases{
{m_n+\eta_k n\over q_k}-1& if $m_n+q_k \in A_k^{*}$ \cr
l \left( n \right)& if $m_n+q_k \not\in A_k^{*}$\cr}
\eqno(5.7)
$$

The function $h_k \left( n \right)$ has some properties collected in 
the following proposition

\vskip .3 truecm\noindent
\Proc{Proposition 5.4}{
The function $h_k \left( n \right)$ verifies
\item{(1)} ${\left( 1+\eta_k \right)n\over q_k}-2 \le h_k \left( n 
\right) \le {\left( 1+\eta_k \right)n \over q_k}-1$ for all $n$.
\item{(2)} If $n>0$ and $n \in A_k^{*}$ then $h_k \left( n \right) \ge
h_k \left( n -1 \right)+1$.
\item{(3)} $h_k \left( n \right) \ge h_k \left( n-1 \right)$ for all 
$n>0$.
\item{(4)} $h_k \left( n+q_k \right) \ge h_k \left( n \right) +1$ for all 
$n$.}

\vskip .3 truecm
\noindent
Now we set $g_k \left( n \right)= \max \left( h_k \left( n \right), 
\left[ {n\over q_k} \right] \right)$ and we state the following 
proposition

\vskip .3 truecm\noindent
\Proc{Proposition 5.5}{
The function $g_k $ 
is non negative and verifies:
\item{(1)} $g_k(0)=0$;
\item{(2)} $g_k \left( n \right) \le {\left( 1+\eta_k 
\right)n\over q_k}$ for all $n$;
\item{(3)} $g_k \left( n_1 \right) + g_k \left( n_2 \right) \le g_k 
\left( n_1 +n_2 \right)$ for all $n_1$ and $n_2$;
\item{(4)} if $n \in A_k$ and $n>0$ then $g_k \left( n \right) \ge g_k 
\left( n-1 \right)+1$.}

\vskip .3 truecm\noindent
The proof of these propositions can be found in [Da]. 

Let $k(n)$ be defined by the condition $q_{k(n)}\le n <q_{k(n)+1}$. 
Note that $k$ is non--decreasing. 

\vskip .3 truecm\noindent
\Proc{ Lemma 5.6 (Davie's lemma)}{Let 
$$
K(n)=n\log 2+\sum_{k=0}^{k(n)}g_{k}(n)\log(2q_{k+1})\; . \eqno(5.8)
$$
The function $K \left( n \right)$ verifies:
\item{(a)} There exists a universal constant $c_{0}>0$ such that 
$$
K(n)\le n\left(\sum_{k=0}^{k(n)}{\log 
q_{k+1}\over q_{k}}+c_{0}\right)\; ;\eqno(5.9) 
$$
\item{(b)} $K(n_{1})+K(n_{2})\le K(n_{1}+n_{2})$ for all $n_{1}$ and $n_{2}$; 
\item{(c)} $-\log |\lambda^n -1| \le K(n)-K(n-1)$. }

\vskip .3 truecm\noindent
\proof
We will apply Proposition 5.4. By (2) we have 
$$
\eqalign{
K(n) &\le n\left[\log 2+\sum_{k=0}^{k(n)}
{(1+\eta_{k})\over q_{k}}\log (2q_{k+1})\right]\cr
&\le n\left[\sum_{k=0}^{k(n)}{\log q_{k+1}\over q_{k}}+\log 2
+\log 2\sum_{k=0}^\infty \left({1\over q_{k}}+{4\over q_{k+1}}\right)
+4\sum_{k=0}^\infty {\log q_{k+1}\over q_{k+1}}\right]\cr}
$$
since $\eta_{k}\le 4q_{k}q_{k+1}^{-1}$. 

By Remark A2.4 the series 
$ \sum {\log q_{k+1}\over q_{k+1}}$ and $\sum q_{k}^{-1}$ are 
uniformly bounded by some constant independent of $\alpha$, thus (a)
follows. 

(b) is an immediate consequence of Proposition 5.5, (3), and the fact 
that $k(n)$ is not decreasing. 

Finally recall that 
$$
-\log |\lambda^n-1| = -\log 2|\sin \pi n\alpha |\in
(-\log\pi\Vert n\alpha\Vert_{\Z},-\log 2\Vert n\alpha\Vert_{\Z})
\; . 
$$
For all $n$ we have either $\Vert n\alpha\Vert_{\Z}>1/4$ or there 
exists some non--negative integer $k$ such that $(8q_{k})^{-1}>
\Vert n\alpha\Vert_{\Z}\ge (8q_{k+1})^{-1}$, thus $n\in A_{k}$, 
which implies by Proposition 5.5, (4), that $g_{k}(n)\ge 
g_{k}(n-1)+1$, and $-\log |\lambda^n-1|\le 2q_{k+1}$. 
Combining these facts together one gets (c). 
\qed

\vskip .3 truecm\noindent
Following [CM] we can now prove Theorem 5.1.

\proof {\it of Theorem 5.1. }
Let $s\left( z\right) =\sum_{n\geq 1}s_nz^n$ be  the 
unique  solution analytic at $z=0$ of the equation 
$s\left( z\right) =z+\sigma \left( s\left( z\right) \right)$, where 
$\sigma (z) = {z^{2}(2-z)\over (1-z)^{2}}= \sum_{n\ge 2}nz^n$. 
The coefficients satisfy 
$$
s_{1}= 1 \; , \; s_{n} = 
\sum_{m=2}^n m\sum_{n_{1}+\ldots +n_{m}= n\, , \, 
n_{i}\ge 1 } s_{n_{1}}\ldots s_{n_{m}}\; . \eqno(5.10)
$$
Clearly there exist two positive constants $c_{1},c_{2}$
such that 
$$
|s_n| \le c_{1}c_{2}^{n}\; . 
$$
From the recurrence relation and 
Bieberbach--De Branges's bound
$|f_{n}|\le n $ for all $n\ge 2$ 
we obtain 
$$
|h_{n}| \le {1\over |\lambda^n -\lambda|}
\sum_{m=2}^n m\sum_{n_{1}+\ldots +n_{m}= n\, , \, 
n_{i}\ge 1 } |h_{n_{1}}|\ldots |h_{n_{m}}|\; . 
$$
We now deduce by induction on $n$ that $|h_{n}|\le s_{n}e^{K(n-1)}$
for $n\ge 1$, where $K\, : \N\rightarrow \R$ 
is defined in (5.8). If we assume this 
holds for all $n'<n$ then the above inequality gives 
$$
|h_{n}| \le {1\over |\lambda^n -\lambda|}
\sum_{m=2}^n m\sum_{n_{1}+\ldots +n_{m}= n\, , \, 
n_{i}\ge 1 } s_{n_{1}}\ldots s_{n_{m}}e^{K(n_{1}-1)+\ldots K(n_{m}-1)}
\; .
$$
But 
$K(n_{1}-1)+\ldots K(n_{m}-1)\le K(n-2)\le K(n-1)+\log |\lambda^n -\lambda|$
and we deduce that 
$$
	|h_{n}|\le e^{K(n-1)}\sum_{m=2}^n m\sum_{n_{1}+\ldots +n_{m}= n\, , \, 
n_{i}\ge 1 } s_{n_{1}}\ldots s_{n_{m}}=s_{n} e^{K(n-1)}\; , 
$$
as required. Theorem 5.1 then follows from the fact that 
$n^{-1}K(n)\le B(\omega )+c_{0}$ for some universal constant 
$c_{0}>0$ (Davie's lemma). 
\qed

\vskip .3 truecm\noindent
{\bf Exercise 5.7} Consider the quadratic polynomial 
$P_{\lambda}(z)=\lambda (z-z^{2})$
(we have conjugated (3.1) by an homothehty so as to eliminate 
a factor $1/2$ in what follows). Its formal linearization 
$H_{\lambda}(z)=\sum_{n=1}^\infty H_{n}(\lambda )z^n$ is given by the 
recurrence 
$$
H_{1}(\lambda )=1\; , \;\;
H_{n}(\lambda )=(1-\lambda^{n-1})^{-1}\sum_{i+j=n}
H_{i}(\lambda )H_{j}(\lambda )\; . 
$$
Define the sequence of positive real numbers 
$(h_{n}(\lambda ))_{n\ge 1}$ by the recurrence 
$$
h_{1}(\lambda )=1\; , \;\; 
h_{n}(\lambda )=|1-\lambda^{n-1}|^{-1}\sum_{i+j=n}
h_{i}(\lambda )h_{j}(\lambda )\; . 
$$
Clearly $|H_{n}(\lambda )|\le h_{n}(\lambda )$. Show that if $\lambda 
=e^{2\pi i \alpha}$, $\alpha\in\R\setminus\Q$ is {\it not} a Brjuno 
number then $\limsup_{n\rightarrow\infty}n^{-1}\log h_{n}(\lambda )
=+\infty $, i.e. the majorant series $\sum_{n\ge 1}h_{n}(\lambda )z^n$
is divergent. Under what assumptions on $\alpha$ one can show that 
there exist two positive constants $c_{0},c_{1}$ and $s>0$ 
such that $h_{n}(\lambda )\le c_{0}c_{1}^n(n!)^s$ for all $n\ge 1$?
[Hint: First show that the set $\{n\ge 0\mid q_{n+1}\ge 
(q_{n}+1)^{2}\}$ is infinite and denote $(q_{i}')_{i\ge 0}$ its 
elements. Then show that the sequence $h_{n}(\lambda )$ is increasing 
and $h_{q_{r+1}'+1}(\lambda )\ge |1-\lambda^{q_{r+1}'}|^{-1}
(h_{q_{r}'+1}(\lambda ))^{\left[{q_{r+1}'\over q_{r}'+1}\right]}$. 
One can also consult [Yo2, Appendice 2, pp. 83--85] and [CM].]

\vskip .3 truecm\noindent
{\bf Exercise 5.8} (Linearization and Gevrey classes, 
see [CM] for solutions and more information. )
Between $\C[[z]]$ and $\C\{z\}$ one has many 
important algebras of ``ultradifferentiable'' power series
(i.e. asymptotic expansions at $z=0$ of functions which are 
``between'' ${\cal C}^\infty$ and $\C\{z\}$). 
Consider two subalgebras $A_{1}\subset A_{2}$ of 
$z\C\left[\left[ z\right]\right]$
closed with respect to the composition of formal series. For example 
Gevrey--$s$ classes, $s>0$ (i.e. series $F(z)=\sum_{n\ge 0}f_{n}z^n$
such that there exist $c_{1}, c_{2}>0$ such that 
$|f_{n}|\le c_{1}c_{2}^n(n!)^s$ for all $n\ge 0$). Let 
$f\in A_{1}$ being such that $f^{\prime}\left( 0 \right)=\lambda\in 
\C^{*}$. We say that $f$ {\it is linearizable in } $A_{2}$ if there 
exists $h_{f}\in A_{2}$ tangent to the identity and such that 
$f \circ h_{f}  = h_{f} \circ R_\lambda$. 
Show that if one requires $A_{2}=A_{1}$, i.e.  the 
linearization $h_{f}$ to be as regular as the given germ $f$, once again the 
Brjuno condition is sufficient. It is 
quite interesting to notice that given any algebra of formal power 
series which is closed under composition (as it should if one whishes 
to study conjugacy problems)  a germ in the 
algebra is linearizable {\it in the same algebra} if the Brjuno 
condition is satisfied. 
If the linearization is allowed to be less regular than the given germ 
(i.e. $A_{1}$ is a proper subset of $A_{2}$) 
one finds new arithmetical conditions, weaker than the Brjuno 
condition. 
Let $(M_{n})_{n\ge 1}$ be a sequence of positive real numbers such 
that: 
\item{0.} $\inf_{n\ge 1} M_{n}^{1/n}>0$; 
\item{1.} There exists $C_{1}>0$ such that $M_{n+1}\le C_{1}^{n+1}M_{n}$ 
for all $n\ge 1$; 
\item{2.} The sequence $(M_{n})_{n\ge 1}$ is logarithmically convex; 
\item{3.} $M_{n}M_{m}\le M_{m+n-1}$ for all $m,n\ge 1$.

\noindent
Let $f= \sum_{n\ge 1}f_{n}z^n\in z\C
\left[ \left[ z \right] \right]$; $f$ belongs to 
the algebra 
$z\C\left[ \left[ z \right] \right]_{(M_{n})}$
if there exist two positive constants $c_{1},c_{2}$ 
such that 
$$
|f_{n}| \le c_{1}c_{2}^nM_{n}\;\; 
\hbox{for all}\; n\ge 1\; . 
$$

\noindent
Show that the condition 3 above implies that 
$z\C\left[ \left[ z \right] \right]_{(M_{n})}$ 
is closed for composition. 
Show that  if $f\in z\C\left[ \left[ z \right] \right]_{(N_{n})}$, 
$f_{1}=e^{2\pi i\alpha}$ 
and $\alpha$ verifies 
$$
	\limsup_{n\rightarrow +\infty} \left( \sum_{k=0}^{k(n)}
	{\log q_{k+1}\over q_{k}} - {1\over n}\log {M_{n}\over N_{n}}\right) <
	+\infty
$$
where $k(n)$ is defined by the condition $q_{k(n)}\le n< q_{k(n)+1}$, 
then the linearization $h_{f}\in z\C\left[ \left[ z \right] \right]_{(M_{n})}$.
(We of course assume that the sequence $(N_{n})_{n\ge 0}$  is asymptotically 
bounded by the sequence $(M_{n})$, i.e. $M_{n}\ge N_{n}$ for all 
sufficiently large $n$).  

\line{}

\vskip .3 truecm
\noindent
{\titsec 5.2 Yoccoz's Theorem }

\vskip .3 truecm\noindent
The main result of Yoccoz can be very simply stated as 
$$
{\cal Y}=\{\alpha\in\R\setminus\Q\mid\, B(\alpha )<+\infty\}=
\hbox{Brjuno numbers}\; , 
$$
but he proves much more than the above:

\vskip .3 truecm\noindent
\Proc{Theorem 5.9}{
\item{(a)} If $B(\alpha )=+\infty$ there exists a non--linearizable 
germ $f\in S_{e^{2\pi i\alpha}}$; 
\item{(b)} If $B(\alpha )<+\infty$ then $r(\alpha )>0$ and 
$$
|\log r(\alpha )+B(\alpha )|\le C\; , \eqno(5.11)
$$
where $C$ is a universal constant (i.e. independent of $\alpha$);
\item{(c)} For all $\varepsilon >0$ there exists $C_{\varepsilon}>0$
such that for all Brjuno numbers $\alpha$ one has 
$$
-B(\alpha )-C\le \log r(P_{e^{2\pi i\alpha}})\le -(1-\varepsilon ) 
B(\alpha )+C_{\varepsilon}\eqno(5.12)
$$
where $C$ is a universal constant (i.e. independent of $\alpha$
and $\varepsilon$.}

\vskip .3 truecm\noindent
The remarkable consequence of (5.11) and (5.12) is that 
the Brjuno function not only identifies the set ${\cal Y}$
but also gives a rather precise estimate of the size of the 
Siegel disks. When $\alpha$ is not a Brjuno number the problem of 
a complete classification of the conjugacy classes of germs in 
$G_{e^{2\pi i\alpha}}$ is open and quite difficult
(perhaps unreasonable) as the following result of Yoccoz shows:

\vskip .3 truecm\noindent
\Proc{Theorem 5.10}{Let $\alpha\in \R\setminus\Q$, 
$B(\alpha )=+\infty$. There exists a set with the power of the 
continuum of conjugacy classes of germs of $G_{e^{2\pi i\alpha}}$, 
each of which does not contain an entire function.}

\vskip .3 truecm\noindent
The proof of the Theorem of Yoccoz (5.9 above) uses a method, 
invented by Yoccoz himself, known as {\it geometric renormalization}. 
Roughly speaking it is a quantitative version of the topological 
construction of Douady--Ghys described in Chapter 4 and which shows 
that the set ${\cal Y}$ is $\hbox{SL}\, (2,\Z )$--invariant. 
Whereas the construction of  non--linearizable germs
$f\in G_{e^{2\pi i\alpha}}$ when $B(\alpha )=+\infty$ and Yoccoz's 
upper bound 
$$
\log r(\alpha )\le C-B(\alpha )
$$
go far beyond the scope of these lectures, it is not too difficult 
to give an idea of how Yoccoz proves Theorem 5.1, i.e. the 
lower bound 
$$
\log r(\alpha )\ge -C-B(\alpha )\; . 
$$
Let $f\in S_{e^{2\pi i\alpha}}$ and let $E\, :\H\rightarrow \D^{*}$
be the exponential map $E(z)=e^{2\pi i z}$. Then $f$ lifts to a map 
$F\, :\H\rightarrow \C$ such that 
$$
\eqalignno{
E\circ F &=f\circ E\; ,  &(5.13)\cr
F\, : \H &\rightarrow\C \;\hbox{is univalent}\; , &(5.14)\cr
F(z) &= z+\alpha+\varphi (z)\; \hbox{where}\, \varphi\;\hbox{is}\,
\Z--\hbox{periodic and}\,\lim_{\IM z\rightarrow +\infty}\varphi (z)=0\; . 
&(5.15)\cr
}
$$
We will denote $S(\alpha )$ the space of univalent functions $F$
verifying (5.13), (5.14) and (5.15). 

\vskip .3 truecm\noindent
{\bf Exercise 5.11} Show that $S(\alpha )$ is compact and that it is 
the universal cover of $S_{e^{2\pi i\alpha}}$. 

\vskip .3 truecm\noindent
{\bf Exercise 5.12} Show that if $f\in S_{1}$ then one has the 
following distorsion estimate: 
$$
\left| z{f'(z)\over f(z)}-1\right|\le {2|z|\over 1-|z|}\; 
\hbox{for all}\, z\in\D\; . \eqno(5.16)
$$

\vskip .3 truecm\noindent
{\bf Exercise 5.13} Use the result of the previous exercise to 
show that if $\varphi$ is as in (5.15) and $z\in\H$ then 
$$
\eqalignno{
|\varphi'(z)| &\le {2\exp (-2\pi\IM z)\over 1-\exp (-2\pi\IM z)}\; , 
&(5.17)\cr
|\varphi (z) | &\le -{1\over\pi}\log (1-\exp (-2\pi\IM z))\; . 
&(5.18)\cr
}
$$

\vskip .3 truecm\noindent
Let $r>0$, $\H_{r}=\H +ir$. It is clear that if $F\in S(\alpha )$
and $r$ is sufficiently large then $F$ is very close to the 
translation $z\mapsto z+\alpha $ for $z\in \H_{r}$. Indeed 
using the compactness of $S(\alpha )$ and Exercise 5.13 one can 
prove the following:

\vskip .3 truecm\noindent
{\bf Exercise 5.14} Let $\alpha\not= 0$. Show that there exists a 
universal constant $c_{0}>0$ (i.e. independent of $\alpha$) such that
for all $F\in S(\alpha )$ and for all $z\in\H_{t(\alpha )}$ where 
$$
t(\alpha )={1\over 2\pi}\log\alpha^{-1}+c_{0}\; , 
\eqno(5.19)
$$
one has 
$$
|F(z)-z-\alpha |\le {\alpha\over 4}\; . \eqno(5.20)
$$
[Hint: Let $\varphi (z)=F(z)-z-\alpha=\sum_{n=1}^\infty
\varphi_{n}e^{2\pi inz}$. If $\IM z>t(\alpha )$ then 
$|F(z)-\alpha -z|\le \sum_{n=1}^\infty |\varphi_{n}|\alpha^n
e^{-2\pi n c_{0}}$ thus $\ldots$.]

\vskip .3 truecm\noindent
Given $F$, the lowest admissible value $t(F,\alpha )$ of 
$t(\alpha )$ such that (5.20) holds for all 
$z\in\H_{t(F,\alpha )}$ represents the height in the upper half plane 
$\H$ at which the strong nonlinearities of $F$ manifest themselves. 
When $\IM z> t(F,\alpha )$, $F$ is very close to the translation 
$z\mapsto z+\alpha $. An example of strong nonlinearity is of course 
a fixed point: if $F(z)=z+\alpha +{1\over 2\pi i}e^{2\pi iz}$, 
$\alpha >0$, then $z=-{1\over 4}+{i\over 2\pi}\log (2\pi\alpha)^{-1}$ 
is fixed and $t(F,\alpha )\ge {1\over 2\pi}\log (2\pi\alpha)^{-1}$. 

The estimates (5.19) and (5.20) are the fundamental ingredient of 
Yoccoz's proof of the lower bound (5.2) together with Proposition 2.15
and the following elementary properties of the conformal capacity.

\vskip .3 truecm\noindent
{\bf Exercise 5.15} Let $U\subset C$ be a simply connected open set, 
$U\not= \C$, and let $z_{0}\in U$. Let $d$ be the distance of $z_{0}$ 
from the complement of $U$ in $\C$ and let $C(U,z_{0})$ denote the 
conformal capacity of $U$ w.r.t. $z_{0}$. Then 
$$
d\le C(U,z_{0})\le 4d\; . \eqno(5.21)
$$

\vskip .3 truecm\noindent
As in the proof of Douady--Ghys theorem we can now construct the first 
return map in the strip $B$ delimited by 
$l=[it(\alpha ),+i\infty [$, $F(l)$ and the segment $[it(\alpha ),
F(it(\alpha ))]$. Given $z$ in $B$ we can iterate $F$ until 
$\RE F^n(z)>1$. If $\IM z\ge t(\alpha )+c$ for some $c>0$ 
then $z'=F^n(z)-1\in B$ and $z\mapsto z'$ is the first return map 
in the strip $B$. Glueing $l$ and $F(l)$ by $F$ one obtains a 
Riemann surface $S$ corresponding to $\hbox{int}\, B$ and 
biholomorphic to $\D^{*}$. This induces a map $g\in S_{e^{2\pi i /
\alpha}}$ which lifts to $G\in S(\alpha^{-1})$. One can then 
show the following (see [Yo2], pp. 32--33)

\vskip .3 truecm\noindent
\Proc{Proposition 5.16}{Let $\alpha\in (0,1)$, $F\in S(\alpha )$ and 
$t(\alpha )>0$ such that if $\IM z\ge t(\alpha )$ then 
$|F(z)-z-\alpha |\le \alpha /4$. There exists $G\in S(\alpha^{-1})$
such that if $z\in \H$, $\IM z\ge t(\alpha )$ and $F^i(z)\in \H$ 
for all $i=0,1,\ldots ,n-1$ but $F^n(z)\not\in \H$ then there exists 
$z'\in\C$ such that 
\item{1.} $\IM z'\ge \alpha^{-1}(\IM z-t(\alpha )-c_{1})$, where 
$c_{1}>0$ is a universal constant;
\item{2.} There exists an integer $m$ such that $0\le m <n$ and 
$G^m(z')\not\in\H$.}

\vskip .3 truecm\noindent
From this Proposition one can conclude the proof as follows. 
Let us recall that $K_{f}=\cap_{n\ge 0}f^{-n}(\D )$ is the 
maximal compact $f$--invariant set containing $0$. Let 
$F\in S(\alpha )$ be the lift of $f\in S_{e^{2\pi i \alpha }}$ and 
let $K_{F}\subset \C$ be defined as the cover of $K_{f}$: 
$K_{F}=E^{-1}(K_{f})$. 
It is immediate to check that 
$$
d_{F}=\sup\{\IM z\,\mid\, z\in\C\setminus K_{F}\}=
-{1\over 2\pi}\log\hbox{dist}\, (0,\C\setminus K_{f})\; . 
$$
Thus by (5.21) one gets 
$$
\exp (-2\pi d_{F})\le C(K_{f},0)\le 4\exp (-2\pi d_{F})\; . 
$$
Theorem 5.1 is therefore equivalent to the lower bound 
$$
\sup_{F\in S(\alpha )} d_{F}\le {1\over 2\pi}B(\alpha )+C
\eqno(5.22)
$$
for some universal constant $C>0$. 

Assume that (5.22) is not true and that there exist 
$\alpha\in\R\setminus\Q\cap (0,1/2)$ with $B(\alpha )<+\infty$ , 
$F\in S(\alpha )$, $z\in \H$ and $n>0$ such that 
$$
\eqalign{
\IM F^n (z) &\le 0\; , \cr
\IM z &\ge {1\over 2\pi}B(\alpha )+C\; . \cr
}
$$

Let us choose $\alpha$, $F$ and $z$ so that $n$ is as small as possible. 
By Proposition 5.16, if $C>c_{0}$, one gets
$$
\eqalign{
\IM z' &\ge \alpha^{-1}[\IM z-t(\alpha )-c_{1}]\cr
&\ge \alpha^{-1}\left[{1\over 2\pi}(B(\alpha )-\log\alpha^{-1})+C
-c_{0}-c_{1}\right]\; . \cr
}
$$
By the functional equation of $B$ one gets 
$$
\IM z'\ge {1\over 2\pi}B(\alpha^{-1})+\alpha^{-1}[C-c_{0}
-c_{1}]\ge {1\over 2\pi}B(\alpha^{-1})+C
$$
provided that $C\ge 2(c_{0}+c_{1})$. But Proposition 5.16 shows that 
this contradicts the minimality of $n$ and we must therefore conclude
that (5.22) holds. \qed

\vskip .3 truecm \noindent
A nice description of the proof of the upper bound $\log r(\alpha 
)\le C-B(\alpha )$ can be found in the Bourbaki seminar of Ricardo 
Perez--Marco [PM1]. 

\line{}

\vskip .3 truecm\noindent
{\titsec 5.3 Some Open Problems}

\vskip .3 truecm\noindent
The first open problem we want to address is whether or not the 
infimum in (5.1) is attained by the quadratic polynomial 
$P_{\lambda}(z)=\lambda z\left( 1-{z\over 2}\right)$:

\vskip .3 truecm
\noindent
{\bf Question 5.17} Does 
$r(\alpha )= \inf_{f\in S_{e^{2\pi i\alpha}}}r(f)= r_{2}( e^{2\pi 
i\alpha})$, i.e. the radius of convergence of the quadratic 
polynomial ?

\vskip .3 truecm\noindent
If this were true then the very precise bound (5.11) would hold also
for $P_{\lambda}$: recalling that $r_{2}(\lambda )=|u(\lambda )|$, 
where $u\, :\D\rightarrow\C$ is the function defined in Section 
3.1, one can ask 

\vskip .3 truecm\noindent
{\bf Question 5.18} Does the function $\alpha\mapsto\log |u(e^{2\pi i
\alpha })|+B(\alpha )\in L^\infty (\S^{1})$?

\vskip .3 truecm
\noindent
Indeed there is a good numerical evidence that much more could be true:

\vskip .3 truecm
\noindent
{\bf Conjecture 5.19} The function $\alpha\mapsto\log |u(e^{2\pi i
\alpha })|+B(\alpha )$ extends to a H\"older $1/2$ function.

\vskip .3 truecm\noindent
We refer to [Ma1] and to [MMY2] and references therein for a 
discussion of Conjecture 5.19. The next Exercises give an idea of how 
to compute approximately but effectively the function 
$\alpha\mapsto\log |u(e^{2\pi i
\alpha })|$ on a computer. More informations can be found
in [He4] (where one can also 
find many problems, most of which are still open) and  [Ma1].

\vskip .3 truecm\noindent
{\bf Exercise 5.20} Let $f \in G_{\lambda}$ be 
linearizable, $\lambda = e^{2 \pi i \alpha}$.
Let $U_{f}$ be the Siegel disk of $f$, $h_{f}$ be the linearization 
of $f$, $z \in U_{f}$, $z = h_{f}(w)$, 
where $w \in \D_{r(f)}$, $|w| = r < r(f)$. Show that  
$$ 
\lim_{m \rightarrow +\infty} {1 \over m} \sum_{j=0}^{m-1}
       \log |f^{j}(z)| =\log r \eqno(5.23) 
$$
[Solution: Since $h_{f}$ conjugates 
$f$ to $R_{\lambda}$ one has 
$f^{j}(z) = f^{j}(h_{f}(w)) =
h_{f}(\lambda^{j} w)$ for all $j \ge 0$ and $w \in \D_{r(f)}$, thus
$$
{1 \over m} \sum_{j=0}^{m-1} \log |f^{j}(z)| = 
{1 \over m} \sum_{j=0}^{m-1} \log |h_{f}(\lambda^{j} w)| \, .
$$
$h_{f}$ has neither poles
nor zeros but $w=0$ thus by the mean property of harmonic 
functions one has $\int_{0}^{1} \log |h_{f}(r e^{2 \pi i x})| dx = \log r $
for all $r \le r(f)$. Finally note that  $w \mapsto \lambda w$ is uniquely 
ergodic on $|w| = r$, and in this case 
Birkhoff's ergodic theorem holds for all initial points, thus
$$ 
\eqalign{
\lim_{m \rightarrow +\infty} {1 \over m} \sum_{j=0}^{m-1} 
       \log |f^{j}(z)| &= 
	   \lim_{m \rightarrow +\infty} {1 \over m} 
       \sum_{j=0}^{m-1} \log | h_{f}(\lambda^{j} w)| \cr 
	   &=
       \int_{0}^{1} \log |h_{f}(r e^{2 \pi i x})| dx = \log r \, .\cr
	   }
$$

\vskip .3 truecm\noindent
{\bf Exercise 5.21} Deduce from the previous exercise that 
for almost every $z \in \partial U_{f}$ with respect to the 
harmonic measure one has 
$$ 
	\lim_{m \rightarrow +\infty} {1 \over m} \sum_{j=0}^{m-1} 
       \log |f^{j}(z)| = \log r(f) \; . \eqno(5.24)
$$

\vskip .3 truecm\noindent
Let us now consider the quadratic polynomial $P_{\lambda}$ once more. 
According to (5.24) to compute $|u(\lambda )|$ one needs to know that 
some point belongs to the boundary of the Siegel disk of $P_{\lambda}$
(and hope $\ldots$).
The critical point cannot be contained in $U_{P_{\lambda}}$ because
$f\vert_{U_{P_{\lambda}}}$ is injective, and from the classical theory of 
Fatou and Julia one knows that $\partial U_{P_{\lambda}}$ is contained in the 
closure of the forward orbit $\{ P_{\lambda}^{k}(1) \,\mid\, k \ge 0\}$ 
of the critical point $z=1$. 
Finally Herman proved if $\alpha$ verifies an arithmetical condition
${\cal H}$, weaker than the Diophantine condition but stronger 
than the Brjuno condition (see, for example, [Yo1] for its precise 
formulation) the critical
point belongs to  $\partial U_{P_{\lambda}}$. 

\vskip .3 truecm\noindent
{\bf Exercise 5.22} If $\alpha \in \,\cd{0}$ 
Herman has also proved that $\partial U_{P_{\lambda}}$ is a 
{\sl quasicircle}, that is the image of $\S^{1}$ under a 
quasiconformal homeomorphism. In this case $h_{P_{\lambda}}$ 
admits a quasiconformal extension to $|w| = r_{2}(\lambda )$ 
and therefore is H\"older continuous [Po]:
$$
	|h_{P_{\lambda}}(w_1) - h_{P_{\lambda}}(w_2)| \le 4 |w_1 - w_2|^{1 - \chi}
	\eqno(5.25)
$$
for all $w_1\, , w_2 \in \partial \D_{r_{2}(\lambda )}$, where $\chi
\in [0,1[$ depends on $\lambda $ is 
the so--called 
Grunsky norm [Po]  associated with 
the univalent function $g(z) = 
r_2(\lambda )/ h_{P_{\lambda}}(r_2(\lambda )/z)$ on $|z|>1$. 
Using this information show that 
$$ 
	\vert {1 \over q_{k}} \sum_{j=0}^{q_{k}-1} \log |P_{\lambda}^{j}(z)|
       - \log r_2(\lambda ) \vert \le {8 \over r_2(\lambda )} ({2 \pi \over 
       q_k})^{1 - \chi} \; , \eqno(5.26)
$$
where $p_{k}/q_{k}$ is a convergent of the continued fraction 
expansion of $\alpha$. Note that (5.26) implies  convergence 
to $\log r_2(\lambda )$ {\it for all} $z \in \partial U_{P_{\lambda}}$,
thus also for the critical point  $z=1$. 
%%%%% fine capitolo 5 %%%%%
\vfill\eject
%\vskip 1. truecm
%%%%% capitolo 6 %%%%%
\noindent
{\tit 6. Small divisors and loss of differentiability}

\vskip .3 truecm
\noindent
In this Chapter we will (very!) briefly illustrate other 
two completely different 
approaches to the problem of linearization of germs of holomorphic 
diffeomorphisms with an indifferent fixed point. 

In the previous chapter we saw how the optimal sufficient condition 
can be obtained by the classical majorant series method as Siegel and 
Brjuno did and that Yoccoz was able to show that it is also necessary 
with his ingenious creation of ``geometric renormalization''. 
Here we will give an idea of two proofs of the Siegel theorem, one 
due to Herman [He1, He2] and the other due essentially 
to Kolmogorov  [K] (see also Arnol'd [Ar3] and Zehnder [Ze2]). 

Herman's method is 
far from giving the optimal number--theoretical condition, the idea 
is simply so original and beautiful that it deserves being known. It 
also illustrates how in one--dimensional small divisor problems the 
problem known as ``loss of differentiability'' does not prevent from 
the application of simple tools like the contraction principle. 
Herman's method can also be extended to (and it is actually described by 
him for) the problem of local conjugacy to rotation of smooth
orientation--preserving diffeomorphisms of the circle. 

The idea of Kolmogorov is do adapt Newton's method for finding the 
roots of algebraic equations so as to apply it for finding the 
solution of the conjugacy equation. This method has been shown by 
R\"ussmann [R\"u] to be adaptable so as to prove the sufficiency
of a condition (slightly stronger than) Brjuno's. The main 
reason for sketching Kolmogorov's argument in this rather limited 
setting is that in the second part of this monograph we will 
illustrate Nash--Moser's implicit function theorem which is 
essentially the abstract and flexible formulation of Kolmogorov's idea. 

\line{}

\vskip .3 truecm\noindent
{\titsec 6.1 Hardy--Sobolev spaces and 
loss of differentiability }

\vskip .3 truecm\noindent
Let $k\in\N$, $r>0$. Following [He2], we introduce the Hardy--Sobolev spaces 
$$
{\cal O}^{k,2}_{r} = \{ f(z)=\sum_{n=0}^\infty f_{n}z^n\,\mid\,
\Vert f\Vert_{{\cal O}^{k,2}_{r}} = (|f_{0}|^{2}+\sum_{n=1}^\infty
n^{2k}|f_{n}|^{2}r^{2n})^{1/2}<+\infty\}\; . 
\eqno(6.1)
$$

\vskip .3 truecm\noindent
{\bf Exercise 6.1} 
\item{(a)} Show that ${\cal O}^{k,2}_{r}$ is a Hilbert space. 
\item{(b)} Show if  $f\in {\cal O}^{k,2}_{r}$, $f(0)=0$,  one has 
$$
\sup_{|z|\le r}|f(z)|\le \sqrt{\zeta (2k)}
\Vert f\Vert_{{\cal O}^{k,2}_{r}}\; ,
$$
where $\zeta$ denotes Riemann's zeta function. 
\item{(c)} Show that if $k\ge 1$ then ${\cal O}^{k,2}_{r}$ is a 
Banach algebra. 
\item{(d)} Show that if $f\in {\cal O}^{k,2}_{r}$ and $\phi$ is 
holomorphic in a neighborhood of $f(\D_{r})$ then $\phi\circ f
\in {\cal O}^{k,2}_{r}$ and on a sufficiently small neighborhood 
$V$ of $f$ in ${\cal O}^{k,2}_{r}$ the map $\psi\in V\mapsto
\phi\circ\psi\in {\cal O}^{k,2}_{r}$ is holomorphic. 

\vskip .3 truecm\noindent
The following very elementary proposition well illustrates 
the phenomenon of ``loss of differentiability'' 
due to the small divisors which already arises at the level of the 
linearized conjugacy equation (6.2). 

\vskip .3 truecm\noindent
\Proc{Proposition 6.2}{Let $0\le\tau\le \tau_{0}$, $\tau_{0}\in\N$, 
$k\ge 1+\tau_{0}$, $f\in {\cal O}^{k,2}_{r}$, $f(0)=0$. 
If $\alpha\in\cd{\tau}$ then the unique solution 
$g := D_{\lambda}^{-1}f$ verifying $g(0)=0$ of 
$$
g\circ R_{\lambda}-g=f \; , \eqno(6.2)
$$
belongs to ${\cal O}^{k-1-\tau_{0},2}_{r}$. Moreover there exists a 
universal constant $C>0$ such that 
$$
\Vert{d^{k-1-\tau_{0}}g\over dz}\Vert_{{\cal O}^{0,2}_{r}}
\le {C\over\gamma}\Vert {d^{k}f\over dz}\Vert_{{\cal O}^{0,2}_{r}}\; , 
$$
where $\gamma = \inf_{n\ge 1}n^{1+\tau}|\lambda^n-1|$. }

\vskip .3 truecm\noindent
\proof
It is a straightforward computation starting from the identity
$g(z)=\sum_{n=1}^\infty (\lambda^n-1)^{-1}f_{n}$. \qed

\vskip .3 truecm\noindent
This Proposition shows that solving the linear equation (6.2)
the small divisors cause the loss of $1+\tau_{0}$ derivatives. 
This loss of differentiability phenomenon is typical of 
small divisors problems and it will be crucial in the discussions 
in the second part of this monograph. The most annoying consequence 
of this phenomenon is the impossibility of using fixed 
points methods to solve conjugacy equations, simply because 
the operator $D_{\lambda}$ is {\it unbounded} if regarded on a {\it 
fixed} Hardy--Sobolev space. 
However under some restriction on $\tau_{0}$ one can actually 
use the contraction principle to solve the conjugacy 
problem, thanks to an ingenious idea of Herman we will shortly 
describe in the next section. 

\line{}

\vskip .3 truecm\noindent
{\titsec 6.2 Herman's Schwarzian derivative trick} 

\vskip .3 truecm\noindent
Let $\Omega$ be a region in the complex plane and $f\, : \Omega
\rightarrow\C$ be holomorphic. 

\vskip .3 truecm\noindent
\Proc{Definition 6.3}{The {\rm Schwarzian derivative} $S(f)$ of $f$ 
is}
$$
\eqalign{
S(f)  & := (\log f')''-{1\over 2}((\log f')')^{2} = 
{f'''\over f'}-{3\over 2}\left({f''\over f'}\right)^2 \cr
& = -2\sqrt{f'}\left({1\over\sqrt{f'}}\right)''\cr
}
\eqno(6.3)
$$

\vskip .3 truecm\noindent
{\bf Exercise 6.4} 
\item{(a)} Prove that the following ``chain rule'' holds: 
$S(f\circ g) = (S(f)\circ g)(g')^{2}+S(g)$. 
\item{(b)} Show that $S(f)\equiv 0$ if and only if 
$f$ is a M\"obius map. 

\vskip .3 truecm\noindent
The idea of Herman is to apply the Schwarzian derivative to the 
conjugacy equation $f\circ h = h\circ R_{\lambda}$: one obtains
$$
\lambda^2(S(h)\circ R_{\lambda})-S(h)=(S(f)\circ h)(h')^{2}\; . \eqno(6.4)
$$
At the r.h.s. appears $h'$, thus one has already lost one derivative 
and this does not seem to lead to anything good. Assuming the r.h.s. 
as given 
one could solve for $S(h)$ but this would cost $1+\tau_{0}$ derivatives
according to Proposition 6.2. However if $\tau_{0}=1$ (which is true 
for almost all $\alpha$ as we saw in Section 4.2) the total loss of 
derivatives is three. The idea now is that applying $S^{-1}$ one 
should recuperate three derivatives and this would imply that the 
map
$$
{\cal O}^{k,2}_{r}\ni h \mapsto
S^{-1}D_{\lambda}^{-1}[(S(f)\circ h)(h')^{2}]
\eqno(6.5)
$$ 
takes its values in ${\cal O}^{k,2}_{r}$ too. 
Note that we have slightly modified the definition of 
$D_{\lambda}$ with respect to the previous Section: here 
$D_{\lambda}f=\lambda^{2}f\circ R_{\lambda}-f$. This does not change 
the conclusions of Section 6.1. 

On a 
disk $\D_{r}$ of sufficiently small radius $f$ is close to 
$R_{\lambda}$ thus $S(f)$ must be small and one can hope to conclude using
a fixed point theorem (the contraction principle, say). This strategy
indeed works: see [He1, He2] for the details. 

The inversion of $S$ is achieved as follows. 
First of all note that it is not restrictive to 
assume $f''(0)=0$, so that $h''(0)=0$
(one can preliminarly conjugate $f$ by the polynomial
$z-{f''(0)\over 2(\lambda -\lambda^{2})}z^{2})$: this 
implies $[(S(f)\circ h)(h')^{2}]_{z=0}=0$.  Let 
$\psi = D_\lambda^{-1}[(S(f)\circ h)(h')^{2}]$. 
If  $\psi$ is small enough (this is always the case 
if one considers a sufficiently small disk) then one can write 
$\psi$ uniquely in the form 
$$
\psi = \psi_{1}'-{1\over 2}\psi_{1}^2 
$$
with $\psi_{1}(0)=0$. 
Now one can easily solve the problem 
$S(h_{1})=\psi=\psi_{1}'-{1\over 2}\psi_{1}^2 $ just looking for 
$h_{1}$ such that $(\log h_{1}')'=\psi_{1}$. This is achieved in three 
steps: 
$$
\eqalign{
\psi_{2}' &=\psi_{1}\; , \cr
\psi_{3} &=e^{c+\psi_{2}}\; \hbox{where}\, c \hbox{ is chosen s.t.}
\, \psi_{3}(0)=1\; , \cr
h_{1} &=\int_{0}^z\psi_{3}(\zeta )d\zeta\; . \cr
}
$$
Then it is immediate to check that $S(h_{1})=\psi$. 

\vskip .3 truecm\noindent
{\bf Exercise 6.5} Let $r>0$, $\psi$ and $h_{1}$ as above. Show that 
if $\psi\in {\cal O}^{0,2}_{r}$ then $h_{1}\in {\cal O}^{3,2}_{r}$. 

\line{}

\vskip .3 truecm\noindent
{\titsec 6.3 Kolmogorov's modified Newton method} 

\vskip .3 truecm\noindent
Here we follow quite closely [St], Volume II, Chapter III, Section 7. We suggest
however the reader to look also at [Ze2] for a complete proof. 

Let $f\in S_{\lambda}$, $\lambda = e^{2\pi i \alpha}$ and assume 
that $\alpha$ is a diophantine number. We want to construct $h$ 
tangent to the identity such that $R_{\lambda}-h^{-1}\circ f\circ 
h=0$. 
Let $\tilde{h}=h-\hbox{id}$ and let us define the composition 
law $\odot$ as
$$
(\hbox{id}+\tilde{h}_{1})\circ (\hbox{id}+\tilde{h}_{2}) = 
\hbox{id}+\tilde{h}_{1}\odot \tilde{h}_{2}\; . \eqno(6.6)
$$
Clearly one expects that 
$$
\tilde{h}_{1}\odot \tilde{h}_{2} =\tilde{h}_{1}+ \tilde{h}_{2}
+\hbox{quadratic terms}\; . 
$$
Let $\tilde{f}=f-R_{\lambda}$ and let us define a second composition 
law $\otimes$ as
$$
R_{\lambda }+\tilde{h}\otimes\tilde{f}=(\hbox{id}+\tilde{h})^{-1}
\circ (R_{\lambda}+\tilde{f})\circ (\hbox{id}+\tilde{h})\; . \eqno(6.7)
$$
Of course one needs $\tilde{h}$ be small so as to assure the existence 
of the inverse in (6.7) but this is not difficult to obtain 
considering a small enough disk since $\tilde{h}=h_{2}z^{2}+\ldots$. 

\vskip .3 truecm\noindent
{\bf Exercise 6.6} Recall Lagrange's Theorem on the inversion of 
analytic functions (see [Di], p. 250): if $h\, :\D_{r}\rightarrow\C$
is holomorphic and tangent to the identity then choosing $r$ small 
enough there exists a unique solution $z=\kappa (w)$ of the equation 
$w=h(z)$. Moreover $\kappa $ is holomorphic in a neighborhood of $0$ and 
is explicitly given by 
$$
\kappa (w)=\sum_{n=1}^\infty {(-1)^n\over n!}{d^{n-1}\over dw^{n-1}}
(h(w))^n\; . 
$$
Get some precise estimate on the size of the domain and of the 
norm of $\kappa$  if $h$ belongs to some 
Hardy--Sobolev space. 

\vskip .3 truecm\noindent
We then define 
$$
{\cal R}(\tilde{h},\tilde{f}) = 
(\hbox{id}+\tilde{h})\circ R_{\lambda}-(R_{\lambda}+\tilde{f})\circ
(\hbox{id}+\tilde{h})\; . \eqno(6.8)
$$
Clearly ${\cal R}(0,0)=0$, ${\cal R}(0,\tilde{f})=-\tilde{f}$ and 
$$
\eqalignno{
{\cal R}(\tilde{h}_{1}\odot \tilde{h}_{2},\tilde{f})
&= (\hbox{id}+\tilde{h}_{1})\circ (\hbox{id}+\tilde{h}_{2})\circ 
R_{\lambda} - (R_{\lambda}+\tilde{f})\circ (\hbox{id}+\tilde{h}_{1})
\circ (\hbox{id}+\tilde{h}_{2})\; , &(6.9) \cr
{\cal R}(\tilde{h}_{2}, \tilde{h}_{1}\otimes \tilde{f})
&= (\hbox{id}+\tilde{h}_{2})\circ R_{\lambda} - 
(\hbox{id}+\tilde{h}_{1})^{-1}\circ (R_{\lambda}+\tilde{f})\circ 
(\hbox{id}+\tilde{h}_{1})\circ (\hbox{id}+\tilde{h}_{2})\; . &(6.10) \cr
}
$$
Comparing (6.9) with (6.10) we have (for $z$ small enough)
$$
(\inf |1+\tilde{h}_{1}'|)|{\cal R}(\tilde{h}_{2},\tilde{h}_{1}\otimes
f)(z)|\le |{\cal R}(\tilde{h}_{1}\odot \tilde{h}_{2},\tilde{f})(z)|\le 
(\sup |1+\tilde{h}_{1}'|)|{\cal R}(\tilde{h}_{2},\tilde{h}_{1}\otimes
f)(z)|\; , 
$$
thus one should get 
$$
C^{-1}\Vert {\cal R}(\tilde{h}_{2},\tilde{h}_{1}\otimes f)\Vert\le 
\Vert {\cal R}(\tilde{h}_{1}\odot \tilde{h}_{2},\tilde{f})\Vert \le 
C \Vert {\cal R}(\tilde{h}_{2},\tilde{h}_{1}\otimes f)\Vert\; . \eqno(6.11)
$$
for a suitably chosen norm 
$\Vert \;\Vert$ and some $C>0$ (see Exercise 6.6). 

Let us now try to solve the equation ${\cal R}(\tilde{h},\tilde{f})=0$ by 
taking a sequence of approximations defined as follows:
\item{(0)} Let $\tilde{h}_{0}=\tilde{g}_{0}=0$, 
$\tilde{f}_{0}=\tilde{f}$: thus ${\cal R}(\tilde{h}_{0},\tilde{f}_{0})
=-\tilde{f}_{0}$; 
\item{(1)} Let $\tilde{f}_{1}=\tilde{g}_{0}\otimes\tilde{f}_{0}=
\tilde{f}_{0}$. 
Choose $\tilde{g}_{1}$ to be the solution of the 
linearized equation 
$$
\partial_{1}{\cal R}(0,0)\tilde{g}_{1}
+\partial_{2}{\cal R}(0,0)\tilde{f}_{1}=0\; , \eqno(6.12)_{1}
$$
where $\partial_{j}$ denotes the partial derivative w.r.t. the 
$j$--th argument. Finally we set 
$\tilde{h}_{1}=\tilde{h}_{0}\odot\tilde{g}_{1}$. 
\item{(i+1)} We choose $\tilde{f}_{i+1}=\tilde{g}_{i}\otimes\tilde{f}_{i}$
and $\tilde{g}_{i+1}$ to be the solution of 
$$
\partial_{1}{\cal R}(0,0)\tilde{g}_{i+1}
+\partial_{2}{\cal R}(0,0)\tilde{f}_{i+1}=0\; , \eqno(6.12)_{i+1}
$$
and we set $\tilde{h}_{i+1}=\tilde{h}_{i}\odot\tilde{g}_{i+1}$. 

It is immediate to check that the linearized 
equations (6.12)$_{i}$
have the form
$$
\tilde{g}_{i}\circ R_{\lambda}-R_{\lambda}\circ\tilde{g}_{i}=
\tilde{f}_{i}\; , \eqno(6.13)
$$
which we studied in Section 6.2. Note that we do {\it not} linearize 
at the point $(0,f)$ since we would get a difference equation
{\it without} constant coefficients: 
$$
\partial_{1}{\cal R}(0,\tilde{f}_{i-1})\tilde{g}_{i}+\partial_{2}{\cal R}(0,
\tilde{f}_{i-1})\tilde{f}_{i}= \tilde{g}_{i}\circ R_{\lambda}-
R_{\lambda }\circ \tilde{g}_{i}-\tilde{f}_{i-1}'\tilde{g}_{i}=0
\; . 
$$
If one could solve (6.13) at each step 
with a bound 
$$
\Vert \tilde{g}_{i}\Vert \le C\Vert \tilde{f}_{i}\Vert\; , \eqno(6.14)
$$
by (6.11) and (6.14) one would have 
$$
\Vert {\cal R}(\tilde{g}_{i}, \tilde{f}_{i})\Vert \le
\sup\Vert d^{2}{\cal R}\Vert (\Vert \tilde{g}_{i}\Vert^{2}+
\Vert\tilde{f}_{i}\Vert^{2}\le C^{3}\Vert {\cal R}(
\tilde{g}_{i-1}, \tilde{f}_{i-1})\Vert^{2}\; . \eqno(6.15)
$$
This would imply the convergence of the iterative scheme to a 
solution of ${\cal R}(\tilde{h}, \tilde{f})=0$ provided that 
one chooses $\Vert \tilde{f}\Vert $ small enough
(i.e. one considers the restriction of $f$ to a small enough
disk $\D_{r}$): indeed iterating  (6.15) one gets
$$
\Vert {\cal R}(\tilde{g}_{i}, \tilde{f}_{i})\Vert \le
(C^{3/2}\Vert {\cal R}(\tilde{g}_{0}, \tilde{f}_{0})\Vert)^{2^{i}}
$$
thus again by (6.11) one has 
$$
\eqalign{
\Vert {\cal R}(\tilde{h}_{i}, \tilde{f})\Vert &= 
\Vert {\cal R}(\tilde{g}_{0}\odot\tilde{g}_{1}\odot
\ldots\odot\tilde{g}_{i}, \tilde{f})\Vert \cr
&\le C \Vert {\cal R}(\tilde{g}_{1}\odot
\ldots\odot\tilde{g}_{i}, \tilde{g}_{0}\otimes \tilde{f})\Vert\cr
&\le C^{i}\Vert{\cal R}(\tilde{g}_{i}, \tilde{f}_{i})\Vert\cr
&\le (C^{2}\Vert {\cal R}(\tilde{g}_{0}, 
\tilde{f}_{0})\Vert)^{2^{i}}\cr
}
$$

\vskip .3 truecm\noindent
{\bf Exercise 6.7} Assuming the estimates above show that the 
sequence $\tilde{h}_{n}$ converges thus by continuity of ${\cal R}$ 
one gets the desired result.

\vskip .3 truecm\noindent
{\bf Exercise 6.8} Use the above scheme to give an alternative proof 
of Koenigs--Poincar\'e theorem.

\vskip .3 truecm\noindent
The above discussion shows how to prove the existence of the 
linearization disregarding the problem of loss of differentiability 
due to small divisors. This makes impossible to get an estimate 
like (6.14) {\it unless one regularizes the r.h.s.}. The simplest
method of regolarization, which is adapted to the analytic case, 
is to {\it consider restrictions} of the domains: 

\vskip .3 truecm\noindent
{\bf Exercise 6.9} Show that if $f\in{\cal O}^{0,2}_{r}$, $f(0)=0$, 
$k\in\N$, for all $\delta >0$ one has 
$$
\Vert f\Vert_{{\cal O}^{k,2}_{re^{-\delta}}}\le 
\left({k\over\delta}\right)^ke^{-k}\Vert f\Vert_{{\cal O}^{0,2}_{r}}
\; . \eqno(6.16)
$$

\vskip .3 truecm\noindent
Combining the above given discussion with a suitable choice of 
restrictions (i.e. a sequence $(\delta_{n})_{n\ge 0}$ such that 
$\sum_{n=0}^\infty\delta_{n}<+\infty$) one can indeed prove Siegel's 
Theorem following the iteration method. 
%%%%% fine capitolo 6 %%%%%
\vfill\eject
%%%%% capitolo 7: sistemi hamiltoniani %%%%%
\noindent
\centerline{\tit Part II. Implicit Function Theorems and KAM Theory}

\vskip 2. truecm
\noindent
{\tit 7. Hamiltonian Systems and Integrable Systems}

\vskip .3 truecm
\noindent
In this Chapter  we will very briefly recall some 
well--known 
facts on symplectic manifolds and Hamiltonian systems. 
Very good references are [AKN] and [AM]. 

\line{}

\vskip .3 truecm
\noindent
{\titsec 7.1 Symplectic Manifolds and Hamiltonian Systems }

\vskip .3 truecm
\noindent
\Proc{Definition 7.1}{ A ${\cal C}^\infty$ {\rm symplectic manifold}
is a $2l$--dimensional ${\cal C}^\infty$ manifold $M$ equipped with 
a non--degenerate ${\cal C}^\infty$ two--form (the symplectic form) $\omega$.
A ${\cal C}^\infty$ map $f\, :U\rightarrow M'$ where $U\subset M$ is open
and $M'$ is also symplectic (with symplectic form $\omega'$) is 
{\rm symplectic}
(or {\rm canonical}) if $f^*\omega' =\omega$.}

\vskip .3 truecm\noindent
The simplest (but important) examples of symplectic manifolds are:
\item{$\bullet$} $M=\R^{2l}\ni (p_1,\ldots ,p_l,q_1,\ldots ,q_l)$, 
$\omega= \sum_{i=1}^l dp_{i}\wedge dq_{i}$ ({\it standard symplectic 
structure}).
If $U$ and $V$ are two open sets in $\R^{2l}$ and $f\, :U\rightarrow V$
then $f$ is symplectic if and only if its Jacobian matrix 
$J_f\in\hbox{Sp}\, (l,\R)$,
the Lie group of $2l\times 2l$ real matrices $A$ such that 
$A^T{\cal I}A={\cal I}$,
where ${\cal I}=\left(\matrix{ 0 & -1\cr 1 & 0\cr}\right)$. 
\item{$\bullet$} $M=T^*N$ where $N$ is a ${\cal C}^\infty$ 
Riemannian manifold.
This is the typical situation in classical mechanics. If $(q_1,\ldots ,q_l)$ 
are local coordinates in $N$ and $(p_1,\ldots ,p_l)$ are the corresponding 
local coordinates in the cotangent space at a point, then $\omega=
\sum_{i=1}^l dp_{i}\wedge dq_{i}$.
\item{$\bullet$} $M=\T^{2l}$, $\omega =\sum_{i=1}^ld\theta_i\wedge d
\theta_{i+l}$.

\vskip .3 truecm\noindent
\Proc{Theorem 7.2 (Darboux)}{Each symplectic manifold $M$ has an atlas 
$(U_\alpha ,\varphi_\alpha )_{\alpha\in {\cal A}}$ such that 
on $\varphi_\alpha (U_\alpha )\subset\R^{2l}$ one has 
$\omega =\varphi_\alpha^* \sum_{i=1}^l dp_{i}\wedge dq_{i}$ (the 
standard symplectic structure on $\R^{2l}$). The transition maps 
$\varphi_\alpha\circ\varphi_\beta^{-1}$ are symplectic 
diffeomorphisms, i.e. their Jacobians $J_{\alpha\beta}(x)\in
\hbox{Sp}\, (l,\R)$ for all $x\in \varphi_\beta (U_\alpha\cap U_\beta )$.}

\vskip .3 truecm\noindent
The atlas given by Darboux's Theorem and the corresponding local coordinates 
are called symplectic. 

\vskip .3 truecm\noindent
\Proc{Definition 7.3}{A {\rm Hamiltonian function} on a symplectic 
manifold $(M,\omega )$ is a function $H\in {\cal C}^\infty (M,\R )$. 
The {\rm Hamiltonian vector field} associated to $H$ is the unique 
$X_H\in  {\cal C}^\infty (M,TM)$ such that $i_{X_H}\omega =dH$.}

\vskip .3 truecm\noindent
Note that in symplectic local coordinates a Hamiltonian vector field
takes the form 
$$
X_H = \sum_{i=1}^l -{\partial H\over\partial q_i}{\partial\over\partial
p_i}+{\partial H\over\partial p_i}{\partial\over\partial
q_i}\; ,\eqno(7.1)
$$
and the associated ordinary differential equations are the classical 
Ha\-mil\-ton's equations of the motion 
of a conservative mechanical system with $l$ degrees of freedom: 
$$
\dot{p_{i}}= -{\partial H\over\partial q_{i}}\; , \;\;
\dot{q_{i}} = {\partial H\over\partial p_{i}}\; , \;\; 
1\le i\le l\; . \eqno(7.2)
$$
Clearly the Hamiltonian function is a first integral of (7.2). 
The coordinates $q_i$ are also called ``generalized coordinates'' 
and the $p_i$ their ``conjugate momenta''. 
In many problems arising 
from celestial mechanics the flow is not complete due to the 
unavoidable occurance of collisions, but we will always assume 
completeness of the Hamiltonian flow. 

\vskip .3 truecm
\noindent
\Proc{Definition 7.4}{
The {\rm Poisson 
bracket} of two functions $f,g\in {\cal C}^\infty (M,\R )$ 
defined on an open 
subset of $(M,\omega )$ is
$\{F,G\} := X_{G}F = \omega (X_{F},X_{G}) = -X_{F}G\; , $
thus 
$X_{\{F,G\}}=-[X_{F},X_{G}]\; . $
Two functions $F,G$ are {\rm in involution} if $\{F,G\}=0$, i.~e. 
when their hamiltonian flows commute. }

\vskip .3 truecm
\noindent
{\bf Exercise 7.5} Show that the 
{\it Hamiltonian flow} $\Phi\, :\R\times M\rightarrow M$ is symplectic:
for all $t\in\R$ one has $\Phi (t,\cdot )^*\omega =\omega$. 
[Hint: use Cartan's formula ${d\over dt}\vert_{t=0}\Phi (t,\cdot )^{*}
\omega = d(i_{X_{H}}\omega )+i_{X_{H}}d\omega$, where 
$X_{H}={d\over dt}\vert_{t=0}\Phi (t,\cdot )$ is the Hamiltonian 
vector field associated to $\Phi$.]

\vskip .3 truecm\noindent
The importance of Exercise 7.5 is that to make symplectic coordinate 
changes of a Hamiltonian vector field it is sufficient to change the 
varables in the corresponding Hamiltonian function. This is a simpler
operation, both conceptually and computationally. 

As we will see in the next Section, among the possible orbits 
of Hamiltonian systems, {\it quasiperiodic} orbits are of special 
interest.

\vskip .3 truecm
\noindent
\Proc{Definition 7.6}{A continuous function $F\, :\R\rightarrow\R$
is {\rm quasiperiodic} if there exist $n\ge 2$, $f\, :\T^n\rightarrow
\R$ continuous  and $\nu\in\R^n\setminus\{0\}$ such that $F(t)=
f(\nu_1t,\ldots \nu_nt)$.}

\vskip .3 truecm\noindent
Let ${\cal M}=\{k\in\Z^n\, \mid\, k\cdot\nu=0\}$. Note that 
${\cal M}$ is a $\Z$--module. If $\dim {\cal M}=n$ then $\nu =0$, 
if $\dim {\cal M}=n-1$ then there exists $\alpha\in\R$ and $k\in\Z^n$
such that $\nu=\alpha k$. If $\dim {\cal M}=0$ then $\nu$ is 
called {\it non--resonant}. 

\vskip .3 truecm\noindent
{\bf Exercise 7.7} Show that the closure of any orbit of the linear flow
$\dot{\theta}=\nu$ on $\T^n$ is diffeomorphic to the torus 
$\T^{n-\dim {\cal M}}$. 

\vskip .3 truecm\noindent
{\bf Exercise 7.8} Show that if $\dim {\cal M}\in\{1,\ldots n-1\}$ 
there exists $A\in\hbox{SL}\,(n,\Z )$ such that posing
$\varphi = A\theta$ the linear flow $\dot{\theta}=\nu$ on $\T^n$
becomes $\dot{\varphi}_i = 0$ for $i=1,\ldots ,m$ and 
$\dot{\varphi}_i=\nu_i'$ for $i=m+1,\ldots ,n$ with 
$(\nu_{m+1}',\ldots ,\nu_n')\in\R^{n-m}$ non--resonant.

\vskip .3 truecm\noindent
{\bf Exercise 7.9} Show that if $\nu$ is non--resonant 
then the Haar measure on $\T^n$ is uniquely ergodic
(see [Mn] for its definition) for the linear flow 
$\dot{\theta}=\nu$ on $\T^n$. 

\line{}

\vskip .3 truecm\noindent
{\titsec 7.2 Integrable Systems} 

\vskip .3 truecm
\noindent
An especially interesting example of 
symplectic manifold is  $M={\Bbb 
R}^l\times {\Bbb 
T}^l$ which can be identified with the cotangent bundle of the 
$l$--dimensional torus ${\Bbb T}^l = {\Bbb R}^l/(2\pi {\Bbb Z})^l$. This 
manifold has a natural symplectic structure defined by the closed 
$2$--form 
$\omega = \sum_{i=1}^l dJ_{i}\wedge d\vartheta_{i}$ where 
$(J_{1},\ldots J_{l},\vartheta_{1},\ldots \vartheta_{l})$ are 
coordinates  on ${\Bbb R}^l\times {\Bbb T}^l$. 

\vskip .3 truecm\noindent
\Proc{Definition 7.10}{Let $U$ denote an open connected subset of 
${\Bbb R}^l$. Whenever an  Hamiltonian system can 
be reduced by a symplectic change of coordinates to a function 
$H\, : U\times {\Bbb T}^l\rightarrow {\Bbb R}$ which {\rm does not 
depend on the angular variables} $\vartheta$ one says that the system 
is {\rm completely canonically integrable} and the variables $J$ are 
called {\rm action variables}. }

\vskip .3 truecm
\noindent
Note that in this case Hamilton's 
equations (7.2) take the particularly simple form 
$$
\dot{J_{i}} = - {\partial H\over\partial \vartheta_{i}} = 0\; , \;\;\;
\dot{\vartheta_{i}} = {\partial H\over\partial J_{i}}\; , \;\;\;
i=1,\ldots ,l
$$
and the flow leaves invariant the $l$--dimensional torus $J=$ 
constant. The motion is therefore bounded and quasiperiodic (or 
periodic). 

Being completely canonically integrable is 
a stronger requirement than integrability by            
quadratures or complete integrability
(see [AKN] for their discussion). In the latter case one             
requires the existence of $l$ independent first integrals in involution   
but their joint level--set may well be non compact                
(this is already the case in the two body problem for non negative        
energy values) and the flow does not need to be quasiperiodic             
(scattering states).                                                      

The main risult in the theory of completely canonically 
integrable systems is the celebrated 

\vskip .3 truecm\noindent
{\bf Theorem 7.11 (Arnol'd--Liouville)}{\it
Let $H\in {\cal C}^\infty (M,\R )$ and assume that  
$F_1,\ldots ,F_l\in {\cal C}^\infty (M,\R )$ are 
$l$ first integrals in involution for the Hamiltonian 
flow associated to $H$. Let $a\in\R^l$
be such that $M_a=\{m\in M\,\mid \, F_i(m)=a_i\,\forall i=1,\ldots ,l\}$
is not empty and assume that the $l$ functions 
$F_1,\ldots ,F_l$ are independent\footnote{$^{1}$}{\rm
As usual $F_1,\ldots ,F_l$ are independent if $dF_{1}\wedge\ldots\wedge
dF_{l}\not= 0$.}
in a neighborhood of $M_{a}$. 
Then if $M_a$ is compact and connected it is
diffeomorphic to the $l$--torus. Moreover there exists an invariant 
open subset $V$ of $M$ which contains $M_a$ and is symplectically
diffeomorphic to $ U\times {\Bbb T}^l$, where $U$ is an open subset of 
$\R^l$. }

\vskip .3 truecm\noindent
Arnol'd--Liouville's Theorem thus assures that the existence of 
sufficiently many first integrals together with the compactness and 
connectedness of their level set guarantees complete canonical integrability.

\line{}

\vskip .3 truecm\noindent
{\titsec 7.3 Examples of completely canonically integrable 
systems}

\vskip .3 truecm\noindent
In this section we will briefly describe some examples of 
completely canonical integrable systems.

\vskip .3 truecm\noindent
{\bf Example 7.12: Harmonic oscillators.}
Let $M=\R^{2l}$ with the standard symplectic structure, 
$S\in \hbox{GL}\, (2l,\R )$ be symmetric and positive definite.
Consider the Hamiltonian system $H(x)={1\over 2}x^TSx$. 
This is completely integrable. Indeed if $J$ is a symplectic  
matrix which diagonalizes  $S$, in the variables 
$y=J^{-1}x$ the Hamiltonian will be 
$$
H(y)=\sum_{i=1}^{2l}
{\lambda_i y_i^2+\lambda_{i+l}y_{i+l}^2\over 2}\; , 
$$ 
where $\lambda_i\, ,
 \, i=1,\ldots ,2l$ are the eigenvalues of $S$. 
Then the functions $F_i={\lambda_i y_i^2+\lambda_{i+l}y_{i+l}^2\over 2}$,
$i=1,\ldots ,l$, 
are independent first integrals in involution and their common level 
set is compact and connected (since $\lambda_i>0$ for all $i$). The
symplectic transformation to action--angle variables is
$$
y_i=\sqrt{2J_i\sqrt{\lambda_{i+l}/\lambda_i}}\cos\chi_i\; , \;\;\; 
y_{i+l}=
\sqrt{2J_i\sqrt{\lambda_{i}/\lambda_{i+l}}}\sin\chi_i\; , \;\;\;
i=1,\ldots ,l\; .
$$

\vskip .3 truecm\noindent
{\bf Example 7.13: The two body problem.}
The Hamiltonian ${\cal H}\,:T^{*}({\Bbb R}^{3}\setminus \{ 0\})\mapsto
{\Bbb R}$
of the two-body problem in the center
of mass frame is (we have assumed $G=1$, where $G$ is the universal 
gravitational constant)
$$
{\cal H}( p, q)={1\over 2\mu}\Vert  p\Vert^{2}-{m_{0}m\over
\Vert  q\Vert}
$$
where $\mu=m_{0}m/(m_{0}+m)$ is the reduced mass of the system.

It is well-known that for negative energy the solutions 
are ellipses with one focus at the origin (i.~e. the center of mass).
These are called {\it keplerian orbits}. 
The shape and the position of the ellipse in space  are 
determined from the knowledge of the major semiaxis $a$, 
the eccentricity $e$, the angle of inclination $i$ of its plane w.r.t. the 
horizontal plane $q_{3}=0$, the argument of perihelion $\omega$ and 
the longitude of the ascending node $\Omega$. The position of the 
planet along the ellipse is determined by the mean anomaly $l$,
which is proportional to the area swept by the position 
vector $q$ of the planet starting from the perihelion. 

The systems admits $5$ {\it independent} first integrals: the total energy
${\cal H}$, the three components of the angular momentum $q\wedge  
p$ and one of the components of the 
Laplace--Runge--Lenz vector $A = p\wedge q\wedge 
p-{m_{0}m q\over
\Vert q\Vert}$. Among these integrals one can choose three 
integrals in involution and construct the completely canonical 
transformation to action--angle variables.  The other two integrals are 
responsible for the  proper complete degeneration of the Kepler 
problem: one can choose action--angle variables so that the 
Hamiltonian depends only on one of the actions. 
Indeed the Delaunay action--angle variables $(L,G,\Theta ,l,g,\theta )$
are related to the orbital elements as follows:
$$
L =\mu \sqrt{(m_{0}+m)a}\; ,\;\; G= L\sqrt{1-e^{2}}\; , \;\;
\Theta = G\cos i\; , \;\; l\; , \;\; g=\omega \; , \;\;
\theta = \Omega\; .
$$
Note that $G$ is the modulus of angular momentum
$ q\wedge  
p$, thus $\Theta$ is its projection along the $q_{3}$--axis. One has
the obvious limitation $|\Theta |\le G$. The
new Hamiltonian reads ${\cal H}=-{\mu^{3}(m_{0}+m)^{2}\over 2L^{2}}$.

The relation among Delaunay variables and the original 
momentum--position $( p, q)$ variables is much more 
subtle and will not be discussed here. 

The two--body problem is the modelization of the motion of
a planet around the Sun. But the 
Delaunay variables are not suitable for 
the description of the orbits of the planets of the solar system
since they are singular for circular orbits ($e=0$, thus $L=G$
anf the argument of the perihelion $g$ is not defined)
and for horizontal orbits ($i=0$ or $i=\pi$, thus 
$G=\Theta$ and the longitude of the ascending node $\theta$ is not 
defined). 
But all the planets of the solar system have almost circular orbits 
(with the exception of Mercury and Mars) and small inclinations.

Poincar\'e solved the problem first introducing a new set of 
action--angle variables $(\Lambda, H, Z, \lambda, h, \zeta)$: 
$\Lambda = L$, $H=L-G$, $Z=G-\Theta$, $\lambda = l+g+\theta$, 
$h=-g-\theta$, $\zeta = -\theta$
($\lambda$ is called the mean longitude, $-h$ is the longitude of the 
perihelion) then considering the couples $(H,h)$ and $(Z,\zeta )$ as 
polar symplectic coordinates: 
$$
\xi_{1} = \sqrt{2H}\cos h\; , \;\;\; \eta_{1} = \sqrt{2H}\sin h\; , 
\;\;\;\;\;\;\;  \xi_{2}=\sqrt{2Z}\cos\zeta \; , \;\;\;
\eta_{2}=\sqrt{2Z}\sin\zeta\; . 
$$
The variables $(\Lambda , \xi,  \lambda, \eta )$ are called 
{\it Poincar\'e variables}. They are well defined also in the case of 
circular ($H=0$) or horizontal ($Z=0$) orbits. 

\vskip .3 truecm\noindent
{\bf Example 7.14: Motion of a ``heavy'' particle on a surface of revolution.}
Let $S\subset \R^3$ be a surface of revolution with the Riemannian metric
induced by its embedding into $\R^3$ and assume that 
$x_3$ is its symmetry axis. Let $f\in{\cal C}^\infty
(\R,\R )$. If the surface never meets the $x_3$ axis 
then it is diffeomorphic to the cylinder $S\approx\R\times\T^1$
and its cotangent bundle will be $\R^3\times\T^1$. 
If $(p_1,p_2,q_1,q_2)$ are symplectic coordinates the Hamiltonian of 
a (heavy) point mass constrained to move on $S$ is 
$H(p,q)={1\over 2}\Vert p\Vert^2+f(q_1)$. Here $f$ is the ``weight''
and $p_2$ (which corresponds to the projection of the angular momentum of 
the particle along the $x_3$--axis) is an independent integral of the
motion. Complete integrability is assured if 
the curve $\{(p_1,q_1)\in\R^2\; \mid\; 
H(p_1,a,q_1,q_2)=E\}$ is closed for some value of $a$ and $E$.
%%%%% fine capitolo 7 %%%%%
\vfill\eject
%%%%% capitolo 8: sistemi quasi integrabili %%%%%
\noindent
{\tit  8. Quasi--integrable Hamiltonian Systems} 

\vskip .3 truecm
\noindent
The importance of completely canonically integrable Hamiltonian 
systems is due both to the fact that their flows can be studied in 
great detail and that many problems in mathematical physics can be 
considered as {\it perturbations} of integrable systems. The most 
famous example is given by the motion of the planets in the Solar 
System (see [Ma2] and references therein for an introduction). If the 
(weak) mutual attraction between the planets is neglected the system 
decouples into several independent Kepler problems and it is 
completely integrable. Exactly this problem gave origin in the 18th 
century to ``perturbation theory'' whose modern formulation is mainly 
due to the monumental work of Henri Poincar\'e [P]. The goal of 
pertubation theory is to understand the dynamics of a ``perturbed''
system which is close to a well--understood one (usually
an integrable system).

\line{}

\vskip .3 truecm
\noindent
{\titsec 8.1 Quasi--integrable Systems} 

\vskip .3 truecm\noindent
Following Poin\-ca\-r\'e [P], the {\it fundamental problem of dynamics} is the 
study of quasi--in\-te\-gra\-ble Hamiltonian systems: let 
$\varepsilon_{0}>0$, 

\vskip .3 truecm\noindent
\Proc{Definition 8.1}{A {\rm quasi--integrable} Hamiltonian 
system is a function ${\cal H}\in{\cal C}^\infty ((-\varepsilon_0,
\varepsilon_0)\times M,\R )$ such that the Hamiltonian function 
$H={\cal H}(0,\cdot )\, : M\rightarrow\R$ 
is completely canonically integrable.}

\vskip .3 truecm\noindent
Using the canonical transformation 
to action--angle variables associated to 
${\cal H}(0,\cdot )$, Hamiltonians 
${\cal H}\, :(-\varepsilon_0,
\varepsilon_0)\times U\times {\Bbb T}^l\mapsto {\Bbb R}$ (smooth or analytic)
of the form 
$$
{\cal H}(\varepsilon, J,\chi ) = 
h_{0}(J)+\varepsilon f(J,\chi )\; , \eqno(8.1)
$$
where $f\in{\cal C}^\infty (U\times\T^l,\R)$,
are typical examples of quasi--integrable Hamiltonian systems. 

The most ambitious program would be to prove that quasi--integrable Hamiltonian
systems are indeed integrable: i.e. to show that there exists a one--parameter
family $V_\varepsilon$ of open connected invariant subsets of $M$ which 
are symplectically diffeomorphic to $U_\varepsilon\times\T^l$
where $U_{\varepsilon}\subset\R^l$ is open, connected 
and such that if $(\tilde{J},\tilde{\chi})$ are the coordinates in
$U_\varepsilon\times\T^l$ one has 
${\cal H}(\varepsilon,\cdot ,\cdot )\mid_{V_\varepsilon}=
h_{\varepsilon}( \tilde{J})$ for some smooth one--parameter family of 
smooth function $h_{\varepsilon}\, : U_{\varepsilon}\times\T^l
\rightarrow\R$. 

In general this is asking too much: a result of Poincar\'e 
shows that in general quasi--integrable Hamiltonian systems are  
not completely integrable (in addition to [P], Tome I, Chapitre V, 
see [BFGG] for  a nice discussion of the consequences of this problem 
and a related result of Fermi). 

\vskip .3 truecm\noindent
\Proc{Theorem 8.2 (Poincar\'e)}{Consider a 
quasi--integrable Hamiltonian of the 
form (8.1), $l\ge 2$. Assume that  
the two following genericity assumptions are satisfied:
(1) {\rm non--degeneracy}: $\det \left({\partial^2 h_0\over
\partial J_i\partial J_k}\right)\not= 0$ on $U$;
(2) {\rm generic perturbations}:
for all $J\in U$ and for all $k\in \Z^l\setminus \{0\}$
either the $k$--th Fourier coefficient $\hat{f}_k(J)$ of $f$ 
does not vanish or there exists $k'\in \Z^l\setminus \{0\}$
parallel to $k$ such that  $\hat{f}_{k'}(J)\not= 0$. Then 
the system is not a smooth one--parameter family of 
completely canonically integrable Hamiltonians. }

\vskip .3 truecm
\noindent
One can also recall the following 
theorem of Markus and Meyer [MM]

\vskip .3 truecm\noindent
\Proc{Theorem 8.3}{Generically hamiltonian systems are neither 
completely canonically integrable nor ergodic}

\vskip .3 truecm
\noindent
{\bf Exercise 8.4} Prove Poincar\'e's Theorem following these lines. 
Using the notations introduced above, if the system were completely
canonically integrable then the new actions $\tilde{J}$ would be 
a system of $l$ independent first integrals of the Hamiltonian 
flow of ${\cal H}$ in involution. Writing them explicitly
in terms of the old local coordinates $(J,\chi )$ one has
$$
\tilde{J} = J+\varepsilon \tilde{J}_{1}(J,\chi )
+{\cal O}(\varepsilon^{2})\; , \eqno(8.2)
$$
for some smooth function $\tilde{J}_{1}\, : \, 
V_{\varepsilon}\rightarrow U_{\varepsilon}$. Imposing that 
$\{\tilde{J},{\cal H}\}={\cal O}(\varepsilon^{2})$ leads 
to the system of linear partial differential equations
$$
\sum_{i=1}^l{\partial h_{0}\over\partial J_{i}}
{\partial \tilde{J}_{1j}\over\partial\chi_{i}} = 
{\partial f\over \partial\chi_{j}}\; , \;\;\;
j=1, \ldots ,l.\eqno(8.3)
$$
Using Fourier series try to find a smooth solution to these equations
$\ldots$. 

\line{}

\vskip .3 truecm\noindent
{\titsec 8.2 Constant Coefficients Linear PDE on $\T^n$
and Loss of Differentiability.}

\vskip .3 truecm\noindent
The (very) short sketch of the proof of Poincar\'e's Theorem led us to 
consider the general constant coefficients 
linear partial differential 
equation on $\T^n$
$$
D_\mu u:= \mu\cdot\partial u=v\; , 
\eqno(8.4)
$$
where $\mu\in\R^n$, 
$\partial u=(\partial_1 u,\ldots ,\partial_n u)$, 
is the gradient of $u$, $v\in {\cal C}^{0,\infty} (\T^n,\R^m)$
(i.e. $v\in {\cal C}^{\infty} (\T^n,\R^m)$ {\it and}
$\int_{\T^n} v(\theta )d\theta =0$). 
Indeed for all fixed value of $J$ the equation (8.3) is a special 
case of (8.4) with $n=m=l$. 

It is easy to check (see Appendix A3 for a detailed 
discussion of the case $n=2$) that $D_\mu$ is hypoelliptic
\footnote{$^{1}$}{A constant coefficients linear partial differential 
operator $P$ is hypoelliptic if all $u$ such that 
$Pu=v$ are ${\cal C}^\infty$ on all open 
sets where $v$ is ${\cal C}^\infty$ (see [H1], p.109).}
if and only if $\mu$ is a diophantine vector, i.e. there exist
two constants $\gamma >0$ and $\tau\ge n-1$ such that 
$$
|\mu\cdot k |\ge \gamma |k|^{-\tau}\;\forall k\in\Z^n\setminus
\{0\}\; , \eqno(8.5)
$$
where $k=(k_1,\ldots k_n)$, $|k|=|k_1|+\ldots +|k_n|$. 

\vskip .3 truecm\noindent
{\bf Exercise 8.5} Prove that almost all $\mu\in\R^n$
is diophantine of exponent $\tau >n-1$.

\vskip .3 truecm\noindent
{\bf Exercise 8.6} Have a look to the book of 
Y. Meyer [Me]. Among many interesting things one finds
the following theorem (Proposition 2, p. 16): Let ${\cal R}$
be a real algebraic number field and let $n$ be its degree over 
$\Q$. Let $\sigma$ be the $\Q$--isomorphism of ${\cal R}$
such that $\sigma ({\cal R})\subset\R$ and let $\mu_1,\ldots \mu_n$
be any basis of ${\cal R}$ over $\Q$. Then $(\sigma (\mu_1),
\ldots ,\sigma (\mu_n))\in\R^n$ is diophantine of exponent
$\tau =n-1$. Try to prove it if you remember a tiny bit of Galois 
theory. Apply it to $(1,\sqrt{2},\sqrt{3},\sqrt{6})$ and
$(1,2^{1/3},2^{2/3})$. There exist also higher dimensional 
generalizations of Roth's Theorem quoted in Exercise 4.9: see, for 
example, the Subspace Theorem [Sch1, Sch2]. 

\vskip .3 truecm\noindent
In addition to knowing that $u\in {\cal C}^{\infty} (\T^n,\R^m)$ one has the 
following more precise estimate: 

\vskip .3 truecm\noindent
\Proc{Proposition 8.7}{Let $\Vert\;\Vert_k$ denote the 
${\cal C}^k$ norm. If $\mu$ is diophantine with exponent 
$\tau $ then for all $r>\tau +n-1$ and for all $i\in \N$
there exists a positive constant $A_i$ such that }
$$
\Vert u\Vert_i\le A_i\Vert v\Vert_{i+r}\; . 
\eqno(8.6)
$$

\vskip .3 truecm\noindent
\proof
Let $u(\theta )=\sum_{k\in\Z^n}\hat{u}_ke^{2\pi i k\cdot\theta}$, 
where obvoiusly one has 
$$
\hat{u}_k=\int_{\T^n}u(\theta )e^{-2\pi i k\cdot\theta}
d\theta\; . 
$$
Then the ${\cal C}^k$--norm can be equivalently 
given in terms of Fourier coefficients: for all $i\in\N$
there exists a positive constant $B_i$
such that 
$$
B_i^{-1}\sup_{k\in\Z^n}[
(1+|k|)^i|\hat{u}_k|]\le \Vert u\Vert_i\le B_i
\sup_{k\in\Z^n}[
(1+|k|)^{i+n+1}|\hat{u}_k|]\; .
\eqno(8.7)
$$
Comparing the Fourier coefficients of $u$
with those of $v$ one has 
$$
\hat{u}_k={\hat{v}_k\over 2\pi ik\cdot\mu}\; \;\; \forall
k\in\Z^n\setminus\{0\}\; ,
\eqno(8.8)
$$
The desired estimates are an easy consequence of (8.7), 
the assumption that $\mu$ is diophantine and of the 
elementary fact $\sum_{k\in\Z^n\setminus\{0\}} |k|^{-\delta}<+\infty$
for all $\delta> n$. \qed

\vskip .3 truecm\noindent
The fact that one needs $r$ more derivatives to bound the norms 
of $u$ in terms
of those of $v$ is what is called the ``loss of differentiability''. 
As we have already seen in Chapter 6 this 
is a typical phenomenon associated to small divisors. 
The 
analogue in the analytic case would be the necessary restriction of 
the domain to control the maximum norm of $u$ in terms of $v$
by means of Cauchy's estimates as we did in Section 6.3. 

In both cases (smooth and analytic) these are not artefacts of 
the methods used but a concrete manifestation of the unboundedness of 
the linear operator $D_\mu^{-1}$. The main consequence of this 
fact is that one cannot use Banach spaces techniques to 
study semilinear equations like $D_\mu u=v+\varepsilon f(u)$, 
where $\varepsilon$ is some small parameter. These semilinear 
equations are however typical of perturbation theory.

\line{}

\vskip .3 truecm\noindent
{\titsec 8.3 KAM Theory, Nekhoroshev Theorem, Arnol'd Diffusion}

\vskip .3 truecm\noindent
Despite Theorem 8.2, most results on quasi--integrable systems 
 have been obtained under the assumption of non--degeneracy 
(i.e. the hessian matrix of $h_{0}$ is non degenerate thus the 
frequency map $J\mapsto \nu_{0}(J)={\partial h_{0}\over\partial 
J}(J)\in {\Bbb R}^l$ is a local diffeomorphism) but accepting the 
fact that one cannot hope for integrability on open sets. 

The general picture 
is provided by KAM [Ar1,Ga, Bo, Yo1] and Nekhoroshev [Ne, Lo] 
theorems: if $\varepsilon$ is 
sufficiently small, most initial conditions (w.r.t. Lebesgue measure) 
lie on invariant $l$--dimensional lagrangian tori carrying 
quasiperiodic motions with Diophantine frequencies. The action 
variables corresponding to these orbits will remain $\epsilon$--close 
to their initial values for all times. 
The complement of this set is open and dense and it is connected if 
$l\ge 3$. It contains a connected ($l\ge 3$) web ${\cal R}$ of 
resonant zones corresponding to ${\Bbb Z}^l$--linearly dependent
frequencies: $\cup_{k\in {\Bbb Z}^l}\{J\in U\, , \nu_{0}(J)\cdot 
k = 0\} \times {\Bbb T}^l$. 
Motion along these resonances cannot be excluded (see [Ar2] for an 
explicit example), resulting in a 
variation of ${\cal O}(1)$ of the actions in a finite 
time\footnote{$^{2}$}{It is conjectured [AKN, p. 189] that 
generically quasi--integrable hamiltonians with more than two degrees 
of freedom are topologically unstable}. But if the hamiltonian is 
{\it analytic} and $h_{0}$ is {\it steep} (for example convex or quasi 
convex) then this variation is very slow: it takes a time at least 
${\cal O} \left(\exp\left({1\over\varepsilon^a}\right)\right)$
to change the actions of ${\cal O}(\varepsilon^b)$, where $a$ and $b$ 
are two positive constants. Moreover each invariant torus has a 
neighborhood filled in with trajectiories which remain close to it for 
an even longer time. Indeed, if $h_{0}$ is quasi--convex one can prove 
[GM] that all trajectories starting at a distance of order $\rho 
<\rho^{*}$ from a Diophantine $l$--torus of exponent $\tau$ will 
remain close to it for a time ${\cal 
O}\left(\exp\left(\exp\left({\rho^{*}\over\rho}\right)^{1/\tau 
+1}\right)\right)$. 

One of the consequences of KAM theorem [P\"o] is the existence, for 
sufficiently small values of $\varepsilon$, of a {\it Cantor set} 
$N_{\varepsilon}$ of values of the frequencies $\nu$ for which the 
Hamiltonian system (8.1) has smooth invariant tori with linear flow. 
Moreover there exists a homeomorphism 
$F_{\varepsilon}\,:N_{\varepsilon}\times {\Bbb T}^l\rightarrow 
U\times {\Bbb T}^l$ $\varepsilon$--close to the identity, Whitney 
smooth w.r.t. the first factor and analytic w.r.t. the second
(if the Hamiltonian (8.1) is analytic) which 
transforms Hamilton's equations into $\dot{\nu} = 0$, 
$\dot{\varphi} = \nu$. 
This foliation into invariant tori is thus parametrized 
over a Cantor set and hence nowhere dense. It exhibits 
the phenomenon of ``anisotropic differentiability'' since it is 
much more regular tangentially to these tori than transversally to them
(see also [BHS]). 
%%%%% fine capitolo 8 %%%%%
\vfill\eject
%\vskip 1. truecm
%%%%% capitolo 9: Nash-Moser %%%%%
\noindent
{\tit 9. The Inverse Function Theorem of Nash and Moser}

\vskip .3 truecm
\noindent
The Inverse Function Theorem for Banach spaces 
is one of the  extremely useful standard tools in the study of a 
variety of non--linear problems, ranging from the good position of 
the Cauchy problem for ordinary differential equations to non--linear 
elliptic equations. Unfortunately the ``loss of differentiability''
typical of small divisors problems prevents from its use 
(with some remarkable exceptions however, see Section 6.2 and [He2]). 
In the analytic case, Kolmogorov suggested the use of a modified 
Newton method to overcome this difficulty but in the differentiable 
case the need of an Inverse Function Theorem in Fr\'echet spaces 
has also other sources: its origin is the solution of the embedding 
problem for Riemannian manifolds by Nash [N]. Later Moser discovered 
how to adapt Kolmogorov's idea to the differentiable case creating 
a theory with a wide spectrum of applications [AG, Gr, H2, Ha, Ni, Ser, St,
SZ, Ze1]: to geometry, to the study of foliations and deformations 
of complex and CR structures, to free boundary problems, etc.. 
In all these cases a non--linear partial 
differential equation is solved using a rapidly convergent iterative 
algorithm introducing at each step of the iteration a smoothing of 
the approximate solution. 

In this Chapter we will follow the presentation of [Ha] very closely.

\line{}

\vskip .3 truecm
\noindent
{\titsec 9.1 Calculus in Fr\'echet Spaces} 

\vskip .3 truecm\noindent
\Proc{Definition 9.1}{A {\rm Fr\'echet space} is a locally convex 
topological vector space 
(lctvs)
which is complete, Hausdorff and metrizable.}

\vskip .3 truecm\noindent
{\bf Exercise 9.2} Show that a lctvs $X$ is Hausdorff if and only if $x\in X$, 
$\Vert x\Vert_i=0\;\forall i\in{\cal I}$ then $x=0$ (where $(\Vert\cdot
\Vert_i)_{i\in
{\cal I}}$ is the collection of seminorms giving the topology of $X$). 
Show that $X$ 
is metrizable if and only if ${\cal I}$ is countable. 

\vskip .3 truecm\noindent
{\bf Exercise 9.3} Show that $\R^\infty$ (space of all sequences 
of real numbers), 
${\cal C}^\infty ( M)$ (where $M$ is a smooth compact manifold), 
${\cal A}(\C )$
(entire functions) are Fr\'echet spaces (thus the exercise asks you
to define suitable seminorms). Show that ${\cal C}_0(\R )$ (continuous 
functions with compact support)
with the usual topology ($f_n\rightarrow f$ if and only if 
there exists a compact interval 
$I$ such that $\hbox{supp}\, f_n\subset I$ for all sufficiently large $n$, 
$\hbox{supp}\, f\subset I$ and $f_n$ converges uniformly to $f$ on $I$)
is a lctvs but it is not a Fr\'echet space since it is not metrizable. 

\vskip .3 truecm\noindent
{\bf Exercise 9.4} Prove that Hahn--Banach Theorem holds in Fr\'echet spaces:
if $X$ is a Fr\'echet space and $x$ is a non--zero vector in 
$X$ then there exists 
a continuous linear functional $l\, :X\rightarrow\R \,(\hbox{or}\, \C )$ 
such that 
$l(x)=1$. This allows to introduce quite straightforwardly 
$X$--valued analytic 
functions [Va]. A function $x\, :\Omega\rightarrow X$, where 
$\Omega$ is a region in $\C$, is analytic if and only if for all 
$l\in X^{*}$ the function $l\circ x$ is analytic. Show that this is 
equivalent to asking that, for all $z_{0}\in\Omega$, $x$ has a 
convergent power series expansion at $z_{0}$: $x(z)=\sum_{n=0}^\infty
(z-z_{0})^nx_{n}$. 

\vskip .3 truecm\noindent
{\bf Exercise 9.5} Extend the theory of Riemann's integration, including 
the fundamental theorem of calculus, to continuous $X$--valued functions 
on $[a,b]\subset \R$.

\vskip .3 truecm\noindent
\Proc{Definition 9.6}{Let $X,Y$ be two Fr\'echet spaces, $U\subset X$
be open, $f\, :U\rightarrow Y$ be continuous. The {\rm derivative}
of $f$ at $x\in U$ in the direction of $h\in X$ is
$$
Df(x)\cdot h := \lim_{t\rightarrow 0}
{f(x+th)-f(x)\over t}\; . \eqno(9.1)
$$
$f$ is ${\cal C}^1$ on $U$ if and only if $Df$ exists for all $x\in U$
and for all $h\in X$ {\rm and} $Df\, :U\times X\rightarrow Y$ is continuous.}

\vskip .3 truecm\noindent
{\bf Remark 9.7} In the case of Banach spaces this definition of ${\cal C}^1$
is weaker than the usual one.

\vskip .3 truecm\noindent
{\bf Exercise 9.8} Prove that the composition of ${\cal C}^1$ maps is 
${\cal C}^1$ and that the chain rule holds:
$D(g\circ f)(x)\cdot h=Dg(f(x))\cdot (Df(x)\cdot h)$. 

\vskip .3 truecm\noindent
{\bf Exercise 9.9} Define higher order derivatives and ${\cal C}^k$
maps between Fr\'echet spaces.

\vskip .3 truecm\noindent
{\bf Exercise 9.10} Let $f\, :{\cal C}^\infty([a,b])\rightarrow 
{\cal C}^\infty([a,b])$, $f(x)=P(x,x',\ldots ,x^{(n)})$,
where $P\in\R [X_0,\ldots X_n]$, is ${\cal C}^\infty$. Is there a 
nice formula for $Df(x)\cdot h$ ? [Hint: start from monomials 
like $(x^{(i)})^k$.]

\vskip .3 truecm\noindent
The following examples show why the extension 
of the inverse function theorem to Fr\'echet spaces is not a straightforward
generalization of the inverse function theorem in Banach spaces but needs
some extra assumption. 

The map $x\mapsto f(x)=\sin x$, where $x\in X=L^{2}([0,1 ])$, is 
of class ${\cal C}^{1}$ according to Definition 9.6. Its derivative 
$Df(0)=$identity but $f$ is not invertible: $f(0)=0$ and the 
functions $x_{n}(\xi )=\pi\chi_{[0,1/n]}(\xi )$, where $\chi_{[0,1/n]}$
denotes the characteristic function of the interval $[0,1/n]$, 
converge to $x=0$ but $f(x_{n})=0$ for all $n$. 

Another example is obtained taking $X={\cal C}^\infty ([-1,1])$ and 
considering the map
$f\, :X\rightarrow X$ defined as $f(x)(\xi )=
x(\xi )-\xi x(\xi )x'(\xi )$ for all $\xi\in [-1,1]$.
Then it is immediate to check that $f$ is smooth and 
$Df(x)\cdot h = h-\xi x' h-\xi x h'$, thus $f(0)=0$ and 
$Df(0)=$identity. But $f$
is not invertible: the sequence 
$x_n(\xi)={1\over n}
+{\xi^n\over n!}\rightarrow 0$ in ${\cal C}^\infty ([-1,1])$
as $n\rightarrow +\infty$ but one can check that it does not 
belong to $f({\cal C}^\infty ([-1,1]))$ for all $n\ge 1$. 
[Hint: use the fact that if $x\in {\cal C}^\infty ([-1,1])$
one can take its Taylor series at $0$ at any finite order and apply 
$f$. ]

An even more interesting counterexample (see [Ha] for details)
is the following:
let $M$ be a compact manifold, $X={\cal C}^\infty (M,TM)$ be the 
Fr\'echet space of smooth vector fields on $M$, $\hbox{Diff}^\infty
(M)$ be the group of smooth diffeomorphisms of $M$
(it is a Fr\'echet manifold, it's not very hard to figure out what 
this means, otherwise look in [Ha]). 
Then the usual exponential map 
$$
\eqalign{
\exp \, :  {\cal C}^\infty (M,TM) &\rightarrow \hbox{Diff}^\infty
(M)\cr
v &\mapsto \exp (v)\cr}
$$
clearly verifies $\exp (0)=\hbox{id}_M$ and $D\exp (0)=$identity, but 
the exponential map is not invertible in general. 
Note that this would have meant that any diffeomorphism extends to 
a one parameter flow. 

For example 
a diffeomorphism of $\S^1$ without fixed points is the exponential of 
a vector field only if it is conjugate to a rotation. But there 
exist [Yo3] diffeomorphisms of $\S^1$
arbitrarily close to the identity which are 
not conjugate to a rotation. 

What goes wrong in all these examples is that although the derivative
of the map is the identity at the origin it fails to be invertible 
at nearby points. Indeed in the second example above one
has $Df(1/n)\cdot\xi^k=\left(1-{k\over n}\right)\xi^k$, thus 
$Df(1/n)\xi^n=0$. 

Thus {\it one has to require the invertibility of $Df$ on a neighborhood 
explicitly} and this is usually difficult to be checked. 
In a Banach space (with the usual definition of derivative of a map 
instead of Definition 9.6) this is not needed. 

\line{}

\vskip .3 truecm\noindent
{\titsec 9.2 Tame Maps and Tame Spaces} 

\vskip .3 truecm\noindent
\Proc{Definition 9.11}{A {\rm graded Fr\'echet space} $X$ is a 
Fr\'echet space with a collection of seminorms $(\Vert\;
\Vert )_{n\in \N}$ which define
the topology and are increasing in strength }
$$
\Vert x\Vert_0\le\Vert x\Vert_1\le\Vert x\Vert_2\le\ldots\;\;
\forall x\in X\; .
$$

\vskip .3 truecm\noindent
\Proc{Definition 9.12}{Let $X,Y$ be graded Fr\'echet spaces, 
$U\subset X$ open, $P\, :U\rightarrow Y$ be a continuous 
map. $P$ is {\rm tame} if for all $x_0\in U$ there exists a 
neighborhood $V\subset U$ of $x_0$ and a non negative integer $r$
such that for all $i\in\N$ there exists $C_i>0$ such that
$$
\Vert P(x)\Vert_i\le C_i(1+\Vert x\Vert_{i+r})\;\;
\forall x\in V\; . \eqno(9.2)
$$
A ${\cal C}^k$ tame map is a ${\cal C}^k$ map $P$ such that 
$D^jP$ is tame for all $0\le j\le k$.}

\vskip .3 truecm\noindent
The most typical example of a tame operator between Fr\'echet 
spaces is given by nonlinear partial differential operators on 
compact manifolds. If $P\, :{\cal C}^\infty (M)\rightarrow 
{\cal C}^\infty (M)$ is a smooth function of $x\in {\cal C}^\infty (M)$
and its partial derivatives of degree at most $r$ then we say 
that the degree of $P$ is $r$ and this will be the 
``loss of differentiability'' in (9.2). 
The proof of this fact is given in [Ha] and uses 
Hadamard's inequalities for functions $x\in {\cal C}^\infty (M)$:
for all $n\in N$ and for all integer $k$ such that 
$0\le k\le n$ there exists $C_{k,n}>0$ such that 
$$
\Vert x\Vert_k\le C_{k,n}\Vert x\Vert_n^{k/n}\Vert x\Vert_0^{1-k/n}
\; \forall x\in {\cal C}^\infty (M) \;\hbox{and}\;
 . \eqno(9.3)
$$

\vskip .3 truecm\noindent
{\bf Exercise 9.13} Prove that the composition of two ${\cal C}^k$
tame maps is a $ {\cal C}^k$ tame map.

\vskip .3 truecm\noindent
\Proc{Definition 9.14}{A graded Fr\'echet space $X$ is {\rm tame}
if it admits {\rm smoothing operators}, i.e. a one--parameter family 
$S(t)\, :X\rightarrow X$, $t\in [1,+\infty )$, of continuous linear 
operators such that there exists a non negative integer $r$ 
and positive real constants $(C_{n,k})_{n,k \in \N}$
such that for all $x\in X$ and for all $t\in [1,+\infty )$ and
for all $k\in\{0,1,\ldots ,n\}$ one has}
$$
\eqalign{
\Vert S(t)x\Vert_n &\le C_{n,k}t^{n-k}\Vert x\Vert_k\cr
\Vert x- S(t)x\Vert_k &\le C_{k,n}t^{k-n}\Vert x\Vert_n\cr
}\eqno(9.4)
$$

\vskip .3 truecm\noindent
{\bf Exercise 9.15 (convolution with regularizing kernels)} 
Let $\psi\in {\cal C}_0^\infty (\R^n )$ and assume that 
$\psi\ge 0$, $\psi\equiv 1$ near $0$. 
Let $\varphi$ be the Fourier transform of $\psi$:
$\varphi (\xi ) = \int_{\R^n}\psi (\eta )e^{-2\pi i \xi \eta }d\eta$. 
Let $t\ge 1$, $\varphi_t (\xi )=t^n\varphi (t\xi )$. Define 
$S(t)f=\varphi_t\star f$, where $f\in {\cal C}^\infty
(\T^n )$. Show that $(S(t))_{t\ge 1}$ is a family 
of smoothing operators on ${\cal C}^\infty
(\T^n )$.

\vskip .3 truecm\noindent
{\bf Exercise 9.16} Prove that in a tame Fr\'echet space Hadamard's 
inequalities hold:
$$
\Vert x\Vert_l\le C(k,n)\Vert x\Vert_k^{1-\alpha}\Vert x\Vert_n^\alpha
\;\;\forall k\le l\le n\;\; , \;\; l=(1-\alpha )k+\alpha n\; .
$$
[Hint: use (9.4) with $t=\Vert x\Vert^{1/(n-k)}_{n}\Vert 
x\Vert^{-1/(n-k)}_{k}$.]

\line{}

\vskip .3 truecm\noindent
{\titsec 9.3 The Nash--Moser Theorem} 

\vskip .3 truecm\noindent
We can finally state Nash--Moser's [N,M]
implicit and inverse function theorems.

\vskip .3 truecm\noindent
\Proc{Theorem 9.17 (implicit function)}{
Let $X,Y,Z$ be three tame Fr\'echet spaces, $U\subset X\times Y$ open, 
$\Phi\,:U\rightarrow Z$ a tame ${\cal C}^r$ map, $2\le r\le \infty$.
Let $(x_0,y_0)\in U$. Assume that there exists a neighborhood
$V_0$ of $(x_0,y_0)$ and a continuous $z$--linear tame map 
$L\, :V_0\times Z\rightarrow Y$, $((x,y),z)\mapsto
L(x,y)\cdot z$, such that if $(x,y)\in V_0$ then $D_y\Phi (x,y)$ is invertible
with inverse $L(x,y)$. Then $x_0$ has a neighborhood $W$ on which 
$\Psi \in {\cal C}^r(W,Y)$ is defined and such that
$\Psi (x_0)=y_0$ and for all $x\in W$ one has 
$(x,\Psi (x))\in U$ and $\Phi (x,\Psi (x))=\Phi (x_0,y_0)$.}

\vskip .3 truecm\noindent
\Proc{Theorem 9.18 (inverse function)}{
Let $X,Y$ be two tame Fr\'echet spaces, $U\subset X$ open, 
$\Phi\,:U\rightarrow Y$ a tame ${\cal C}^r$ map, $2\le r\le \infty$.
Let $x_0\in U$, $y_0=\Phi (x_0)$. Assume that there exists a neighborhood
$V_0$ of $x_0$ and a continuous $y$--linear tame map 
$L\, :V_0\times Y\rightarrow X$, $(x,y)\mapsto
L(x)\cdot y$, such that if $x\in V_0$ then $D\Phi (x)$ is invertible
with inverse $L(x)$. Then $x_0$ has a neighborhood $V\subset V_0$ 
and $y_0$ has a neigborhood $W$ such that $\Phi\, : V\rightarrow
W$ is a tame ${\cal C}^r$ diffeomorphism.}

\vskip .3 truecm\noindent
{\bf Exercise 9.19} Show that the two previous theorems are equivalent.

\vskip .3 truecm\noindent
We refer the reader to [Ha] for the proofs of Theorems 9.17 and 9.18. 
The main idea of the proof is to use a modified Newton's method for 
finding the root of the equation $\Phi (x)=y$. It makes use of the 
smoothing operators $S(t)$ to guarantee convergence. Here we will 
content ourselves with a brief sketchy description of the argument.

Without loss of generality we can assume $x_0=y_0=0$. An algorithm 
for constructing a sequence $x_j\in X$ which will converge to a solution 
$x$ of  $\Phi (x)=y$ (for small enough $y$) is the following: 
fix a sequence $t_j=e^{(3/2)^j}$, so that $t_{j+1}=t_j^{3/2}$, and let
$$
\eqalign{
x_0 &= 0 \;\;\;\;\; (\hbox{initial guess}\, )\; , \cr
x_j &= \ldots \;\;\;\;\; (j\hbox{--th guess}\, )\; , \cr
z_j &= y-\Phi (x_j) \;\;\;\;\; (j\hbox{--th error}\, )\; , \cr
\Delta x_j &= S(t_j)L(x_j)z_j \;\;\;\;\; (j\hbox{--th correction}\, )
\; , \cr
x_{j+1} &= x_j+\Delta x_j\;\;\;\;\; (j+1\hbox{--th guess} \, )\; . \cr
}
$$
The idea to show convergence of this algorithm is the following: 
let 
$$
R(x,h)=\int_0^1D^2\Phi (x+th)(h,h)dt
$$
denote the quadratic integral remainder in Taylor's formula. Since 
$$
\Phi (x+h)=\Phi (x)+D\Phi (x)\cdot h +R(x,h)
$$
one gets 
$$
\eqalign{
z_{j+1} &= y-\Phi (x_{j+1}) = y-\Phi (x_j+\Delta x_j) \cr
&= z_j-D\Phi (x_j)S(t_j)L(x_j)z_j-R(x_j,\Delta x_j)\; .\cr
}
$$
Using the identity $z_j=D\Phi (x_j)L(x_j)z_j$ we find 
$$
z_{j+1} = D\Phi (x_j)[I-S(t_j)]L(x_j)z_j+R(x_j,\Delta x_j)\; . 
$$
The first term tends to zero very rapidly since $S(t_j)\rightarrow I$
as $j\rightarrow +\infty$ and the second term is quadratic. 

This short description of the idea of the proof makes also clear why one 
needs the assumption $\Phi$ at least of class ${\cal C}^2$
(in Banach spaces ${\cal C}^1$ is enough). 
%%%%% fine capitolo 9 %%%%%
\vfill\eject
%%%%% capitlo 10 %%%%%
\noindent
{\tit 10. From Nash--Moser's Theorem to KAM:
Normal Form of Vector Fields on the Torus}
\par
\vskip .3 truecm
\noindent
Following Herman we will prove in this Chapter a normal form theorem 
for vector fields on the torus which can be considered as the basic 
KAM theorem in higher dimension (without taking the symplectic 
structure into account). The proof will be an application of 
Nash--Moser's Theorem. For a proof of KAM theorem see, for example, 
[Bo].

Let $\hbox{Diff}^\infty (\T^l, 0)$ denote the group of ${\cal 
C}^\infty$ diffeomorphisms $f$ of the torus $\T^l$ homotopic to the 
identity and such that $f(0)=0$. 
This space can be identified to an open subset of the tame
Fr\'echet space ${\cal C}^\infty(\T^l,\R^l,0)=\{u\in{\cal C}^\infty
(\T^l ,\R^l)\, , \, u(0)=0\}$: $u$ corresponds to a diffeomorphism 
$f$ if and only if for all $\chi\in\T^l$ one has 
$\hbox{id}+\partial u(\chi )\in\hbox{GL}\,(l,\R)$. In this case one has 
$f=\hbox{id}_{\T^l}+u$. 

Since the tangent bundle of $\T^l$ is 
canonically isomorphic to $\T^l\times\R^l$ one can also identify
the space of ${\cal C}^\infty$ vector fields on the torus with ${\cal 
C}^\infty (\T^l,\R^l)$. 

Let $\mu\in\R^l$ be Diophantine with exponent $\tau$ and constant 
$\gamma$. We will denote $R_{\mu}$ the translation by $\mu$
on the torus $\T^l$: $R_{\mu}(\chi_{1},\ldots ,\chi_{l})=
(\chi_{1}+\mu_{1},\ldots ,\chi_{l}+\mu_{l})$. 

\vskip .3 truecm\noindent
{\bf Exercise 10.1} Show that the map 
$$
\eqalign{
I\, :\hbox{Diff}^\infty (\T^l, 0) &\rightarrow \hbox{Diff}^\infty (\T^l, 
0)\cr
f &\mapsto I(f)=f^{-1}\cr
}
$$
is a tame ${\cal C}^\infty$ map. Its derivative is 
$$
DI(f)\cdot h = -[(\partial f)^{-1}\cdot h]\circ f^{-1}\; . \eqno(10.1)
$$

\vskip .3 truecm\noindent
The following statements (and proof) are taken from ([Bo], pp. 
139--141) and [He5]. 

\vskip .3 truecm\noindent
{\bf Theorem 10.2}{\it Let $\mu\in\R^l$. The map
$$
\eqalign{
\Phi_{\mu}\, :\hbox{Diff}^\infty (\T^l, 0)\times\R^l &\rightarrow
\hbox{Diff}^\infty (\T^l)\; , \cr
(f,\nu ) &\mapsto R_{\nu}\circ f\circ R_{\mu}\circ f^{-1}\; , \cr}
\eqno(10.2)
$$
is a tame ${\cal C}^\infty$ map. Moreover, if $\mu$
is a diophantine\footnote{$^1$}{\rm In this situation $\mu$ is 
diophantine if there exist two constants $\gamma >0$ and 
$\tau \ge l$ such that $|\mu\cdot k+p|\ge \gamma |k|^{-\tau }$
for all $k\in\Z^l\setminus\{ 0\}$ and for all $p\in\Z$. }
vector then $\Phi_{\mu}$ is a tame 
${\cal C}^\infty$ local diffeomorphism near $f=
\hbox{id}_{\T^l}$, $\nu = 0$.}

\vskip .3 truecm\noindent
The meaning of the second part is that when $\mu$ is diophantine 
the diffeomorphisms of the torus $\T^l$ conjugate to the translation 
$R_{\mu}$ by a diffeomorphism close to the identity form a 
Fr\'echet submanifold of codimension $l$ of $\hbox{Diff}^\infty 
(\T^l)$ which is transverse in $\mu$ to the space $\R^l$ of the 
translations on the torus.

\vskip .3 truecm\noindent
{\bf Exercise 10.3} Guess the statement of for vector 
fields equivalent to Theorem 10.2.

\vskip .3 truecm\noindent
Here is the solution:

\vskip .3 truecm\noindent
{\bf Theorem 10.3}{\it Let $\mu\in\R^l$. The map
$$
\eqalign{
\Psi_{\mu}\, :\hbox{Diff}^\infty (\T^l, 0)\times\R^l &\rightarrow
{\cal C}^\infty (\T^l, \R^l)\; , \cr
(f,\nu ) &\mapsto \nu + f_{*}\mu = \nu +\partial f\circ f^{-1}\cdot
\mu\; , \cr}
\eqno(10.3)
$$
is a tame ${\cal C}^\infty$ map. Moreover, if $\mu$
is a diophantine vector (see (8.5) ) then $\Psi_{\mu}$ is a tame 
${\cal C}^\infty$ local diffeomorphism near $f=
\hbox{id}_{\T^l}$, $\nu = 0$.}

\vskip .3 truecm\noindent
\proof
First of all note that $\Psi_{\mu}(\hbox{id}_{\T^l}, 0)=\mu$. 
The first assertion is an immediate consequence of Exercises 9.13
and 10.1. Moreover, using (10.1), one easily checks that
$$
\eqalign{
D\Psi_{\mu}(f,\nu )\cdot (\Delta f,\Delta\nu ) &=
\Delta \nu + (\partial\Delta f)\circ f^{-1}\cdot\mu\cr
&+\partial^{2}f\circ f^{-1}\cdot (-(\partial f)^{-1}\circ f^{-1}\cdot
\Delta f\circ f^{-1},\mu )\cr
&= \Delta \nu + [(\partial\Delta f)\cdot\mu\cr
&+\partial^{2}f\cdot (-(\partial f)^{-1}\cdot
\Delta f,\mu )]\circ f^{-1}\cr}
\eqno(10.4)
$$
(we recall that here one has $\Delta f\in{\cal C}^\infty (\T^l,\R^l 
,0)$, $\Delta\nu\in\R^l$).

If one introduces $u$, writing $\Delta f=\partial f\cdot u$ then 
one gets
$$
(\partial\Delta f)\circ f^{-1}\cdot\mu=[\partial^{2}f\cdot u\cdot\mu
+\partial f\cdot\partial u\cdot\mu]\circ f^{-1}
$$
and
$$
\partial^{2}f\circ f^{-1}\cdot (-(\partial f)^{-1}\circ f^{-1}\cdot
\Delta f\circ f^{-1},\mu )=-[\partial^{2}f\cdot u\cdot\mu]\circ 
f^{-1}\;.
$$
Therefore
(10.4) simplifies considerably and becomes
$$
D\Psi_{\mu}(f,\nu )\cdot (\partial f\cdot u ,\Delta\nu )
= \Delta\nu + (\partial f\cdot \partial u\cdot\mu )\circ f^{-1}\; . 
\eqno(10.5)
$$
To prove the second assertion we will apply Theorem 9.18 to 
$\Phi = \Psi_{\mu}$ at the points $x_{0}=
(\hbox{id}_{\T^l}, 0)$ and $y_{0}=\mu$. We must just check that 
$D\Psi_{\mu}(f,\nu )$ is invertible {\it for all} 
$(f,\nu )$ in a neighborhood of $(\hbox{id}_{\T^l}, 0)$. This leads 
us to the equation 
$$
\Delta\nu + (\partial f\cdot\partial u \cdot\mu )\circ f^{-1} = w\; . 
\eqno(10.6)
$$
Composing on the 
right with $f$ and multiplying both sides by $(\partial f)^{-1}$ one gets 
$$
\mu\cdot\partial u = (\partial f)^{-1}\cdot [w\circ f-\Delta 
\nu]\; ,
\eqno(10.7)
$$
i.e.\ an equation of the form (8.4) with 
$v= (\partial f)^{-1}\cdot [w\circ f-\Delta 
\nu]$. This clarifies why one needs the term $\nu$ in the definition 
(10.3) of $\Psi_{\mu}$: indeed one fixes it so as to assure that 
$v\in{\cal C}^{0,\infty}(\T^l,\R^l)$, i.e.~it has zero 
average on the torus $\T^l$. One can also check that the map 
$(f,w)\mapsto\Delta\nu$ is tame. 

Proposition 8.7 allows to conclude since it shows that the map 
$(f,\nu )\mapsto u=D_{\mu}^{-1}v$ is tame. 
\qed
%%%%% fine capitolo 10 %%%%%
\vfill\eject
%%%%% appendici %%%%%
\centerline{\tit Appendices}

\vskip 2. truecm
\noindent
{\titsec  A1. Uniformization, Distorsion and Quasi--conformal maps}

\vskip .3 truecm
\noindent

In this appendix we recall some elementary and less elementary facts 
from the theory of conformal and quasi--conformal maps of one complex 
variable. 

\line{}

\vskip .3 truecm
\noindent
{\bf A1.1} A nonempty connected open set is called a {\it region}. 

\vskip .3 truecm\noindent
\Proc{Theorem A1.1 (The Maximum Principle)}{If $f(z)$ is analytic and 
non--constant in a region $\Omega$ of the complex plane $\C$, then 
its absolute value $|f(z)|$ has no maximum in $\Omega$.}

\vskip .3 truecm\noindent
\proof It is an easy consequence of the fact that non--constant 
analytic functions map open sets onto open sets. \qed

\vskip .3 truecm\noindent
The maximum principle implies that if $f$ is defined and continuous on 
a compact set $K$ and analytic in the interior of $K$ then the 
maximum of $|f(z)|$ on $K$ is assumed on the boundary of $K$. 
Another easy consequence is the following

\vskip .3 truecm
\noindent
{\bf Exercise A1.2 (Schwarz's Lemma, automorphisms of the disk)}
{\it Schwarz's Lem\-ma}: If $|f(z)|$ is analytic for 
$|z|<1$ and satisfies the conditions $|f(z)|\le 1$, $f(0)=0$, then 
$|f(z)|\le |z|$ and $|f'(0)|\le 1$. If $|f(z)|=|z|$ for some 
$z\not= 0$ or if $|f'(0)|=1$ then $f(z)=cz$ with $c\in \C$, 
$|c|=1$. {\it Automorphisms of the disk}: Show that if $f\, : \D
\rightarrow\D$ is an automorphism of the disk and $f(0)=0$ then 
$|f'(0)|=1$ and $f$ is a rotation. Deduce from this that 
the group of automorphisms of the unit disk $\D$ 
is 
$$
\hbox{Aut}\,(\D ) = \{ z\mapsto T(z)={az+b\over 
\overline{b}z+\overline{a}}\, , \; a,b\in \C\, , \; 
|a|^{2}-|b|^{2}=1\} ]\; .
$$

\line{}

\vskip .3 truecm\noindent
{\bf A1.2} A mapping $f$ of a region $\Omega$ into $\C$ is called {\it 
conformal} if it is holomorphic and injective. Such maps are also 
called {\it univalent}. Since an analytic map is injective if and only if 
$f'\not= 0$ if $f$ is univalent in $\Omega$ its derivative never 
vanishes, i.e. it has no critical points inside $\Omega$. 
The most important result of the theory of conformal maps is certainly 
the 

\vskip .3 truecm
\noindent
\Proc{Theorem A1.3 (Riemann Mapping Theorem)}{Given any simply 
connected region $\Omega$ which is not the whole plane, and a point 
$z_{0}\in \Omega$ there exists a unique conformal map (the Riemann 
map)
$f\,:\D\rightarrow \Omega$ such that $f$ is onto, $f(0)=z_{0}$ and 
$f'(0)>0$.}

\vskip .3 truecm\noindent
{\bf Exercise A1.4} Drop the requirement $f'(0)>0$. Then $f$ is not 
unique but the number $|f'(0)|$ does not depend on $f$. It is called
[Ah1] 
the {\it conformal capacity} of $\Omega$ with respect to $z_{0}$ and it will 
be denoted $C(\Omega ,z_{0})$. [Hint: use Exercise A1.3]

\vskip .3 truecm\noindent
One should not think to an arbitrary simply connected region as the 
``potato'' of PDEs but as a rather irregular object. For 
example the boundary needs not to be locally connected. 

\line{}

\vskip .3 truecm\noindent
{\bf A1.3}
Assume that  the 
region $\Omega$ is bounded and $\partial \Omega$ is a closed Jordan 
curve (i.e. $\partial \Omega = \gamma ([0,1])$, where 
$\gamma\,:\, [0,1]\rightarrow \C$ is continuous, $\gamma (0)=\gamma 
(1)$ and $\gamma (t_{1})=\gamma (t_{2})$ if and only if $t_{1}=t_{2}$
or $t_{1}=0$, $t_{2}=1$). In this case the Riemann map $f\, : 
\D\rightarrow\Omega$ has a nice boundary behaviour: 

\vskip .3 truecm\noindent
\Proc{Theorem A1.5 (Caratheodory)}{ A Riemann map $f\, : 
\D\rightarrow\Omega$ extends to a homeomorphism of $\overline{\D}$ 
onto $\overline{\Omega}$ if and only if $\partial\Omega$ is a closed 
Jordan curve.}

\vskip .3 truecm
\noindent
The topic of the boundary behaviour of conformal maps is very rich and 
it's an active research area: we refer to [Po] for more informations 
and references. We will only need two other results: the first extends 
Caratheodory's theorem dropping the assumption that the restriction of 
$f$ to the boundary of the disk is injective.

\vskip .3 truecm
\noindent
\Proc{Theorem A1.6}{Let $f\, : \D\rightarrow\Omega$ be a Riemann map. 
The following four conditions are equivalent:
\item{(i)} $f$ has a continuous extension to $\overline{\D}$; 
\item{(ii)} $\partial\Omega$ is a continuous curve, i.e. 
$\partial\Omega = \{\varphi (\zeta )\, , \zeta\in \T\}$ with 
$\varphi$ continuous; 
\item{(iii)} $\partial\Omega$ is locally connected; 
\item{(iv)} $\C\setminus\Omega$ is locally connected.}

\vskip .3 truecm
\noindent 
Our second result, due to Fatou, applies to all $f\, : \D\rightarrow\C$
holomorphic {\it and bounded} (thus we are dropping the assumption of 
$f$ being injective in $\D$). 

\vskip .3 truecm\noindent
\Proc{Theorem A1.7 (Fatou)}{Let $f\, : \D\rightarrow\C$ be holomorphic 
and bounded. Then $f$ has a non--tangential limit at almost all points 
$\zeta\in \T=\partial\D$. Moreover if $f$ is not identically zero then 
$\varphi (\zeta ) =\lim_{r\rightarrow 1-}f(r\zeta )$ (which belongs to 
$L^\infty (\T )$) is not zero almost everywhere.}

\line{}

\vskip .3 truecm\noindent
{\bf A1.4} Another fundamental result is the celebrated 

\vskip .3 truecm\noindent 
\Proc{Theorem A1.8 (Uniformization Theorem)}{The only simply connected 
Riemann surfaces, up to biholomorphic equivalence, are the Riemann 
sphere $\overline{\C}=\C\cup\{\infty\}$, the complex plane $\C$ and 
the unit disk $\D$.}

\vskip .3 truecm\noindent
{\bf Exercise A1.9} Prove that the group of automorphisms of the 
Rie\-mann sphere is the group $\hbox{PGL}\, (2,\C)$ acting by homographies: 
if $g=\abcd\in \hbox{PGL}\, (2,\C)$ then $z\mapsto g\cdot z = 
{az+b\over cz+d}$. The group of automorphisms of the complex plane is 
simply the affine group. 

\line{}

\vskip .3 truecm\noindent
{\bf A1.5} One can consider univalent functions $f$ on regions of $\Cbar$
with values in $\Cbar$: in this case  $f$ must be meromorphic and 
injective. Here are some elementary properties: 

\vskip .3 truecm\noindent
{\bf Exercise A1.10} \item{(a)} If $f$ is univalent on a region 
$\Omega\subset \Cbar$ then $f$ is analytic except for at most a single
simple pole and $f'$ never vanishes. 
\item{(b)} If $f\, : \Omega\rightarrow \Omega'$ is onto and univalent 
then $f^{-1}\, : \Omega'\rightarrow\Omega$ is also univalent.
\item{(c)} A univalent map is a homeomorphism.
\item{(d)} A univalent map preserves angles between curves and their 
orientation (that's why they're called conformal!).
\item{(e)} The composition of univalent maps is univalent; $f$ is 
univalent if and only if $1/f$ is univalent.

\vskip .3 truecm\noindent
{\bf Exercise A1.11} Prove that if $f\, : \Omega\rightarrow\Cbar$ is 
univalent and $A\subset \Omega$ is measurable then 
$$
\hbox{Area}\,(f(A)) = \int_{A}|f'(x+iy)|^{2}dxdy\; . 
$$

\vskip .3 truecm\noindent
{\bf Exercise A1.12 (The Area Formula)} Show that if $f\, 
:\D\rightarrow\C$ is univalent, letting $f(z)=\sum_{n=0}^\infty 
f_{n}z^n$, one has 
$$
\hbox{Area}\,(f(\D)) = \pi \sum_{n=1}^\infty n|f_{n}|^{2}\; . 
$$
[Hint: First consider the disk $\D_{r}$ of radius $r<1$. If $f=u+iv$ 
then $\hbox{Area}\,(f(\D)) = \int_{\partial\D_{r}}udv={i\over 
2}\int_{\partial\D_{r}}fd\overline{f}$. Then let $r\rightarrow 1-$.]

\vskip .3 truecm\noindent
Let $S_{1}$ denote the collection of functions $f$ univalent in $\D$ 
and such that $f(0)=0$, $f'(0)=1$, thus $f(z)=z+f_{2}z^{2}+\ldots$. 
With $\Sigma_{1}$ we will denote all functions $g(\zeta )=\zeta +g_{0}+
g_{1}\zeta^{-1}+\ldots $ univalent in the outer disk 
$\E =\{\zeta\in \Cbar\, , |\zeta |>1\}$. Clearly if $f\in S_{1}$ then 
$g(\zeta )=1/f(\zeta^{-1})$ belongs to $\Sigma_{1}$ and omits $0$. 
Conversely, if $g\in \Sigma_{1}$ and $g(\zeta)\not= 0$ for all 
$\zeta\in\E$ then $f(z)=1/g(z^{-1})$ belongs to $S_{1}$. 
One of the most important results on univalent functions is the 
object of the following exercise:

\vskip .3 truecm\noindent
{\bf Exercise A1.13 (Area Theorem)} If $g\in\Sigma_{1}$ then 
$|g_{1}|\le \sum_{n=1}^\infty n|g_{n}|^{2}\le 1$. Prove that equality 
holds if and only if $g(\zeta )=\zeta+g_{0}+g_{1}\zeta^{-1}$ with 
$|g_{1}|=1$. [Hint: show that area$(\C\setminus g(\E )) = \pi 
(1-\sum_{n=1}^\infty n|g_{n}|^{2} )$.]

\vskip .3 truecm\noindent
One of the main consequences is the following 
apparently innocent bound: if $f\in S_{1}$ then 
$$
|f_{2}|\le 2\; , \eqno(A1.1)
$$
as one can easily check applying the Area Theorem to 
$g(\zeta )=\sqrt{f(\zeta^{-2})}$. However this estimate will 
have many important consequences as we will see soon. 
A (much harder and for a long time conjectural) result is the 
celebrated [DeB]

\vskip .3 truecm\noindent
\Proc{Theorem A1.14 (Bieberbach--De Branges)}{ If $f\in S_{1}$ then 
$|f_{n}|\le n$ for all $n$.}

\vskip .3 truecm\noindent
Using (A1.1) one can easily show that the image of $\D$ through a 
univalent map cannot be too small:

\vskip .3 truecm\noindent
\Proc{Theorem A1.15 (Koebe $1/4$--Theorem)}{If $f\in S_{1}$ then 
$\D_{1/4}\subset f(\D)$.}

\vskip .3 truecm\noindent 
\proof 
Let $w\in \D$ and assume $w\notin f(\D)$. Then 
$$
\tilde{f}(z)={wf(z)\over w-f(z)}=z+(f_{2}+w^{-1})z^{2}+\ldots
$$
belongs to $S_{1}$. Applying (A1.1) to both $f$ and $\tilde{f}$ one gets
$|w|^{-1}\le |f_{2}|+|f_{2}+w^{-1}|\le 4$. \qed

\vskip .3 truecm\noindent
The {\it Koebe function} $f(z)=z(1-z)^{-2}=\sum_{n=1}nz^n$ maps the 
unit disk $\D$ conformally onto  
$\C\setminus (-\infty , -1/4)$. Therefore 
Bieberbach--De Branges' Theorem  and Koebe $1/4$--Theorem are optimal. 
If $f$ is univalent and analytic in $\D$, given any $z_{0}\in \D$
the {\it Koebe transform } of $f$ at $z_{0}$
$$
\eqalign{
K_{z_{0},f}(z) &= {f\left({z+z_{0}\over 1+\overline{z_{0}}z}\right)-f(z_{0})
\over (1-|z_{0}|)^{2}f'(z_{0})} \cr
&= 
z+\left[ {(1-|z_{0}|)^{2}f''(z_{0})\over 
2f'(z_{0})}-\overline{z_{0}}\right]z^{2}+\ldots\cr}
\eqno(A1.2)
$$
belongs to $S_{1}$. This is a very useful tool in order to transfer 
the information at $0$ to information at any point of the disk. 
Applying systematically this idea, from (A1.1) one deduces the 
following important distorsion estimates:

\vskip .3 truecm\noindent
{\bf Exercise A1.16 (Koebe distortion theorems)} If $f$ maps $\D$
conformally into $\C$ then $\forall z\in \D$ one has:
$$
\eqalignno{
\left| (1-|z|^{2}){f''(z)\over f'(z)}-2\overline{z}\right| &\le 4
\; , \; &(A1.3)\cr
|f'(0)|{|z|\over (1+|z|)^{2}} \le |f(z)-f(0)| &\le 
|f'(0)|{|z|\over (1-|z|)^{2}}\; , \;  &(A1.4)\cr
|f'(0)|{1-|z|\over (1+|z|)^{3}} \le |f'(z)| &\le 
|f'(0)|{1+|z|\over (1-|z|)^{3}}\; , \;  &(A1.5)\cr
{1\over 4} (1-|z|^{2})|f'(z)| \le \hbox{dist}\, (f(z),\partial 
f(\D ))&\le (1-|z|^{2})|f'(z)|\; . \; &(A1.6)\cr}
$$

\vskip .3 truecm\noindent
These estimates show that the growth of $f$ as $z$ approaches 
$\partial\D$ cannot be faster than $(1-r)^{-2}$, where $r=|z|$. The 
next theorem (see [Po]) for a proof) shows that the average growth is 
much lower than $(1-r)^{-2}$. 

\vskip .3 truecm\noindent 
\Proc{Theorem A1.17}{Let $f$ map $\D$ conformally into $\C$. Then 
$f(\zeta )=\lim_{r\rightarrow 1} f(r\zeta )\not=\infty$ exists for 
almost all $\zeta\in\T=\partial\D$ and for $0\le r<1$ one has}
$$
{1\over 2\pi}\int_{0}^{2\pi}|f(re^{i\theta}-f(0)|^{2/5}dt\le 
5|f'(0)|^{2/5}\; . \eqno(A1.7)
$$

\line{}

\vskip .3 truecm\noindent
{\bf A1.6} We conclude our brief introduction to univalent functions with the 
proof of a fundamental property of $S_{1}$. In order to do this we 
recall ([Re], p. 163)

\vskip .3 truecm\noindent 
\Proc{Lemma A1.18 (Hurwitz)}{If a sequence $(f_{n})_{n\in\N}$ of 
functions holomorphic in a region $\Omega\subset\C$ converges 
uniformly on compact subsets of $\Omega$ to a non--constant 
holomorphic function $f\, :\Omega\rightarrow \C$ then the following 
statements hold:
\item{(a)} if all the images $f_{n}(\Omega )$ are contained in a 
fixed set $A$ then $f(\Omega )\subset A$;
\item{(b)} if all the maps $f_{n}\,:\Omega\rightarrow \C$ are 
injective then so is $f\,:\Omega\rightarrow \C$; 
\item{(c)} if all the maps $f_{n}\,:\Omega\rightarrow \C$ are
locally biholomorphic, then so is $f\,:\Omega\rightarrow \C$.}

\vskip .3 truecm\noindent
This is the ingredient we missed for the proof of the following 

\vskip .3 truecm\noindent
\Proc{Theorem A1.19}{$S_{1}$ endowed with the topology of uniform 
convergence on compact subsets of $\D$ is a compact topological 
space.}

\smallskip\noindent
\proof
Any sequence $(f_{n})_{n\in\N}\subset S_{1}$ is equicontinuous and 
uniformly bounded on compact subsets of $\D$ by Koebe distortion 
theorems. Limit functions are in $S_{1}$ because they are univalent 
by Hurwitz's lemma and the normalisation $|f_{n}'(0)|$ for all 
$n\in\N$. \qed

\vskip .3 truecm\noindent
{\bf Exercise A1.20} 
Prove that Theorem A1.19 is equivalent to the following (see [Mc]): 
the space of {\it all} univalent maps $f\, : \D\rightarrow
\overline{\C}$ is compact up to post--composition with 
automorphisms of $\overline{\C}$. This precisely means 
that any sequence of univalent maps contains a subsequence 
$f_{n}\, : \D\rightarrow
\overline{\C}$ such that $M_{n}\circ f_{n}$ converges to a univalent 
map $f$, uniformly on compact subsets of $\D$, for some sequence of 
M\"obius transforms $M_{n}\in\hbox{PGL}\,(2,\overline{\C})$. 

\line{}

\vskip .3 truecm\noindent
{\bf A1.7} Let $f\,:\C\rightarrow\C$ be a ${\cal C}^1$ 
orientation--preserving diffeomorphism. Then given any point $z_0
\in\C$ one has 
$$
f(z)=f(z_0)+f_z(z_0)(z-z_0)+f_{\bar{z}}(z_0)(\bar{z}-\bar{z}_0)
+\hbox{o}\,(|z-z_0|)\; , \eqno(A1.8)
$$
where 
$$
f_z={1\over 2}\left({\partial f\over\partial x}-i
{\partial f\over\partial y}\right)\; , \;\;\;
f_{\bar{z}}={1\over 2}\left({\partial f\over\partial x}+i
{\partial f\over\partial y}\right)\; , \;\;\;
(z=x+iy)\; . \eqno(A1.9)
$$
Note that if $f$ is analytic in $z_0$ then $f_{\bar{z}}(z_0)=0$
(Cauchy--Riemann). The Jacobian determinant of $f$ is 
$J=|f_z|^2-|f_{\bar{z}}|^2$. Since $f$ is orientation--preserving 
one has $J>0$, thus $|f_z|>|f_{\bar{z}}|$.

\vskip .3 truecm\noindent
\Proc{Definition A1.21}{The {\rm dilatation} of $f$ in $z_0$ is }
$$
D_f(z_0) := {|f_z(z_0)|+|f_{\bar{z}}(z_0)|\over
|f_z(z_0)|-|f_{\bar{z}}(z_0)|}\ge 1\; . \eqno(A1.10)
$$

\vskip .3 truecm\noindent
Note that if $f$ is conformal then $D_f=1$.

Here is a geometric interpretation of the meaning of 
the dilatation: the differential $df(z_0)$ maps a circle
in the tangent space $T_{z_0}\C$ into an ellipse in 
$T_{f(z_0)}\C$. The dilatation measures the distorsion 
since it is the ratio of the major semiaxis amd the minor
semiaxis. Indeed applying  (A1.8) to an infinitesimal
circle $\Delta z=\varepsilon e^{i\theta}$ centered at $z_{0}$
one finds an infinitesimal ellipse centered at $f(z_{0})$
with  major semiaxis 
%% 
 % is mapped into 
 % $|f_{\bar{z}}(z_0)\overline{\Delta z}+f_z(z_0)\Delta z|=\delta$. 
 % However 
 % $$
 % |f_{\bar{z}}(z_0)\overline{\Delta z}+f_z(z_0)\Delta z|^2
 % = |\Delta z|^2(|f_z(z_0)|^2+|f_{\bar{z}}(z_0)|^2)+
 % 2|f_z(z_0)||f_{\bar{z}}(z_0)|\Delta z\cdot\overline{\Delta z}\; ,
 % $$
 % where 
 % $\Delta z\cdot\overline{\Delta z}=(\Delta x,\Delta y)\cdot (\Delta x,
 % -\Delta y)=(\Delta x)^2-(\Delta y)^2$. Therefore one has 
 % $$
 % \delta=(\Delta x)^2[|f_z|^2+|f_{\bar{z}}|^2+2|f_z||f_{\bar{z}}|]
 % +(\Delta y)^2[|f_z|^2+|f_{\bar{z}}|^2-2|f_z||f_{\bar{z}}|]
 % $$
 % and this is an infinitesimal ellipse with
 %%
$[|f_z|-|f_{\bar{z}}|]^{-1}\varepsilon$ and minor semiaxis 
$[|f_z|+|f_{\bar{z}}|]^{-1}\varepsilon$.

The {\it maximal dilatation} of $f$ on $\C$ is 
$D_f=\sup_{z\in\C}D_f(z)$. If $D_f<+\infty$, let 
$\kappa_f=(D_f-1)/(D_f+1)$. Then one has 
${|f_{\bar{z}}|\over |f_z|}\le \kappa_f<1$. 

\vskip .3 truecm\noindent
\Proc{Definition A1.22}{$f$ is {\rm quasiconformal}
if $D_f<+\infty$, i.e. $\kappa_f<1$.}

\vskip .3 truecm\noindent
Clearly if $f$ is conformal then $D_f=1$, $\kappa_f=0$.

We want now to extend the notion of quasiconformal map to 
homeomorphisms. We will follow the geometric approach 
outlined in [Ah2].

\vskip .3 truecm\noindent
{\bf Exercise A1.23} Given two rectangles $R_1$ and $R_2$
respectively with sides $a_1\le b_1$ and $a_2\le b_2$, 
show that there exists a conformal map of $R_1$ onto $R_2$
which maps vertices on vertices if and only if 
${a_1\over b_1}={a_2\over b_2}$.

\vskip .3 truecm\noindent
\Proc{Definition A1.24}{A {\rm quadrilateral}
$Q(z_1,z_2,z_3,z_4)$ is a Jordan domain in $\C$
with four distinguished boundary points $z_1,z_2,z_3,z_4$.
Its {\rm modulus} $M(Q)$ is the ratio $a/b$ of the lengths
$a<b$ of the sides of any rectangle $R$ which is the 
conformal image of $Q$ and whose vertices are image 
of the distinguished points.}

\vskip .3 truecm\noindent
Note that the modulus of a quadrilateral is a conformal invariant.
Thus one can use its variation under a homeomorphism to 
measure the lack of conformality of a map.

\vskip .3 truecm\noindent
\Proc{Definition A1.25}{Let $f\,:\C\rightarrow\C$
be an orientation--preserving homeomorphism. Its 
{\rm maximal dilatation} is $D_f:=\sup_{Q\subset\C\,,\,
Q\,\hbox{quadrilateral}}{M(f(Q))\over M(Q)}$.
Let $\kappa_f=(D_f-1)/(D_f+1)$. $f$ is a {\rm quasiconformal
homeomorphism} if $D_f,+\infty$, i.e. $\kappa_f<1$.}

\vskip .3 truecm\noindent
{\bf Exercise A1.26} Prove that $f$ is conformal if and only if 
$D_f=1$.

\vskip .3 truecm\noindent
\Proc{Theorem A1.27}{A quasiconformal homeomorphism $f$ with 
maximal dilatation $D_f$ is almost everywhere differentiable 
and at each point $z_0$ where $f$ is differentiable one has }
$$
{|f_{\bar{z}}(z_0)|\over |f_z(z_0)|}\le\kappa_f\; . 
$$

\vskip .3 truecm\noindent
Perhaps the most useful result in the theory of quasiconformal 
maps is the following existence theorem also known as 
Measurable Riemann Mapping Theorem [Ah2, p.98]

\vskip .3 truecm\noindent
\Proc{Theorem A1.28}{Let $\mu$ be a complex--valued measurable 
function with $\Vert\mu\Vert_\infty <1$. There exists a quasiconformal 
mapping $f$ such that 
$$
f_{\bar{z}}=\mu (z)f_z\;\hbox{almost everywhere}\eqno(A1.11)
$$
and $f$ leaves the points $0,1,\infty$ fixed. }

\vskip .3 truecm\noindent
The equation (A1.11) is also known as {\it Beltrami equation}.

Quasiconformal maps have been introduced in the subject of holomorphic 
dynamics by Dennis Sullivan and Adrien Douady and have rapidly become a 
standard tool. What we will need in Chapter 3 is the following

\vskip .3 truecm\noindent
\Proc{Theorem A1.29 (Douady--Hubbard: stability of the quadratic 
polynomial)}{Let $P_\lambda (z)=\lambda\left(z-{z^2\over 2}\right)$
and let $F(z)=P_\lambda (z)+\psi (z)$ where $\psi$ is 
holomorphic and bounded in the disk $\D_3$, $\psi (z)=
\sum_{n=2}^\infty\psi_nz^n$ (i.e. $\psi(0)=\psi'(0)=0$).
Assume that $\sup_{z\in\D_3}|\psi (z)|<10^{-2}$.
Then there exists a quasiconformal homeomorphism $h$
such that on the disk $\D_2$ one has $h^{-1}Fh=P_\lambda$. If 
$\psi$ is small enough then $h$ is near the identity in the ${\cal C}^0$
topology.}

\vskip .5 truecm
\vskip 2. truecm\noindent
\noindent
{\titsec A2. Continued Fractions}

\vskip .3 truecm
\noindent
In this appendix we recall some elementary facts on standard real continued 
fractions (we refer to [MMY], and references therein, for more general 
continued fractions). 

We will consider the iteration of the Gauss map
$$
A : (0,1) \mapsto [0,1] \; , \eqno(A2.1)
$$
defined by
$$
A (x) = {\ 1\ \over x} -
\left[\ {\ 1\ \over x}\ \right]\;. \eqno(A2.2)
$$
$A$ is piecewise analytic with branches 
$$
A(x)=x^{-1}-n \; \; \hbox{if}\; {1\over n+1}<x\le{1\over n}\; ,
n\ge 1\; .
$$

\vskip .3 truecm\noindent
{\bf Exercise A2.1} Prove that $A^*(\rho (x)dx)=\rho (x)dx$
where  $\rho (x)=[(1+x)\log 2]^{-1}$, i.e. $\rho$
is an invariant probability density for the Gauss map.

\vskip .3 truecm\noindent
Let
$$
G = {\sqrt{5}+1\over 2}\; , \; 
g = G^{-1} = {\sqrt{5}-1\over 2} \; . \;\;\;
$$
To each $x \in \Bbb R \setminus \Bbb Q$ we associate a continued fraction
expansion by iterating $A$ as follows. Let
$$
\eqalign{x_0 & = x - [x] \; , \cr
            a_0 & = [x] \; , \cr} \eqno(A2.3)
$$
then one obviously has
$x = a_0 +  x_0$. 
We now define inductively for all $n \ge 0$
$$
\eqalign{x_{n+1} & = A(x_n) \; , \cr
            a_{n+1} & = \left[ {1 \over x_n} \right]  \ge 1\; ,\cr}
            \eqno(A2.4)
$$
thus
$$
x_{n}^{-1} = a_{n+1} + x_{n+1} \; . \eqno(A2.5)
$$
Therefore we have
$$
x=a_0 + x_0=a_0+{1\over a_1 + 
     x_1}= \ldots =a_0 + \displaystyle{1 \over a_1
     + \displaystyle{1  \over a_2 + \ddots +
     \displaystyle{1 \over a_n + x_n}}}\; , 
     \eqno(A2.6)
$$
and we will write
$$
x=[a_0,a_1,\ldots ,a_n,
     \ldots] \;. \eqno(A2.7)
$$
The nth-convergent is defined by
$$
{p_n \over q_n} = [a_0,a_1,\ldots ,
                     a_n] =
                     a_0 + \displaystyle{1 \over a_1
     + \displaystyle{1  \over a_2 + \ddots +
     \displaystyle{1 \over a_n }}}
     \;. \eqno(A2.8)
$$

\vskip .3 truecm
\noindent
{\bf Exercise A2.2} Show that  the numerators $p_n$ and denominators
$q_n$ are recursively determined by
$$
p_{-1}=q_{-2}=1 \;\;,\;\;\;p_{-2}=q_{-1}=0 \;\;,\eqno(A2.9)
$$
and for all $n \ge 0$ one has 
$$
\eqalign{ p_n &= a_n p_{n-1} +  p_{n-2} \; , \cr
            q_n &= a_n q_{n-1} + q_{n-2} \; . \cr}
           \eqno(A2.10)
$$

\vskip .3 truecm\noindent
{\bf Exercise A2.3} Show that for all $n\ge 0$ one has 
$$
\eqalignno{x &= {p_n + p_{n-1}  x_n \over q_n + q_{n-1}
       x_n } \; , &(A2.11) \cr
              x_n &= - {q_n x -p_n \over q_{n-1} x - p_{n-1}}
              \; , &(A2.12) \cr
              q_n p_{n-1} - p_n q_{n-1} &= (-1)^n \;\; .  &(A2.13) \cr}
$$

\vskip .3 truecm
\noindent
Note that $q_{n+1}>q_n>0$ and that the sequence of the numerators 
$p_n$ has the same constant sign of $x$. Equation (A2.13) implies also that 
for all $k\ge 0$ and for all $x\in\R\setminus\Q$ one has 
${p_{2k}\over q_{2k}}<x<{p_{2k+1}\over q_{2k+1}}$.

Let
$$
\beta_n = \Pi_{i=0}^n x_i = (-1)^n (q_n x - p_n)\quad\hbox{for\ }n\ge 0,\quad
  \hbox{and\ }\be_{-1}=1\;\; .\eqno(A2.14)
$$
Then $x_n=\beta_n\beta_{n-1}^{-1}$ and $\beta_{n-2}=a_n\beta_{n-1}+
\beta_{n}$. 

\vskip .3 truecm\noindent
\Proc{Proposition A2.4} {For all
$x \in \Bbb R \setminus \Bbb Q$ and for all $n \ge 1$ one has
\item{(i)}\qquad $ \left|q_n x - p_n\right|
={\dst 1\over\dst q_{n+1}+q_nx_{n+1}}$,
so that ${\dst 1\over\dst 2}<\beta_nq_{n+1}<1$\ ;
\item{(ii)}$\beta_n\le g^n$ and $q_n\ge{\dst1
\over\dst 2}G^{n-1}$\ .
}
\par
\vskip .3 truecm\noindent
\proof
Using (A2.11) one has 
$$
\eqalign{
|q_nx-p_n| &= \left|q_n{p_{n+1}+p_nx_{n+1}\over q_{n+1}+q_nx_{n+1}}
-p_n\right| = {|q_np_{n+1}-p_nq_{n+1}|\over q_{n+1}+q_nx_{n+1}}\cr
&= {1\over  q_{n+1}+q_nx_{n+1}}\cr
}
$$
by (A2.13). This proves (i). 

Let us now consider $\beta_n=x_0x_1\ldots x_n$. If $x_k\ge g$ for some 
$k\in\{0,1,\ldots ,n-1\}$, then, letting $m=x_k^{-1}-x_{k+1}\ge 1$, 
$$
x_kx_{k+1}=1-mx_k\le 1-x_k\le 1-g=g^2\; .
$$
This proves (ii). \qed

\vskip .3 truecm
\noindent 
{\bf Remark A2.5} Note that from {\it (ii)} it follows that 
$\sum_{k=0}^\infty {\log q_{k}\over q_{k}}$ and 
$\sum_{k=0}^\infty {1\over q_{k}}$
are always convergent and their sum is uniformly bounded. 

\vskip .3 truecm
\noindent
For all integers  $k\ge 1$, the iteration of the Gauss map $k$ times leads
to the following partition of 
$(0,1)$; $\sqcup_{a_{1},\ldots ,a_{k}}
I(a_{1},\ldots ,a_{k})$, where $a_{i}\in {\Bbb N}$, 
$i=1,\ldots ,k$, and 
$$
I(a_{1},\ldots ,a_{k}) = \cases{
\left( {p_{k}\over q_{k}},{p_{k}+p_{k-1}\over q_{k}+q_{k-1}}\right)
&if $k$ is even\cr
\left( {p_{k}+p_{k-1}\over q_{k}+q_{k-1}}, {p_{k}\over q_{k}}\right)
&if $k$ is odd\cr}
$$
is the branch of $A^k$ determined by the fact that all points 
$x\in I(a_{1},\ldots ,a_{k}) $ have the 
first $k+1$ partial quotients exactly equal to $\{0,a_{1},\ldots ,a_{k}\}$. 
Thus 
$$
I(a_{1},\ldots ,a_{k}) = \left\{ x\in (0,1)\,\mid\; x=
{p_{k}+p_{k-1}y\over q_{k}+q_{k-1}y}\; , \; y\in (0,1)\right\}\; . 
$$
Note that ${dx\over dy}= {(-1)^{k}\over (q_{k}+q_{k-1}y)^{2}}$ is 
positive (negative) if $k$ is even (odd). It is immediate to check 
that any rational number $p/q\in (0,1)$, $(p,q)=1$, is the endpoint 
of exactly two branches of the iterated Gauss map. Indeed $p/q$ can 
be written as $p/q=[\bar{a}_{1},\ldots ,\bar{a}_{k}]$ with $k\ge 1$ 
and $\bar{a}_{k}\ge 2$ in a unique way and it is the left (right) 
endpoint of $I(\bar{a}_{1},\ldots ,\bar{a}_{k})$ and 
the right (left) endpoint of $I(\bar{a}_{1},\ldots ,\bar{a}_{k}-1,1)$ 
if $k$ is even (odd). 

\vskip 2. truecm\noindent
\noindent
{\titsec  A3. Distributions, Hyperfunctions, Formal Series. 
Hypoellipticity and Diophantine Conditions. }

\vskip .3 truecm
\noindent
{\bf A3.1} We follow here [H1], Chapter 9 but we also 
recommend [Ph], especially the first few chapters, for a nice 
introduction to hyperfunctions and their applications. 

Let $K$ be a non empty compact subset of ${\Bbb R}$.  
A {\it hyperfunction with support in} $K$ is a linear functional $u$ 
on the space ${\cal O}(K)$ of functions analytic in a neighborhood of $K$
such that for all neighborhood $V$ of $K$ there is a constant 
$C_{V}>0$ such that 
$$
|u(\varphi )| \le C_{V}\sup_{V}|\varphi|\; , \;\;\;
\forall \varphi\in {\cal O}(V)
\; . 
$$
We denote by $A'(K)$ the space of hyperfunctions with support in $K$. 
It is a Fr\'echet space: a seminorm is associated to each 
neighborhood $V$ of $K$. 

Let ${\cal O}^{1} (\overline{\Bbb C}\setminus K)$ denote the complex 
vector space of functions holomorphic on $(\overline{\Bbb C}\setminus K)$
and vanishing at infinity. One has the following 

\vskip .3 truecm\noindent
\Proc{Proposition A3.1}{ The spaces $A'(K)$ and 
${\cal O}^{1}(\overline{\Bbb C}\setminus K)$ are canonically isomorphic. 
To each $u\in A'(K)$ corresponds $\varphi\in 
{\cal O}^{1}(\overline{\Bbb C}\setminus K)$ given by
$$
\varphi (z) = u(c_{z})\; , \; \forall z\in {\Bbb C}\setminus K\; , 
$$
where $c_{z}(x) = {1\over \pi}{1\over x-z}$. Conversely to each 
$\varphi \in {\cal O}^{1}(\overline{\Bbb C}\setminus K)$ 
corresponds the hyperfunction 
$$
u(\psi ) = {i\over 2\pi} \int_{\gamma} \varphi (z)\psi (z) dz\; , \; 
\forall \psi \in A
$$
where $\gamma$ is any piecewise ${\cal C}^{1}$ path  winding around 
$K$ in the positive direction. We will also use the notation 
$$
u(x) = {1\over 2i}[\varphi (x+i0)-\varphi (x-i0)]
$$
for short.}

\vskip .3 truecm\noindent
\proof                                                                    
It is very easy: note that the function $x\mapsto c_{z}(x)$ is analytic in a               %
neighborhood of $K$ for all $z\notin K$. Then it is immediate to check    
applying Cauchy's formula that these two correspondences are              
surjective and are the inverse  of one another. \qed                      

\line{}

\vskip .3 truecm
\noindent
{\bf A3.2} Let ${\Bbb T}^{1}= {\Bbb R}/{\Bbb Z}\subset
\C/\Z$. A {\it hyperfunction on} $\T$ is a linear funtional $U$ on 
the space ${\cal O}({\Bbb T}^{1})$ of 
functions analytic in a complex 
neighborhood of ${\Bbb T}^{1}$ such that for all neighborhood $V$ 
of $\T$ there exists $C_{V}>0$ such that 
$$
|U(\Phi )| \le C_{V}\sup_{V}|\varphi|\; , \;\;\;
\forall \Phi\in {\cal O}(V)
\; . 
$$
We will denote  
$A'({\Bbb T}^{1})$. 
the Fr\'echet space of  hyperfunctions with support in $\T$. 
For $U\in A'(\T )$, let $\hat{U}(n) := U(e_{-n})$ with 
$e_{n}(z)=e^{2\pi i n z}$. 

\vskip .3 truecm\noindent
{\bf Exercise A3.2} Show that the doubly infinite sequence 
$(\hat{U}(n))_{n\in {\Bbb Z}}$ satisfies 
$$
|\hat{U}(n)|<C_{\varepsilon}e^{2\pi |n|\varepsilon}\; . 
$$
for all $\varepsilon >0$ and for all $n\in \Z$ with a suitably chosen
$C_{\varepsilon}>0$. Conversely show that any such sequence is the Fourier 
expansion of a unique hyperfunction with support in $\T$. 

\vskip .3 truecm
\noindent
Let ${\cal O}_{\Sigma}$ denote the complex vector space of 
holomorphic functions $\Phi\, : {\Bbb C}\setminus {\Bbb R}\rightarrow 
{\Bbb C}$, $1$--periodic, bounded at $\pm i \infty$ and such that 
$\Phi (\pm i\infty) := \lim_{\IM z\rightarrow \pm\infty }\Phi (z)$
exist and verify $\Phi (+i\infty) = -\Phi (-i\infty )$. 

\vskip .3 truecm\noindent
{\bf Exercise A3.3} Show  that the spaces $A'({\Bbb T}^{1})$ and 
${\cal O}_{\Sigma}$ are canonically isomorphic. Indeed to
each $U\in A'({\Bbb T}^{1})$ corresponds $\Phi\in 
{\cal O}_{\Sigma}$ given by
$$
\Phi (z) = U(C_{z})\; , \; \forall z\in {\Bbb C}\setminus K\; , 
$$
where $C_{z}(x) = \cotg \pi ( x-z)$. Conversely to each 
$\Phi\in {\cal O}_{\Sigma}$
corresponds the hyperfunction 
$$
U(\Psi ) = {i\over 2} \int_{\Gamma} \Phi (z)\Psi (z) dz\; , \; 
\forall \Psi \in A({\Bbb T}^{1})
$$
where $\Gamma$ is any piecewise ${\cal C}^{1}$ path  winding around a 
closed interval $I\subset {\Bbb R}$ of length $1$
in the positive direction. We will also use the notation 
$$
U(x) = {1\over 2i}[\Phi (x+i0)-\Phi (x-i0)]
$$
for short. 

\vskip .3 truecm
\noindent
The nice fact is that the following diagram commutes: 
$$
\diagram{
A'([0,1]) & 
\hfl{}{} & {\cal O}^{1}(\overline{\Bbb C}\setminus [0,1]) 
    \cr
\vfl{\sum_{Z}}{}&&\vfl{}{\sum_{Z}} 
    \cr
A'({\Bbb T}^{1}) & 
\hfl{}{} & {\cal O}_{\Sigma}
    \cr}
$$
the horizontal lines are the above mentioned isomorphisms, $\sum_{Z}$ 
is the sum over integer translates: $(\sum_{\Z}\varphi 
)(z)=\sum_{n\in\Z}\varphi (z-n)$. 

\line{}

\vskip .3 truecm\noindent
{\bf A3.3} As we have seen in A3.2  
periodic distributions and hyperfunctions are naturally identified 
with the two following subspaces of the complex vector space of {\it formal}
Fourier series
$$
\eqalign{
\varphi\in{\cal D}'(\T ) &\Leftrightarrow \varphi (\theta )=\sum_{-\infty}^{+\infty}
\hat{\varphi}(n)e^{2\pi i n\theta}\;\hbox{and there exists} M>0\, ,\; r>0\cr
&\phantom{\Leftrightarrow} \hbox{such that}\, |\hat{\varphi}(n)|\le 
M|n|^{+r}\;
\forall n\in\Z^*\; , \cr
\varphi\in{\cal A}'(\T ) &\Leftrightarrow \varphi (\theta )=\sum_{-\infty}^{+\infty}
\hat{\varphi}(n)e^{2\pi i n\theta}\;\hbox{and for all}\,\varepsilon >0 \hbox{ there exists} 
C_\varepsilon>0\cr
&\phantom{\Leftrightarrow} \hbox{such that}\, |\hat{\varphi}(n)|\le C_\varepsilon 
\exp (2\pi |n|\varepsilon )\;
\forall n\in\Z\; .\cr
}
$$
Let us now consider the following linear first--order difference equation on 
$\T^1$
$$
f_g(\theta+\alpha )-f_g(\theta )=g(\theta )
$$
where $\alpha\in\R\setminus\Q$. A necessary condition for the existence of a solution is 
that $\int_0^{2\pi}g(\theta)d\theta =0$. Thus we introduce the zero--mean Dirac delta 
function on $\T^1$
$$
\delta_{\T,0}(\theta) = \sum_{n\in\Z\, , n\not= 0} e^{2\pi i n \theta}
$$
and we note that the corresponding $f_\delta$ plays the role of a fundamental solution since
$$
f_g=f_\delta\odot g=\sum_{n=-\infty}^{+\infty}\hat{f}_\delta(n)\hat{g}(n)e^{2\pi i n \theta}
={1\over 2\pi}\int_0^{2\pi}f_\delta (
\theta-\theta_1 )g(\theta_1 )d\theta_1
$$
when the integral makes sense. 

On the other hand one clearly has 
$$
f_\delta (\theta) = \sum_{n\not= 0}{e^{2\pi i n \theta}\over 
e^{2\pi i n \alpha}-1}
$$
as a formal power series. We have the following elementary Proposition which can also 
be taken as an equivalent definition of diophantine numbers

\vskip .3 truecm\noindent
\Proc{Proposition A3.4 }{$f_\delta$ is a distribution if and only if $\alpha\in\hbox{CD}$. 
$f_\delta$ is a hyperfunction if and only if the denominators $q_n$ of the convergents 
of $\alpha$ verify $\lim_{n\rightarrow +\infty}{log q_{n+1}\over q_n}=0$.}

\vskip .3 truecm\noindent
The proof is immediate and it is left as an Exercise. 

Note that the above discussion carries over easily to the linear PDE on the two--dimensional 
torus $\T^2$
$$
(\partial_{\theta_1}+\alpha\partial_{\theta_2})f=g
$$
(which is associated to the linear flow $\dot{\theta_1}=1$, $\dot{\theta_2}=\alpha$). 
In this case one has $\delta_{\T^2,0}= \sum_{n\in\Z^2, n\not= 0}
e^{2\pi i (n_1\theta_1+n_2\theta_2)}$ and the fundamental solution is 
$f_\delta =  \sum_{n\in\Z^2, n\not= 0} {e^{2\pi i (n_1\theta_1+n_2\theta_2)}\over 
2\pi i (n_1+n_2\alpha )}$. 
Then Proposition A3.4 holds also in this case, showing that the operator 
$\partial_{\theta_1}+\alpha\partial_{\theta_2}$ is hypoelliptic if and only if 
$\alpha$ is diophantine. 

\vfill\eject
%%%%% fine appendici %%%%%
%%%%% bibliografia %%%%%
\noindent
{\tit References}

\vskip .3 truecm 
\item{[AG]} S. Alinhac, P. G\'erard ``Op\'erateurs 
pseudo--diff\'erentiels et th\`eor\'eme de Nash--Moser''
Savoirs Actuels, CNRS Editions (1991)
%\item{[Ah1]} L.V. Ahlfors ``Complex Analysis'' 3rd Ed. McGraw--Hill 
%(1979)
\item{[Ah1]} L.V. Ahlfors ``Conformal Invariants: Topics in Geometric 
Function Theory'' McGraw--Hill (1973)
\item{[Ah2]} L.V. Ahlfors ``Lectures on Quasiconformal Mappings''
Van Nostrand (1966)
\item{[AM]} R. Abraham, J. Marsden ``Foundations of Mechanics''
Benjamin Cummings, New York (1978)
\item{[Ar1]} V.~I. Arnol'd ``Small denominators and problems of 
stability of motion in classical celestial mechanics'' 
Russ. Math. Surv. {\bf 18} (1963), 85--193. 
\item{[Ar2]}  V.~I. Arnol'd ``Instability of dynamical systems with 
several degrees of freedom'' Sov. Math. Dokl. {\bf 5} (1964), 581--585.
%\item{[Ar]} V. I. Arnold ``On the mappings of the circumference onto 
%itself'' Translations of the Amer. Math. Soc. {\bf 46}, 2nd series, 
%213--284 (1961)
\item{[Ar3]} V.~I. Arnol'd ``Geometrical Methods in the Theory of 
Ordinary Differential Equations'' Springer--Verlag (1983)
\item{[AKN]}  V.~I. Arnol'd, V.~V. Kozlov and A.~I. Neishtadt
``Dynamical Systems III'', 
Springer--Verlag (1988). 
%\item{[BM]} Berretti A and Marmi S 1994
%Scaling near Resonances and Complex Rotation Numbers for the
%Standard Map {\it Nonlinearity}  {\bf 7}  603--621
%\item{[Be]} B. Benzaghou ``Alg\`ebres de Hadamard''
%Bull. Soc. Math. France {\bf 98} (1970), 209--252
\item{[Be]} A. Beardon ``Iteration of Rational Functions''
Springer--Verlag (1991)
\item{[BFGG]} G. Benettin, G. Ferrari, L. Galgani and A. Giorgilli
``An Extension of the Poincar\'e--Fermi Theorem on the Nonexistence 
of Invariant Manifolds in Nearly Integrable Hamiltonian Systems''
Il Nuovo Cimento {\bf 72B} (1982) 137
\item{[BHS]} H.~W. Broer, G.~B. Huitema and M.~B. Sevryuk
``Quasi--Periodic Motions in Families of Dynamical Systems''
Springer--Verlag (1996)
%\item{[Bo1]} E. Borel ``Le\c cons sur les fonctions monog\`enes
%uniformes d'une variable complexe'' Gauthier--Villars, Paris
%(1917).
%\item{[BG]} A. Berretti and G. Gentile 
\item{[Bo]} J.~B. Bost ``Tores invariants des syst\`emes dynamiques 
hamiltoniens'' S\'e\-mi\-nai\-re Bourbaki {\bf 639}, Ast\'erisque {\bf 
133--134} (1986), 113--157.
%\item{[BPV]} N. Buric, I. Percival and F. Vivaldi
%``Critical Function and Modular Smoothing'' {\it Nonlinearity}
%{\bf 3} (1990), 21--37.
\item{[Br]} A. D. Brjuno ``Analytical form of differential equations''
{\it Trans. Moscow Math. Soc.} {\bf 25} (1971), 131-288; {\bf 26}
(1972), 199-239.
\item{[CG]} L. Carleson and T. Gamelin ``Complex Dynamics'' 
Universitext, Sprin\-ger--Verlag, Berlin Heidelber New York (1993)
\item{[CM]} T. Carletti and S. Marmi ``Linearization of Analytic and 
Non--Analytic Germs of Diffeomorphisms of $(\C ,0)$''
Bulletin de la Societ\'e Math\'ematique de France (1999)
\item{[Da]} A.M. Davie ``The critical function for the 
semistandard map'' {\it Nonlinearity} {\bf 7} (1994), 219 - 229.
%\item{[Da2]} A. M. Davie ``Renormalisation for analytic 
%area--preserving maps'' University of Edinburgh preprint, (1995).
\item{[DeB]} L. de Branges ``A proof of the Bieberbach conjecture''
{\it Acta Math.} {\bf 154} (1985), 137--152.
\item{[Di]} J. Dieudonn\'e ``Calcul Infinit\'esimal''
Hermann, Paris (1980)
\item{[Do]} A. Douady ``Disques de Siegel et anneaux de Herman''
S\'eminaire Bourbaki n. 677, {\it Ast\'erisque} {\bf 152--153}
(1987), 151--172
%\item{[E1]} Ecalle J 1981 and 1985 
%{\it Les fonctions r\'esurgentes et leurs applications. Tomes I, II, III},
%(Publications math\'ematiques d'Orsay  {\bf 81--05},  {\bf 81--06}, 
% {\bf 85--05})
\item{[Fa]} K. Falconer ``Fractal Geometry. Mathematical Foundations
and Applications'' John Wiley and Sons (1990)
%\item{[Fal]} G. Faltings ``Diophantine Approximation on Abelian 
%Varieties'' {\it Ann. of Math.} {\bf 133} (1991) 549--576
%\item{[Gr]} J. B. Garnett ``Bounded Analytic Functions'' Academic Press,
%New York, (1981).
\item{[Ga]} G. Gallavotti ``Quasi--Integrable Mechanical Systems''
in Ph\'enom\`enes critiques, sy\-st\`e\-mes al\'eatoires, th\'eories de 
jauge, Part I, II, Les Houches 1984, North--Holland, Amsterdam (1986)
539--624
\item{[GM]} A. Giorgilli and  A. Morbidelli ``Invariant KAM tori and 
global stability for Hamiltonian systems'' ZAMP {\bf 48} (1997), 
102--134.
\item{[Gr]} M.L. Gromov ``Smoothing and Inversion of Differential 
Operators'' {\it Math. USSR Sbornik} {\bf 17} (1972), 381--434
\item{[Ha]} R.S. Hamilton ``The Inverse Function Theorem of
Nash and Moser'' {\it Bull. A.M.S.} {\bf 7} (1982), 65--222.  
\item{[H1]} L. H\"ormander ``The Analysis of Linear Partial 
Differential Operators I'' Grundlehren der mathematischen 
Wissenschaften {\bf 256}, Springer--Verlag, Ber\-lin, Heidelberg, 
New York, Tokyo (1983)
\item{[H2]}  L. H\"ormander ``The boundary problem of physical
geodesy'' {\it Arch. Rat. Mech. Anal.} {\bf 62} (1976), 
1--52
\item{[He1]} M.~R. Herman ``Examples de fractions rationelles ayant 
une orbite dense sur la sphere de Riemann'' Bulletin de la 
Societ\'e Math\'ematique de France {\bf 112} (1984), 93--142
\item{[He2]} M.~R. Herman ``Simple proofs of local 
conjugacy theorems for diffeomorphisms of the circle
with almost every rotation numbers'' Bull. Soc. Bras. Mat.
{\bf 16} (1985) 45--83
%\item{[He1]} M. R. HERMAN - {\it Sur les courbes invariantes par les 
%diff\'eomorphismes de l'anneau}, Ast\'erisque {\bf 103--104} 
%(1983).
%\item{[He2]} M.~R. HERMAN - {\it Some open problems in dynamical systems},
%ICM1998.
\item{[He3]} M.~R. Herman ``Are there critical points one the 
boundary of singular domains?'' Commun. Math. Phys. {\bf 99} (1985)
593--612.
\item{[He4]} M.~R. Herman ``Recent results and some open questions 
on Siegel's linearization theorem 
of germs of complex analytic diffeomorphisms of ${\bf C}^n$ near a fixed 
point''  Proc. VIII Int. Conf. Math. Phys. Mebkhout and Seneor eds. 
(Singapore: World Scientific) (1986), 138--184.
\item{[He5]} M.~R. Herman ``D\'emonstration du th\'eor\`eme
des courbes translat\'ees par dif\-f\'e\-o\-mor\-phi\-smes de l'anneau; 
d\'emonstration du th\'eor\`eme des tores invariants'' manuscripts
(1980) and ``Abstract methods in small divisors: implicit function
theorems in Fr\'echet spaces'', lectures given at the CIME conference on 
Dynamical Systems and Small Divisors, Cetraro 1998
%\item{[L1989]} J. LASKAR - {\it A numerical experiment on the chaotic
%behaviour of the Solar System}, Nature {\bf 338} (1989), 237--238.
\item{[HL]} G.~H. Hardy and J.~E. Littlewood ``Notes on the theory of 
series (XXIV): a curious power series'' Proc. Cambridge Phil. Soc. 
{\bf 42} (1946), 85--90
\item{[HW]} G.H. Hardy and E.M. Wright ``An introduction to the theory of 
numbers'' Fifth Edition, Oxford Science Publications (1990)
\item{[K]} A.~N. Kolmogorov ``On the persistence of conditionally 
periodic motions under a small perturbation of the Hamilton function''
Dokl. Akad. Nauk SSSR {\bf 98} (1954) 527--530 (in Russian: English
translation in G. Casati and J. Ford, editors, Stochastic Behavior in
Classical and Quantum Hamiltonian Systems, Lecture Notes in Physics 
{\bf 93} (1979) 51--56 Springer--Verlag).
%\item{[KH]} A. Katok and B. Hasselblatt ``Introduction to the modern 
%theory of dynamical systems'' Encyclopedia of Mathematics and its 
%Applications {\bf 54}, Cambridge University Press, (1995)
\item{[L]} S. Lang ``Introduction to Diophantine Approximation''
Addison--Wesley (1966)
\item{[Lo]} P. Lochak ``Canonical perturbation theory via 
simultaneous approximations'', Russ. Math. Surv. {\bf 47} (1992), 
57--133. 
%\item{[Ma1982]} J.N. Mather ``Existence of quasi--periodic orbits for 
%twist homeomorphisms of the annulus'' Topology {\bf 21} 457--467 
%(1982)
\item{[Ma1]} S. Marmi ``Critical Functions for Complex Analytic Maps''
{\it J. Phys. A: Math. Gen.} {\bf 23} (1990), 3447--3474
\item{[Ma2]} S. Marmi ``Chaotic Behaviour in the Solar System 
(Following J. Laskar)'' 
S\'eminaire Bourbaki n. 854, November 1998, 
to appear in Ast\'erisque
\item{[Me]} Y. Meyer ``Algebraic Numbers and Harmonic Analysis''
North--Holland Mathematical Library {\bf 2} (1972)
\item{[MM]} L. Markus and K.~R. Meyer ``Generic Hamiltonian Systems 
are neither integrable nor ergodic'' Memoirs of the A.M.S. {\bf 144}
(1974)
\item{[MMY]} S. Marmi, P. Moussa and J.--C. Yoccoz
``The Brjuno functions and their regularity properties''
{\it Commun. Math. Phys.} {\bf 186} (1997), 265-293
\item{[MMY2]} S. Marmi, P. Moussa and J.--C. Yoccoz
``Complex Brjuno Functions'' preprint SPhT Saclay, France, 
71 pages (1999)
\item{[Mc]} C. T. McMullen ``Complex Dynamics and Renormalization''
Ann. of Math. Studies, Princeton University Press (1994)
\item{[Mn]} R. Ma\~ne ``Ergodic Theory and Differentiable Dynamics''
Springer--Verlag (1987)
\item{[M]} J. Moser ``A rapidly convergent iteration method and 
nonlinear differential equations'' {\it Ann. Scuola Norm. Sup. Pisa}
{\bf 20} (1966) 499-535
\item{[N]} J. Nash ``The embedding problem for Riemannian manifolds''
{\it Ann. of Math.} {\bf 63} (1956) 20--63
\item{[Ne]} N.N. Nekhoroshev ``An exponential estimate for the time 
of stability of nearly integrable Hamiltonian systems'' {\it Russ. Math. 
Surveys} {\bf 32} (1977), 1--65.
\item{[Ni]} L. Niremberg ``An abstract form of the non--linear 
Cauchy--Kowalewskaya theorem'' {\it J. Diff. Geom.} {\bf 6} (1972) 
561--576
\item{[PM1]} R. P\'erez--Marco ``Solution compl\`ete au probl\`eme 
de Siegel de lin\'earisation d'une application holomorphe au 
voisinage 
d'un point fixe (d'apr\`es J.--C. Yoccoz)'' 
S\'eminaire Bourbaki n. 753, {\it Ast\'erisque} {\bf 206}
(1992), 273--310
%\item{[PM2]} R. P\'erez--Marco ``Non linearizable holomorphic dynamics 
%having an uncountable number of symmetries'' {\it Invent. Math.}
%{\bf 119} (1995), 67--127
\item{[Ph]} F. Pham (Editor) ``Hyperfunctions and Theoretical Physics''
Lecture Notes in Mathematics {\bf 449} Springer--Verlag (1975)
\item{[P]} H. Poincar\'e  ``Les M\'ethodes Nouvelles de la M\'ecanique 
Celeste'', tomes I--III Paris Gauthier--Villars (1892, 1893, 1899).
\item{[Po]} Ch. Pommerenke ``Boundary Behaviour of Conformal Maps''
Grundlehren der Mathematischent Wissenschaften {\bf 299}, 
Springer--Verlag (1992)
\item{[P\"o]} J. P\"oschel ``Integrability of Hamiltonian systems on 
Cantor sets'' Comm. Pure Appl. Math. {\bf 35} 653--696 (1982)
%\item{[P\"o2]} J. P\"oschel ``On Invariant Maniflods of Complex Analytic 
%Mappings Near Fixed Points'' {\it Exp. Math.} {\bf 4} (1986), 97--109
\item{[Re]} R. Remmert ``Classical Topics in Complex Function 
Theory'' Graduate Texts in Mathematics {\bf 172}, Springer--Verlag 
(1998)
%\item{[Ris]} E. Risler ``Lin\'earisation des perturbations 
%holomorphes des rotations et applications'' M\'emoires de la  
%Soc. Math. France {\bf 77} (1999)
\item{[R\"u]} H. R\"ussmann ``Kleine Nenner II: Bemerkungen zur 
Newtonschen Methode'' Nachr. Akad. Wiss. G\"ottingen Math. Phys. Kl
(1972) 1--20
\item{[S]} C.L. Siegel ``Iteration of analytic functions''
Annals of Mathematics $\bf{43}$ (1942), 807-812.
\item{[Sch1]} W.M. Schmidt ``Diophantine Approximation''
Lecture Notes in Mathematics, {\bf 785}, Springer--Verlag (1980)
\item{[Sch2]} W.M. Schmidt ``Diophantine Approximations and Diophantine
Equations''
Lecture Notes in Mathematics, {\bf 1467}, Springer--Verlag (1991)
\item{[Ser]} F. Sergeraert ``Un th\'eor\`eme de fonctions 
implicites sur certains espaces de Fr\'echet et quelques applications''
{\it Ann. Scient. \`Ec. Norm. Sup.} {\bf 5} (1972), 599--660.  
%\item{[Sim]} B. Simon ``Almost periodic Schr\"odinger operators IV.
%The Maryland model'' Ann. Phys. {\bf 159} (1985) 157--183
%\item{[Th]} V. Thilliez 
%``Quelques propri\'et\'es de quasi--analyticit\'e''
%Gazette des Mathematiciens, {\bf 70} (1996), 49--68
\item{[ST]} J. Silverman and J. Tate ``Rational Points on 
Elliptic Curves'' Undergraduate Texts in Mathematics, 
Springer--Verlag (1992)
%\item{[St]} E.M. Stein ``Singular integrals and differentiability properties 
%of functions'' Princeton University Press (1970)
\item{[SZ]} D. Salomon and E. Zehnder ``KAM theory in confuguration 
space'' {\it Comm. Math. Helvetici} {\bf 64}, (1989), 84--132
%\item{[SM]} C.~L. Siegel and  J. Moser ``Lectures on Celestial 
%Mechanic'' Springer--Verlag, Berlin--Heidelberg--New York (1971).
\item{[St]} S. Sternberg ``Celestial Mechanics'' (two volumes) 
W.A. Benjamin, New York (1969).
\item{[Va]}  F.H. Vasilescu ``Analytic Functional Calculus'' 
D. Reidel Publ. Co. (1982)
%\item{[Y2]} J.--C. YOCCOZ - {\it Introduction to Hyperbolic Dynamics} 
%in ``Real and Complex Dynamical Systems'' B. Branner and P. Hjorth eds.
%NATO ASI Series {\bf C464}, Kluwer, Dordrecht (1995), 265--291
%\item{[Yo]} L.S. Young ``Ergodic Theory of Differentiable Dynamical 
%Systems'' in ``Real and Complex Dynamical Systems'' 
%B. Branner and P. Hjorth eds.
%NATO ASI Series {\bf C464}, Kluwer, Dordrecht (1995), 293--336. 
%\item{[Ri]} G. J. Rieger ``Mischung und Ergodizit\"ata bei Kettenbruchen 
%nach n\"achsten genzen'' {\it J. Reine Angew. Math.} {\bf 310} (1979),
%171-181.
%\item{[Vo]} S.M. Voronin ``Analytic Classification of Germs of Conformal 
%Mappings $(\C ,0)\rightarrow (\C ,0)$ with Identity Linear Part''
%{\it Funct. Analys. Appl.} {\bf  ??} (1981), ???
%\item{[We]} A. Weil ``Elliptic Functions according to Eisenstein and 
%Kronecker'' Springer-Verlag, Berlin Heidelberg New York (1976)
%\item{[Wh]} H. Whitney ``Analytic extensions of differentiable functions
%defined in closed sets'', Trans. Amer. Math. Soc. {\bf 36} (1934)
%63--89
\item{[Yo1]} J.--C. Yoccoz ``An introduction to small divisors problems''
in ``From number theory to physics'', 
M. Waldschmidt, P. Moussa, J.M. Luck and C. Itzykson
(editors) Springer--Verlag (1992), 659--679
\item{[Yo2]} J.--C. Yoccoz ``Th\'eor\`eme de Siegel, nombres de Bruno et 
polyn\^omes quadratiques'' {\it Ast\'erisque} {\bf 231} (1995), 3-88. 
\item{[Yo3]} J.--C. Yoccoz, lectures given at the CIME conference on 
Dynamical Systems and Small Divisors, Cetraro 1998, to appear in 
Lecture Notes in Mathematics
\item{[Ze1]} E. Zehnder ``Generalized Implicit Function Theorems
with Applications to some Small Divisor Problems (I and II)''
{\it Commun. Pure Appl. Math.} {\bf 28} (1975) 91--140, 
{\bf 29} (1976) 49--113.
\item{[Ze2]} E. Zehnder ``A simple proof of a generalization of a 
Theorem by C.~L. Siegel'' Lecture Notes in Mathematics {\bf 597}
(1977) 855--866.
%%%%% fine bibliografia %%%%%
\vfill\eject
%%%%% indice analitico %%%%%
\noindent
{\tit Analytical index} 

\vskip .3 truecm\noindent

%{\sevenrm
\item{} action--angle variables pp.49--53 
\item{} adjoint action p. 4
\item{} area formula p. 70
\item{} area theorem p. 70
\item{} Arnol'd--Liouville Theorem p. 50

\vskip .3 truecm\noindent
\item{} Beltrami equation p. 74
\item{} Best Approximation Theorem p. 22
\item{} Bieberbach--De Branges Theorem pp. 31, 70, 71
\item{} Brjuno number pp. 25, 27, 32, 33, 34
\item{} Brjuno function p. 25, 28, 34
\item{} Brjuno Theorem p. 17, 27

\vskip .3 truecm\noindent
\item{} Caratheodory's theorem p. 68
\item{} centralizer p. 4, 8
\item{} completely canonically integrable pp. 49, 50, 51, 53, 54
\item{} conformal map p. 14, 17, 67, 73, 74
\item{} conformal capacity pp. 14, 36, 68
\item{} conjugate p. 4
\item{} continued fractions pp. 21, 22, 23, 25, 39, 75  
\item{} Cremer's Theorem p. 9
\item{} critical point (on the boundary) pp. 20, 39
\item{} cycle p. 12

\vskip .3 truecm\noindent
\item{} Darboux's theorem  p. 47
\item{} Davie's lemmas pp. 29, 30, 32 
\item{} dilatation pp. 73, 74 
\item{} diophantine number, condition, vector  pp. 20, 23, 25, 39, 43, 
55--57, 64, 65, 81
\item{} Douady--Ghys' Theorem pp. 21, 27, 36
\item{} Douady--Hubbard's theorem pp. 16, 74

\vskip .3 truecm\noindent
\item{} Fatou set pp. 12, 13
\item{} Fatou's theorem pp. 20, 69
\item{} Fr\'echet space pp. 58--62, 64

\vskip .3 truecm\noindent
\item{} Gauss map pp. 75, 77
\item{} $G_{\delta}$ set p. 9
\item{} geometric renormalization pp. 34, 40
\item{} germ p. 3
\item{} Gevrey class p. 32
\item{} Grunsky norm p. 39

\vskip .3 truecm\noindent
\item{} Hamiltonian pp. 48, 53, 57
\item{} Hardy--Sobolev spaces pp. 40, 44
\item{} hedgehog p. 14
\item{} hyperfunction pp. 78--80
\item{} hypoelliptic pp. 55, 81
\item{} Hurwitz' Lemma p. 72

\vskip .3 truecm\noindent
\item{} Jarnik's theorem p. 24
\item{} Julia set p. 12

\vskip .3 truecm\noindent
\item{} KAM Theorem pp. 56, 57, 64
\item{} Koebe $1/4$--Theorem p. 71
\item{} Koebe distorsion theorems p. 71 
\item{} Koebe function p. 71
\item{} Koebe transform p. 71
\item{} Koenigs--Poincar\'e Theorem pp. 7-9

\vskip .3 truecm\noindent
\item{} Lagrange's inversion theorem p. 44
\item{} linearizable, linearization pp. 5--10, 13--14, 16, 19, 22, 
27--28, 32--34, 38, 40, 46
\item{} Liouville number p. 24
\item{} Liouville's Theorem p. 23
\item{} loss of differentiability pp. 40--41, 46, 56, 58, 61

\vskip .3 truecm\noindent
\item{} Maximum principle p. 67
\item{} Measurable Riemann Mapping Theorem p. 74

\vskip .3 truecm\noindent
\item{} Nash--Moser Theorem pp. 58, 62--64
\item{} Nekhoroshev Theorem pp. 56--57
\item{} normal family p. 11--12

\vskip .3 truecm\noindent
\item{} orbit p. 4

\vskip .3 truecm\noindent
\item{} Poincar\'e Theorem p. 54
\item{} Poisson bracket p. 48

\vskip .3 truecm\noindent
\item{} quadratic polynomial pp. 16, 37
\item{} quadrilateral pp. 73--74
\item{} quasicircle p. 39
\item{} quasiconformal map pp. 16, 39, 73--74
\item{} quasi--integrable system pp. 53--54, 56
\item{} quasiperiodic function p. 49

\vskip .3 truecm\noindent
\item{} rational map pp. 11--12 
\item{} region p. 67
\item{} Riemann mapping theorem p. 68
\item{} Roth's Theorem pp. 23--55

\vskip .3 truecm\noindent
\item{} Schwarzian derivative p. 42
\item{} Schwarz's Lemma p. 67
\item{} Siegel--Brjuno Theorem pp. 17, 27  
\item{} Siegel disk p. 14
\item{} smoothing operators pp. 58, 61--62
\item{} spherical metric p. 11
\item{} stable point pp. 13--14
\item{} symmetry p. 4
\item{} symplectic manifold p. 47

\vskip .3 truecm\noindent
\item{} tame map pp. 61--62, 64--66
\item{} tame Fr\'echet space pp. 61--62

\vskip .3 truecm\noindent
\item{} ultradifferentiable power series p. 32
\item{} Uniformization theorem p. 69
\item{} uniquely ergodic pp. 38, 49
\item{} univalent function pp. 10,39,  67, 69--72

\vskip .3 truecm\noindent
\item{} Yoccoz's lower bound p. 28
\item{} Yoccoz's proof of Siegel's Theorem pp. 16--20
\item{} Yoccoz's  $u(\lambda )$ function pp. 17--20, 37--39
\item{} Yoccoz's Theorem pp. 33--34
%}
%%%%% fine indice analitico %%%%%
\vfill\eject
%%%%% lista dei simboli %%%%%
\noindent
{\titsec List of symbols}

\vskip .3 truecm
\noindent
\item{} $\hbox{Ad}$: adjoint action 
\item{} $c(\Omega,z_0)$: conformal capacity of $\Omega$ w.r.t. $z_0$ 
\item{} $\C$: complex plane
\item{} $\C^*$: $\C\setminus\{ 0\}$
\item{} $\overline{\C}$: Riemann sphere 
\item{} $\C\{z\}$: ring of convergent power series in one complex variable
\item{} $\C [[z]]$: ring of formal power series in one complex variable
\item{} $\hbox{Cent}$: centralizer 
\item{} $\widehat{\hbox{Cent}}$: formal centralizer
\item{} $D_f(z_0)$, $D_f$: dilatiation of f (maximal)  
\item{} $\D$: open unit disk $\{|z|<1\}$
\item{} $\D_r$: open disk $\{|z|<r\}$ of radius $r>0$. 
\item{} $\E$: outer disk $\{|z|>1\}$
\item{} $[f]$: orbit of $f$ 
\item{} $F(R)$: Fatou set of $R$ 
\item{} $G$: group of germs of holomorphic diffeomorphisms of $(\C, 0)$
\item{} $\hat{G}$: formal analogue of $G$
\item{} $G_\lambda$: elements of $G$ with linear part $\lambda$ 
\item{} $\hat{G}_\lambda$: formal analogue of  $G_\lambda$
\item{} $J(R)$: Julia set of $R$
\item{} $\N$: non--negative integers
\item{} $\Omega$: a region of $\C$
\item{} $\Q$: rational integers
\item{} $R_\lambda$: the germ $R_\lambda (z)=\lambda z$
\item{} $S$: univalent maps on $\D$
\item{} $S_\lambda$: elements of $S$ with linear part $R_\lambda$
\item{} $S_{\S^1}$: elements of $S$ with linear part of unit modulus 
\item{} $\S^1$: unit circle $\{|z|=1\}$
\item{} $u$ Yoccoz's function, see Section 3.1
\item{} ${\cal Y}$: see Chapter 4
\item{} $\Z$: integers
%%%%% fine lista dei simboli %%%%%

\noindent
Stefano Marmi, Dipartimento di Matematica e Informatica, 
Universit\`a di Udine, Via delle Scienze 206, Localit\`a Rizzi, 33100 
Udine, Italy;  marmi@dimi.uniud.it
\vfill\eject
\bye